\documentclass[reqno]{amsart}

\usepackage{caption}
\newenvironment{code}{\captionsetup{type=listing}{}}

\usepackage[draft]{minted}
\usepackage{amssymb}
\usepackage{amsmath}
\usepackage{amsthm}
\usepackage{mathtools}
\usepackage[svgnames]{xcolor}
\usepackage[pdftex]{hyperref}
\usepackage{enumitem}
\usepackage{subcaption}
\usepackage{placeins}
\usepackage{vmargin}
\setmarginsrb {2cm}{2cm}{2cm}{2cm} {1cm}{1cm}{0.5cm}{0.5cm}
\usepackage{tcolorbox}
\tcbset{colback=White}
\usepackage{tikz,tkz-tab}
\usetikzlibrary{matrix}
\usetikzlibrary{trees}
\usepackage{pgfplots}
\pgfplotsset{compat=newest} 
\pgfplotsset{plot coordinates/math parser=false}
\usepgfplotslibrary{polar}
\usepackage{float}
\usepackage{algorithm}
\usepackage{algorithmic}

\graphicspath{{graphicspath/}{Figures/}}

\hypersetup{
  pdftitle={},
  pdfauthor={Stefan Le Coz <stefan.lecoz@math.cnrs.fr>},
  pdfsubject={},
  pdfkeywords={}, 
  colorlinks=true,
  linkcolor=DarkBlue  ,          
  citecolor=DarkRed,        
  filecolor=DarkMagenta,      
  urlcolor=DarkGreen,           
}

\newtheorem{theorem}{Theorem}[section]

\theoremstyle{definition}
\newtheorem{assumption}[theorem]{Assumption}

\theoremstyle{remark}

\DeclarePairedDelimiter{\norm}{\lVert}{\rVert}

\newcommand{\scalar}[2]{\left( #1,#2 \right)}

\newcommand{\dual}[2]{\left\langle #1,#2 \right\rangle}

\newcommand{\R}{\mathbb{R}}

\DeclareMathOperator{\sech}{sech}

\DeclareMathOperator{\dn}{dn}

\renewcommand{\leq}{\leqslant}
\renewcommand{\geq}{\geqslant}

\DeclareMathAlphabet{\mathpzc}{OT1}{pzc}{m}{it}

\usepackage{color}

\input{default.pygstyle.tex}
\input{default-pyg-prefix.pygstyle.tex}

\begin{document}

\title[Numerical Simulations on Nonlinear Quantum Graphs with the {GraFiDi} Library]{Numerical Simulations on Nonlinear Quantum Graphs\\ with the {GraFiDi} Library}

\author[C.~Besse]{Christophe Besse}
\author[R.~Duboscq]{Romain Duboscq}
\author[S.~Le Coz]{Stefan Le Coz}
\thanks{This work is
  partially supported by the project \emph{PDEs on Quantum Graphs} from CIMI Labex ANR-11-LABX-0040}

\address[Christophe Besse and Stefan Le Coz]{Institut de Math\'ematiques de Toulouse ; UMR5219,
  \newline\indent
  Universit\'e de Toulouse ; CNRS,
  \newline\indent
  UPS IMT, F-31062 Toulouse Cedex 9, 
  \newline\indent
  France}
\email[Christophe Besse]{Christophe.Besse@math.univ-toulouse.fr}
\email[Stefan Le Coz]{stefan.lecoz@math.cnrs.fr}

\address[Romain Duboscq]{Institut de Math\'ematiques de Toulouse ; UMR5219,
  \newline\indent
  Universit\'e de Toulouse ; CNRS,
  \newline\indent
  INSA IMT, F-31077 Toulouse, 
  \newline\indent
  France}
\email[Romain Duboscq]{Romain.Duboscq@math.univ-toulouse.fr}


\subjclass[2010]{}

\date{\today}

\begin{abstract}
Nonlinear quantum graphs are metric graphs equipped with a nonlinear Schr\"odinger equation. Whereas in the last ten years they have known considerable developments on the theoretical side, their study from the numerical point of view remains  in its early stages. The goal of this paper is to present the Grafidi library~\cite{Grafidi}, a Python library which has been developed with the numerical simulation of nonlinear Schr\"odinger equations on graphs in mind. We will show how, with the help of the Grafidi library, one can implement the popular normalized gradient flow and nonlinear conjugate gradient flow methods to compute ground states of a nonlinear quantum graph. We will also simulate the dynamics of the nonlinear Schr\"odinger equation with a Crank-Nicolson relaxation scheme and a Strang splitting scheme. Finally, in a series of numerical experiments on various types of graphs, we will compare the outcome of our numerical calculations for ground states with the  existing  theoretical results, thereby illustrating the versatility and efficiency of our implementations in the framework of the Grafidi library. 
\end{abstract}

\maketitle

\section{Introduction}

The nonlinear Schr\"odinger equation
\begin{equation*}
iu_t+\Delta_{\Omega} u+f(u)=0, 
\end{equation*}
where $ u:\R_t\times \Omega_x\to \mathbb C$ is a popular model for wave propagation in Physics. It appears in particular in the modeling of Bose-Einstein condensation and in nonlinear optics. In general, the set $\Omega$ is chosen to be either the full space $\R^d$ (with $d=1$ in general in optics and $d=1,2$ or $3$ for Bose-Einstein condensation), or a subdomain of the full space. For example, in Bose-Einstein condensation, the potential might be chosen in such a way that the condensate is confined in various shapes $\Omega$, e.g. balls or cylinders. In some cases, the shape of $\Omega$ is very thin in one direction, for example in the case of $Y$-junctions (see e.g.~\cite{ToOsDe08}), or in the case of $H$-junctions (see e.g.~\cite{HuTrMa11}). In these cases, it is natural to perform a reduction to a one-dimensional model set on a graph approximating the underlying spatial structure (see e.g.~\cite{SoBaMa17}).

The study of nonlinear quantum graphs, i.e. metric graphs equipped with a nonlinear evolution equation of Schr\"odinger type, is therefore motivated at first by applications in Physics. An overview of  various applications of nonlinear Schr\"odinger equations on metric graphs in physical settings is proposed by Noja in~\cite{No14}.  One may also refer to~\cite{GnWa16,SaSoBaMa13} for  analysis of standing waves in a physical context. The validity of the graph approximation for planar branched systems was considered  by   Sobirov, Babadjanov and Matrasulov~\cite{SoBaMa17}.

The mathematical aspects of nonlinear equations set on metric graphs are also interesting on their own. Among the early studies, one finds the works of Ali Mehmeti~\cite{Al94}, see also~\cite{AlVoNi01}.  Dispersive effects for the Schr\"odinger group have been considered on star graphs~\cite{MeAmNi15} and the tadpole graph~\cite{MeAmNi17}.
In the last ten years, a particular theoretical aspect has attracted considerable interest: the ground states of nonlinear quantum graphs, i.e. the minimizers on  graphs of the Schr\"odinger energy at fixed mass constraint. The literature devoted to ground states on graphs is already too vast to give an exhaustive presentation of the many works on the topic, we refer to Section~\ref{sec:ground-states} for a small sample of relevant examples of the existing results. 
There seem to be relatively few works devoted to the numerical simulations of nonlinear quantum graphs (one may refer e.g. to~\cite{BeMaPe19,KaMaPeXi20,MaPe16} which are mostly theoretical works completed with a numerical section).

In view of the sparsity of numerical tools adapted to quantum graphs, we have developed a Python library, the Grafidi library\footnote{See \url{https://plmlab.math.cnrs.fr/cbesse/grafidi}}, which aims at rendering the numerical simulation of nonlinear quantum graphs simple and efficient.

From a conceptual point of view, the library relies on the finite difference approximation of the Laplacian on metric graphs with vertex conditions described by matrices (see Section~\ref{sec:preliminaries} for details). Inside each of the edges of the graph, one simply uses the classical second order finite differences approximation for the second derivative in one dimension. On the other hand,  for discretization points close to the vertices, the finite differences approximation would involve the value of the function at the vertex, which is not directly available. To substitute for this value, we make use of (again) finite differences approximations of the boundary conditions. As a consequence, the approximation of the Laplacian of a function close to a vertex involves values of the function on each of the edges incident to this vertex. Details are given in Section~\ref{sec:space_disc}.

The basic functionalities of the Grafidi library are presented in Section~\ref{sec:library}. The Grafidi library has been conceived with ease of use in mind and the user should not need to deal with technicalities for most of common uses. A graph is given as a list of edges, each edge being described by the labels (e.g. $A$, $B$, etc.) of the vertices that the edge is connecting and the length of the edge. With this information, the graph-constructor of the library constructs the graph and the matrix of the Laplacian on the graph with Kirchhoff (i.e. default) conditions at the vertices and a default number of discretization points. One may obviously choose to assign other types of vertices conditions, either with one the pre-implemented type ($\delta$, $\delta'$, Dirichlet) or even with a user defined vertex condition for advanced uses.  A function on the graph is then given by the collection of functions on each of the edges. The graph and functions on the graph are easily represented with commands build in the Grafidi library. 

We present in Section~\ref{sec:numerical-methods} the implementation for nonlinear quantum graphs of four numerical methods popular in the simulation of nonlinear Schr\"odinger equations.

The first two methods that we present concern the computation of ground states, i.e. minimizers of the energy at fixed mass. Ground states are ubiquitous in the analysis of nonlinear Schr\"odinger equations: they are the profiles of orbitally stable standing wave solutions and serve as building blocks for the analysis of the dynamics, in particular in the framework of the \emph{Soliton Resolution Conjecture}. The two methods that we implement are the  normalized gradient flow, which was analyzed in details in our previous work~\cite{BeDuLe20}, and the conjugate gradient flow, which was described in~\cite{antoine2017efficient,danaila2017computation} in a general domain. 
The idea behind these two methods is that, since the ground states are minimizers of the energy at fixed mass, they may be obtained at the continuous level by using the so-called continuous normalized gradient flow, i.e. a gradient flow corresponding to the Schr\"odinger energy, projected on the sphere of constant mass.

The next two methods that we present in Section~\ref{sec:numerical-methods} concern the simulation of the nonlinear Schr\"odinger flow on the graph. Numerical schemes for nonlinear Schr\"odinger equations abound, we have selected a Crank-Nicolson relaxation scheme and a Strang splitting scheme, which have both been shown to be very efficient for the simulation of the Schr\"odinger flow (see~\cite{besse2004relaxation, weideman1986split}). As for the methods to compute ground states, thanks to the Grafidi library, the implementation of the time-evolution methods is not more difficult on graphs than it is in the case of a full domain.

To illustrate and validate further the use of the Grafidi library and the numerical methods presented, we have performed a series of numerical experiments in various settings in Section~\ref{sec:ground-states}. As the theoretical literature is mainly devoted to the analysis of ground states, we have chosen to also focus on the calculations of ground states using the normalized and conjugate gradient flows. We distinguish between four categories of graphs: compact graphs, graphs with a finite number of edges and at least one semi-infinite edge, periodic graphs and trees. For each of these types of graphs, we perform ground states calculations. The comparison of the  outcomes of our experiments with the existing theoretical results reveals an excellent agreement between the two.

\section{Space discretization of the Laplacian on graphs}
\label{sec:space_disc}

\subsection{Preliminaries}
\label{sec:preliminaries}
A \emph{metric graph} $\mathcal G$ is a collection of edges $\mathcal E$ and vertices $\mathcal V$. Two vertices can be connected by more than one edge (in which case we speak of \emph{bridge}), and an edge can connect a vertex to himself (in which case we refer to the edge as \emph{loop}). To each edge $e\in\mathcal E$, we associate a length $l_e$ and identify the edge $e$ with the interval $[0,l_e]$ ($[0,\infty)$ if $l_e=\infty$). 

A \emph{function} $\psi$ on the graph is a collection of maps $\psi_e:I_e\to\R$ for each $e\in\mathcal E$. It is natural to define \emph{function spaces} on $\mathcal G$ as direct sums of function spaces on each edge: for $p\in[1,\infty]$ and for $k=1,2$, we define
\[
L^p(\mathcal G)=\bigoplus_{e\in\mathcal E}L^p(I_e),\quad H^k(\mathcal G)=\bigoplus_{e\in\mathcal E}H^k(I_e).
\]
We denote by $\scalar{\cdot}{\cdot}$ the scalar product on $L^2(\mathcal G)$ and by $\dual{\cdot}{\cdot}$ the duality product on $H^1(\mathcal G)$. 
As no compatibility conditions have been given on the vertices yet, a function $\psi\in H^1(\mathcal G)$ has a priori multiple values on each of the vertices. For a vertex $v\in\mathcal V$, we denote by
\[
\psi(v)=(\psi_e(v))_{e\sim v}\in\R^{d_v}
\]
the vector of the values of $\psi$ at $v$, where $e\sim v$ denotes the \emph{edges incident} to $v$ and $d_v$ is the \emph{degree} of $v$, i.e. the number of edges  incident to $v$. In a similar way, for $\psi\in H^2(\mathcal G)$, we denote by
\[
\psi'(v)=(\psi'_e(v)) _{e\sim v}\in\R^{d_v}
\]
the vector of the outer derivatives of $\psi$ at the vertex $v$. For brevity in notation, we shall also note
\[
\psi(\mathcal V)=(\psi(v))_{v\in\mathcal V},\quad \psi'(\mathcal V)=(\psi'(v))_{v\in\mathcal V},
  \]
  the vectors constructed by the values of $\psi$ and $\psi'$ at each
  of the vertices.

\begin{figure}[htpb!]
  \centering
  \begin{tikzpicture}
    \coordinate (A) at (-3,0);
    \coordinate (O) at (0,0);
    \coordinate (B) at (2,-1);
    \coordinate (C) at (1.2,1.5);

    \node[label=left:$A$]  at (A) {$\bullet$};
    \node[label={[shift={(-0.2,-0.9)}]$O$}] at (O) {$\bullet$};
    \node[label=right:$B$] at (B)  {$\bullet$};
    \node[label=right:$C$] at (C) {$\bullet$};
    
    \draw (A) -- (O);
    \draw (B) -- (O);
    \draw (C) -- (O);
  \end{tikzpicture}  
  \caption{3-edges star graph $\mathcal{G}_{3,sg}$}
  \label{fig:3stargraph}
\end{figure}
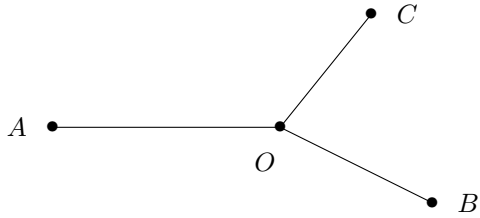

To give an example, we
consider the simple 3-edges star graph $\mathcal{G}_{3,sg}$ drawn on Figure~\ref{fig:3stargraph}.
The degree $d_O$ of the vertex $O$ is $d_O=3$, the set of vertices is $\mathcal{V}=\{O,A,B,C\}$ and the set of edges is
$\mathcal{E}=\{[OA],[OB],[OC]\}$. The vectors $\psi(O)$ and $\psi'(O)$ are given by
\[
\psi(O)=
\begin{pmatrix}
  \psi_{OA}(O)\\ \psi_{OB}(O)\\\psi_{OC}(O)
\end{pmatrix}, \qquad
  \psi'(O)=-
  \begin{pmatrix}
    \partial_{\mathbf{n}_{OA}}\psi_{OA}(O)\\
    \partial_{\mathbf{n}_{OB}}\psi_{OB}(O)\\
    \partial_{\mathbf{n}_{OC}}\psi_{OC}(O)
  \end{pmatrix},
\]
where $\mathbf{n}_{OM}=\mathbf{OM}/\|\mathbf{OM}\|$, $M\in \{A,B,C\}$, is the
inward unit vector, and
\[
  \partial_{\mathbf{n}_{OM}}\psi_{OM}(O)=\lim_{\substack {t\to 0\\t>0}}\frac{\psi_{OM}(O+t \mathbf{n}_{OM})-\psi_{OM}(O)}{t}.
\]

  A \emph{quantum graph} is a metric graph $\mathcal G$ equipped with a \emph{Hamiltonian} operator $H$, which is usually defined in the following way. The operator $H$ is a second order unbounded operator
  \[
H:D(H)\subset L^2(\mathcal G)\to L^2(\mathcal G),
\]
which  is such that for $u\in D(H)\subset H^2(\mathcal G)$ and for each edge $e\in\mathcal E$ we have 
  \begin{equation}
(Hu)_e=-u_e''.\label{eq:def-H}
\end{equation}
The domain $D(H)$  of $H$ is a subset of $H^2(\mathcal G)$ of functions verifying specific vertex compatibility conditions, described in the following way. At a vertex $v\in\mathcal V$, let $A_v,B_v$ be $d_v\times d_v$ matrices. The \emph{compatibility conditions} for $u\in H^2(\mathcal G)$ may then be described as
\[
A_vu(v)+B_vu'(v)=0.
\]
For the full set of vertices $\mathcal V$, we denote by
\[
  A=\operatorname{diag}(A_v,v\in\mathcal V),\quad
B=\operatorname{diag}(B_v,v\in\mathcal V)
\]
the matrices describing the compatibility conditions. The \emph{domain} $D(H)$ of $H$ is then given by
\begin{equation}
  \label{eq:bc}
D(H)=\{ u\in H^2(\mathcal G): Au(\mathcal V)+Bu'(\mathcal V)=0 \}. 
\end{equation}
We will assume that $A$ and $B$ are such that $H$ is self-adjoint, that is at each vertex $v$ the $d_v\times 2d_v$ augmented matrix $(A_v|B_v)$ has maximal rank and the matrix $A_vB_v^*$ is self-adjoint. 
Recall that (see e.g.~\cite{BeKu13}) the boundary conditions at a vertex $v\in\mathcal V$ may be reformulated using three orthogonal and mutually orthogonal operators $P_{D,v}$ ($D$ for Dirichlet), $P_{N,v}$ ($N$ for Neumann) and $P_{R,v}$ ($R$ for Robin) and an invertible self-adjoint operator $\Lambda_v:\mathbb C^{d_v}\to\mathbb C^{d_v}$ such that for each $u\in D(H)$ we have
\[
P_{D,v}u(v)=
P_{N,v}u'(v)=
\Lambda_{v}P_{R,v}u'(v)-P_{R,v}u(v)=0.
\]
The \emph{quadratic form} associated with
$H$ is then expressed as
\[
Q(u)=\frac12\dual{Hu}{u}=\frac12\sum_{e\in\mathcal E}\norm{u_e'}_{L^2}^2+\frac12\sum_{v\in\mathcal V}\scalar{P_{R,v}u}{\Lambda_vP_{R,v}u}_{\mathbb C^{d_v}},
\]
and its \emph{domain} is given by
\[
D(Q)=H^1_D(\mathcal G)=\{u\in  H^2(\mathcal G):P_{D,v}u=0,\;\forall v\in\mathcal V\}.
  \]
  Among the many possible vertex conditions, the \emph{Kirchhoff-Neumann condition} is the most frequently encountered. By analogy with Kirchhoff laws in electricity (preservation of charge and current), it consists  at a vertex $v$ to require:
  \[
u_e(v)=u_{e'}(v),\;\forall e,e'\sim v,\quad \sum_{e\sim v}u_e'(v)=0.
    \]
Another popular vertex condition is the $\delta$ or \emph{Dirac condition} of strength $\alpha_v\in\R$. It corresponds to continuity of the function at the vertex $v$, and a jump condition of size $\alpha_v$ on the derivatives, that is 
      \[
u_e(v)=u_{e'}(v),\;\forall e,e'\sim v,\quad \sum_{e\sim v}u_e'(v)=\alpha_v u(v),
\]
where we slightly change our notation to designate by $u(v)$ the common value of $u$ at $v$. For $\alpha_v=0$, we obviously recover the Kirchhoff-Neumann condition. If $\delta$ conditions are requested on each of the vertices of the graph, the quadratic form and its associated domain $H_D^1(\mathcal G)$ are given by
\[
Q(u)=\frac12\sum_{e\in\mathcal E}\norm{u_e'}_{L^2}^2+\frac12\sum_{v\in\mathcal V}\alpha_v|u(v)|^2,\quad H_D^1(\mathcal G)=\left\{ u\in H^1(\mathcal G):\forall v\in\mathcal V,\forall e,e'\sim v,\;u_{e}(v)=u_{e'}(v) \right\}.
\]

\subsection{Space discretization}

We present here the space discretization of the second order unbounded
operator $H$.
We discretize each edge $e\in\mathcal{E}$ with $N_e\in\mathbb{N}^*$ interior
points (when $e\in\mathcal{E}$ is semi-infinite, we choose a large but
finite length and we add an artificial terminal
vertex with appropriate - typically Dirichlet - boundary condition). We therefore obtain a uniform
discretization  $\{x_{e,k}\}_{0\leq k\leq N_e+1}$  of the edge $e$ that can be
assimilated to 
the interval $I_e = [0,l_e]$, \textit{i.e.} 
\begin{equation*}
x_{e,0} : = 0 < x_{e,1}< \cdots<x_{e,N_e}< x_{e,N_e+1} : = l_e,
\end{equation*}
with $x_{e,k+1} - x_{e,k} = l_e/(N_e+1) : = \delta x_e$ for $0\leq k\leq
N_e$ (see Figure~\ref{fig:disc_edge}). We denote by $v_1$ the vertex
at $x_{e,0}$, by $v_2$ the one at $x_{e,N_e+1}$
and, for any $u\in H^1_D(\mathcal G)$, for all $e\in\mathcal{E}$
and $1\leq k\leq N_e$, 
\begin{equation*}
u_{e,k} : = u_e (x_{e,k}),
\end{equation*}
as well as
\begin{equation*}
  u_{e,v} : =
  \begin{cases}
    u_e(x_{e,0})&\text{ if }v=v_1, \\
    u_e(x_{e,N_e+1})&\text{ if }v=v_2.
  \end{cases}
\end{equation*}
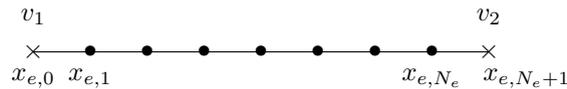
\begin{figure}[htpb!]
  \centering
  \begin{tikzpicture}
  \draw (-3,0) -- (3,0);
  \node[label=above:$v_1$] at (-3,0) {$\times$};
  \node[label=above:$v_2$] at (3,0) {$\times$};
  \node[below] at (-3,-0.1) {$x_{e,0}$};
  \node[below] at (-2.25,-0.1) {$x_{e,1}$};
  \node[below] at (2.25,-0.1) {$x_{e,N_e}$};
  \node[below] at (3,-0.1) {$\hspace{1cm}x_{e,N_e+1}$};
  \foreach \x in {-2.25,-1.5,...,2.25} {
    \node at (\x,0) {$\bullet$};
  };
\end{tikzpicture}  
  \caption{Discretization mesh of an edge $e\in \mathcal{E}$}
  \label{fig:disc_edge}
\end{figure}

We now assume that $N_e\geq 3$ and discretize the Laplacian operator on the interior of $e$, \textit{i.e.} we give an approximation of $Hu(x_{e,k})$ for $1\leq k\leq N_e$. Two cases need to be distinguished: the points closed to the boundary ($k=1,N_e$) and the other points. We shall start with the later.

Note that we do not discretized the Laplacian on the vertices, because, as will appear in a moment, the values of the functions at the vertices are determined in terms of the values at the interior nodes with the boundary conditions. 

For any $2\leq k\leq N_e-1$, the
second order approximation of the Laplace operator by finite differences on
$e\in \mathcal E$ is given by 
\begin{equation*}
H u(x_{e,k}) \approx -\frac{u_{e,k-1} - 2u_{e,k} + u_{e,k+1}}{{\delta x_e}^2}.
\end{equation*}
For the cases $k = 1$ and $k = N_e$ corresponding to the neighboring nodes of the
vertices $v_1$ and $v_2$, the approximation requires $u_{e,v_1}$
and $u_{e,v_2}$. We therefore use the boundary conditions
\[
  A_v u(v)+B_v u'(v)=0, \quad v \in \{v_1,v_2\},
\]
in order to evaluate them. To avoid any order reduction, we use second
order finite differences to approximate the outgoing
derivatives. Therefore, we
need the two closest neighboring nodes and for $-2\leq j\leq 0$, we denote 
\begin{align*}
u_{e,v_1,j} = u_{e}(x_{e,|j|})\quad\textrm{and}\quad u_{e,v_2,j} = u_{e}(x_{e,N_e+j+1}).
\end{align*}
The second order approximation of the outgoing derivative from $e$ at
$v\in \{v_1,v_2\}$ is given by
\begin{equation*}
u_e'(x_{e,v}) \approx (Du_{e,v})_0:=\frac{3 u_{e,v,0} - 4 u_{e,v,-1} + u_{e,v,-2}}{2\delta x_e}.
\end{equation*}
As a matter of fact, to increase precision, we have chosen in
  the implementation of the Grafidi library to use third order finite
  differences approximations for the derivatives at the vertex. This
  is transparent for the user and we restrict ourselves to second order in this presentation to increase readability.
We therefore have the approximation of the boundary conditions 
\begin{equation}
  \label{eq:bc01}
A_v [u_{v,0}] + B_v [Du_{v,0}]= 0,
\end{equation}
where $[u_{v,0}] = (u_{e,v,0})_{e\sim v}$ and $[Du_{v,0}] = ( (Du_{e,v})_0)_{e\sim
  v}$. We define the diagonal matrix $\Lambda \in
\mathbb{R}^{d_v\times d_v}$ with diagonal components by
\[
  \Lambda_{j,j} = \frac{1}{\delta x_{e_j}},\quad e_j \sim v, \quad j=1,\dots,d_{v}.
\]
Therefore, the approximate boundary condition~\eqref{eq:bc01} can be
rewritten as
\begin{equation}
  \label{eq:bc02}
  \left(A_v+\frac32 B_v \Lambda\right) [u_{v,0}] = 2 B_v \Lambda [u_{v,-1}] - \frac{1}{2} B_v
  \Lambda [u_{v,-2}],
\end{equation}
where $[u_{v,-1}] = (u_{e,v,-1})_{e\sim v}$ and $[u_{v,-2}] = (u_{e,v,-2})_{e\sim v}$.
Assuming that $A_v+\frac32 B_v \Lambda$ is invertible (which can be
done without loss of generality, see~\cite{BeDuLe20,BhBoHe20}), this is equivalent to
\begin{equation}
  \label{eq:bc03}
  \begin{aligned}
    [u_{v,0}] &= 2 \left(A_v+\frac32 B_v \Lambda\right)^{-1}B_v \Lambda [u_{v,-1}] - \frac{1}{2}
    \left(A_v+\frac32 B_v \Lambda\right)^{-1}B_v \Lambda [u_{v,-2}].
  \end{aligned}
\end{equation}
Solving the linear system~\eqref{eq:bc02} of size $d_v \times d_v$  allows to compute
the boundary 
values $[u_{v,0}]$ in terms of interior nodes. 
Thus, the value of $u_{e,v_1}$
(resp. $u_{e,v_2}$) depends linearly on the vectors $[u_{v_1,-1}]$ and
$[u_{v_1,-2}]$ (resp. $[u_{v_2,-1}]$ and $[u_{v_2,-2}]$) which take values from every edge connected to
the vertex $v_1$ (resp. $v_2$). It is then possible
to deduce an approximation of the Laplace operator at $x_{e,1}$ and
$x_{e,N_e}$. Indeed, from~\eqref{eq:bc03} there exist $(\alpha_{e,v})_{e\sim v}\in \mathbb{R}^{d_v}$,
for $v\in\{v_1,v_2\}$, which depend on every discretization parameter $\delta
x_e$ corresponding to the edges connected to $v$, such that 
\begin{equation*}
H u(x_{e,1}) \approx \frac{\displaystyle u_{e,2} - 2u_{e,1} + \sum_{e\sim v_1}\alpha_{e,v_1}(4 u_{e,v_1,-1} - u_{e,v_1,-2})}{{\delta x_e}^2},
\end{equation*}
and
\begin{equation*}
H u(x_{e,N_e}) \approx \frac{\displaystyle u_{e,N_e-1} - 2u_{e,N_e} + \sum_{e\sim v_2}\alpha_{e,v_2}(4 u_{e,v_2,-1} - u_{e,v_2,-2})}{{\delta x_e}^2}.
\end{equation*}
Since $(u_{e,v,j})_{-2\leq j\leq 0,v\in\{v_1,v_2\}}$ are interior mesh points
from the other edges, we limit our discretization to the interior mesh
points of the graph.
The approximated values of $u$ at each vertex
will be computed
using~\eqref{eq:bc03}. We denote $[{u}] =
(u_{e,k})_{1\leq k\leq N_e, e\in\mathcal{E}}$ the vector in $\mathbb{R}^{N}$,
with $N = \sum_{e\in\mathcal{E}} N_e$, representing the values of $u$ at each
interior mesh point of each edge of $\mathcal{G}$. We introduce the matrix
$[[{H}]]\in\mathbb{R}^{N\times N}$ corresponding to the
discretization of $H$ on the interior of each edge of the graph, which yields the
approximation 
\begin{equation*}
H u \approx [[{H}]]\, [{u}].
\end{equation*}

  To define discretized integrals on the graph, we proceed in the following way. We use the standard trapezoidal rule on each of the edges: on an edge $I_e$, for a vector $[u]$ (corresponding to a discretized function $u$) we approximate
  \[
\int_{I_e}u_e(x)dx\approx \mathcal I_e([u]):=\delta x_e\left(\sum_{k=0}^{N_e+1}u_{e,k}-\frac{u_{e,0}+u_{e,N_e+1}}{2}\right),
    \]
    where the terminal values $u_{e,0}$, $u_{e,N_e+1}$ are computed with~\eqref{eq:bc03}. The full integral is then approximated by
    \[
\int_{\mathcal G}u(x)dx\approx \sum_{e\in\mathcal E}\mathcal I_e([u]).
\]
This formula defines directly the discretization of $L^p(\mathcal G)$, that we denote $\ell^p(\mathcal G)$. 
      
\begin{figure}[htpb!]
    \includegraphics[width=0.26\textwidth]{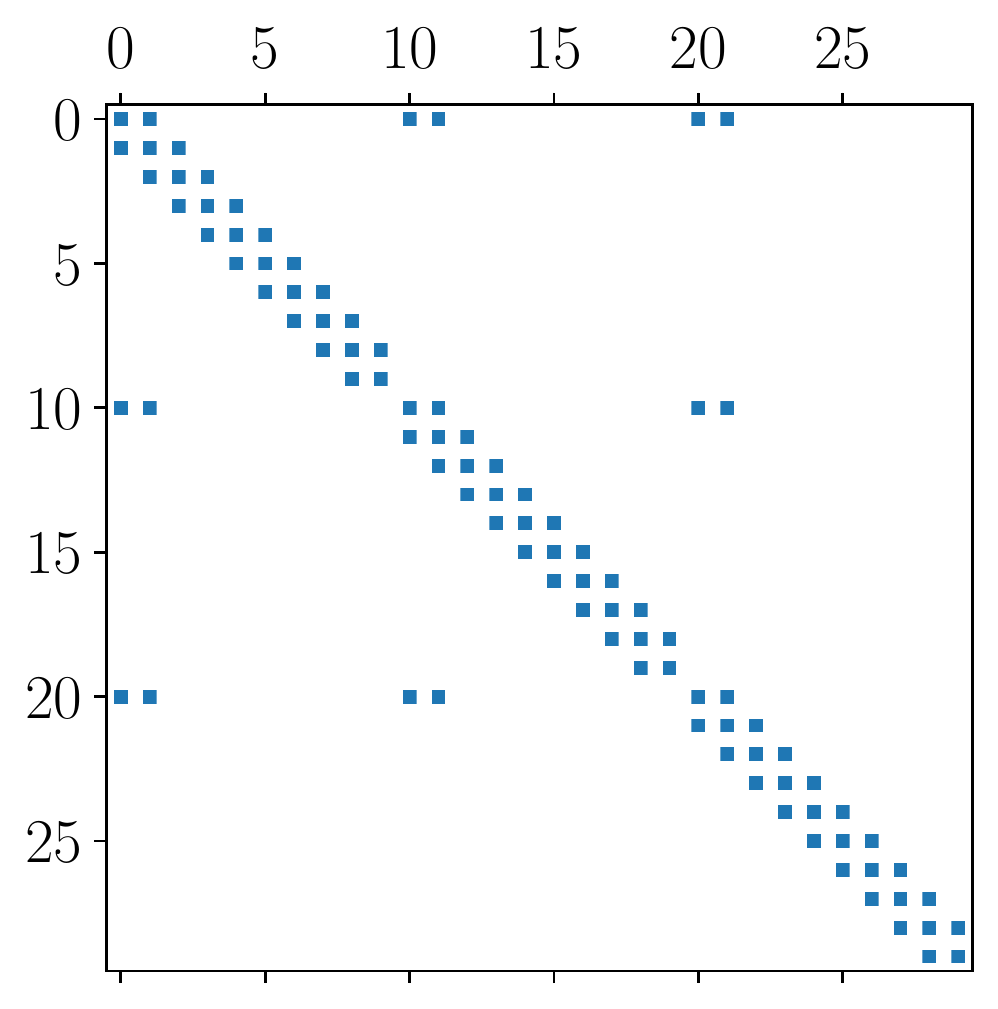}
    \caption{Matrix representation $[[H]]$ of $H$. }
\label{fig:f1f2}
\end{figure}

As an example, we consider the operator $H$ for the graph
$\mathcal{G}_{3,sg}$ of Figure~\ref{fig:3stargraph} with Dirichlet boundary conditions for the
exterior vertices $A$, $B$ and $C$ and Kirchhoff-Neumann conditions for the
central vertex $O$. We plot on Figure~\ref{fig:f1f2} the positions of the non
zero coefficients of the corresponding 
matrix $[[{H}]]$ when the discretization is such that $N_e = 10$, for
each $e\in\mathcal{E}$. The coefficients accounting for the Kirchhoff boundary
condition are the ones not belonging to the tridiagonal component of
the matrix.

\section{Some elements of the Grafidi library}
\label{sec:library}

\subsection{First steps with the Grafidi library}
\label{sec:first-steps}

We  introduce the Grafidi library by presenting some very basic
manipulations on an example: we describe the 
3-edges star graph $\mathcal{G}_{3,sg}$ drawn on Figure~\ref{fig:3stargraph} with $\mathcal{V}=\{O,A,B,C\}$ and 
$\mathcal{E}=\{[OA],[OB],[OC]\}$. We assume that the length of each edge is $10$.
Our goal in this simple example is to draw a
function $u$ that lives on the graph $\mathcal{G}_{3,sg}$, given by
\begin{equation}
  \label{eq:def-u}
  u(x)= 
    \begin{cases}
      e^{-x^2}& \text{for } x\in [OA],\\
      e^{-x^2} &\text{for }  x\in [OB],\\
      e^{-x^2}&\text{for }  x\in [OC].
    \end{cases}
  \end{equation}
The result is achieved using the code given in Listing~\ref{lst:RD}.
\begin{listing}[htbp!]
  \centering





  \begin{Verbatim}[commandchars=\\\{\},linenos=false,frame=single,framesep=2mm,fontsize=\footnotesize]
\PYG{k+kn}{import} \PYG{n+nn}{networkx} \PYG{k}{as} \PYG{n+nn}{nx}
\PYG{k+kn}{import} \PYG{n+nn}{numpy} \PYG{k}{as} \PYG{n+nn}{np}
\PYG{k+kn}{import} \PYG{n+nn}{matplotlib.pyplot} \PYG{k}{as} \PYG{n+nn}{plt}
\PYG{k+kn}{from} \PYG{n+nn}{Grafidi} \PYG{k+kn}{import} \PYG{n}{Graph} \PYG{k}{as} \PYG{n}{GR}
\PYG{k+kn}{from} \PYG{n+nn}{Grafidi} \PYG{k+kn}{import} \PYG{n}{WFGraph} \PYG{k}{as} \PYG{n}{WF}

\PYG{n}{g\PYGZus{}list}\PYG{o}{=}\PYG{p}{[}\PYG{l+s+s2}{\PYGZdq{}O A \PYGZob{}\PYGZsq{}Length\PYGZsq{}:10\PYGZcb{}\PYGZdq{}}\PYG{p}{,} \PYG{l+s+s2}{\PYGZdq{}O B \PYGZob{}\PYGZsq{}Length\PYGZsq{}:10\PYGZcb{}\PYGZdq{}}\PYG{p}{,} \PYG{l+s+s2}{\PYGZdq{}O C \PYGZob{}\PYGZsq{}Length\PYGZsq{}:10\PYGZcb{}\PYGZdq{}}\PYG{p}{]}
\PYG{n}{g\PYGZus{}nx} \PYG{o}{=} \PYG{n}{nx}\PYG{o}{.}\PYG{n}{parse\PYGZus{}edgelist}\PYG{p}{(}\PYG{n}{g\PYGZus{}list}\PYG{p}{,}\PYG{n}{create\PYGZus{}using}\PYG{o}{=}\PYG{n}{nx}\PYG{o}{.}\PYG{n}{MultiDiGraph}\PYG{p}{())}

\PYG{n}{g} \PYG{o}{=} \PYG{n}{GR}\PYG{p}{(}\PYG{n}{g\PYGZus{}nx}\PYG{p}{)}

\PYG{n}{fun} \PYG{o}{=} \PYG{p}{\PYGZob{}\PYGZcb{}}
\PYG{n}{fun}\PYG{p}{[(}\PYG{l+s+s1}{\PYGZsq{}O\PYGZsq{}}\PYG{p}{,} \PYG{l+s+s1}{\PYGZsq{}A\PYGZsq{}}\PYG{p}{,} \PYG{l+s+s1}{\PYGZsq{}0\PYGZsq{}}\PYG{p}{)]}\PYG{o}{=}\PYG{k}{lambda} \PYG{n}{x}\PYG{p}{:} \PYG{n}{np}\PYG{o}{.}\PYG{n}{exp}\PYG{p}{(}\PYG{o}{\PYGZhy{}}\PYG{n}{x}\PYG{o}{**}\PYG{l+m+mi}{2}\PYG{p}{)}
\PYG{n}{fun}\PYG{p}{[(}\PYG{l+s+s1}{\PYGZsq{}O\PYGZsq{}}\PYG{p}{,} \PYG{l+s+s1}{\PYGZsq{}B\PYGZsq{}}\PYG{p}{,} \PYG{l+s+s1}{\PYGZsq{}0\PYGZsq{}}\PYG{p}{)]}\PYG{o}{=}\PYG{k}{lambda} \PYG{n}{x}\PYG{p}{:} \PYG{n}{np}\PYG{o}{.}\PYG{n}{exp}\PYG{p}{(}\PYG{o}{\PYGZhy{}}\PYG{n}{x}\PYG{o}{**}\PYG{l+m+mi}{2}\PYG{p}{)}
\PYG{n}{fun}\PYG{p}{[(}\PYG{l+s+s1}{\PYGZsq{}O\PYGZsq{}}\PYG{p}{,} \PYG{l+s+s1}{\PYGZsq{}C\PYGZsq{}}\PYG{p}{,} \PYG{l+s+s1}{\PYGZsq{}0\PYGZsq{}}\PYG{p}{)]}\PYG{o}{=}\PYG{k}{lambda} \PYG{n}{x}\PYG{p}{:} \PYG{n}{np}\PYG{o}{.}\PYG{n}{exp}\PYG{p}{(}\PYG{o}{\PYGZhy{}}\PYG{n}{x}\PYG{o}{**}\PYG{l+m+mi}{2}\PYG{p}{)}

\PYG{n}{u} \PYG{o}{=} \PYG{n}{WF}\PYG{p}{(}\PYG{n}{fun}\PYG{p}{,}\PYG{n}{g}\PYG{p}{)}
\PYG{n}{\PYGZus{}} \PYG{o}{=} \PYG{n}{WF}\PYG{o}{.}\PYG{n}{draw}\PYG{p}{(}\PYG{n}{u}\PYG{p}{)}
\end{Verbatim}

  \caption{Simple Python example to draw a function on a 3-star graph $\mathcal{G}_{3,sg}$.}
  \label{lst:RD}
\end{listing}
We now describe each part of this simple example. The
functionalities of the Grafidi library rely on the following Python libraries:
networkx, numpy and matplotlib, which we first import.
The networkx library is mandatory and should be imported
after starting Python. Depending on the desire to make drawings and to make linear algebra operations,
it is recommended to import matplotlib and numpy. 
We then need to import the Grafidi library. It is made of two main classes: Graph and WFGraph, which we choose to import respectively as \verb+GR+ and \verb+WF+.

We then begin by creating a variable \verb+g_nx+, which is an instance of the
\verb+classes.multidigraph.MultiDiGraph+ of the networkx class. This
choice is motivated by the need of the description of a directed graph and the
possibility of multiple edges connecting the same  two nodes. Observe here that we have to choose an arbitrary orientation of the non-oriented graph for numerical purposes.
We
choose to describe the metric graph in the  Python list \verb+g_list+.
We identify each vertex by a Python string. Each element of \verb+g_list+ corresponds to an edge connecting two vertices. The length of each edge of the graph is defined with the keyword
\verb+Length+.

\begin{figure}[!htbp]
  \centering
  \includegraphics[width=.31\textwidth]{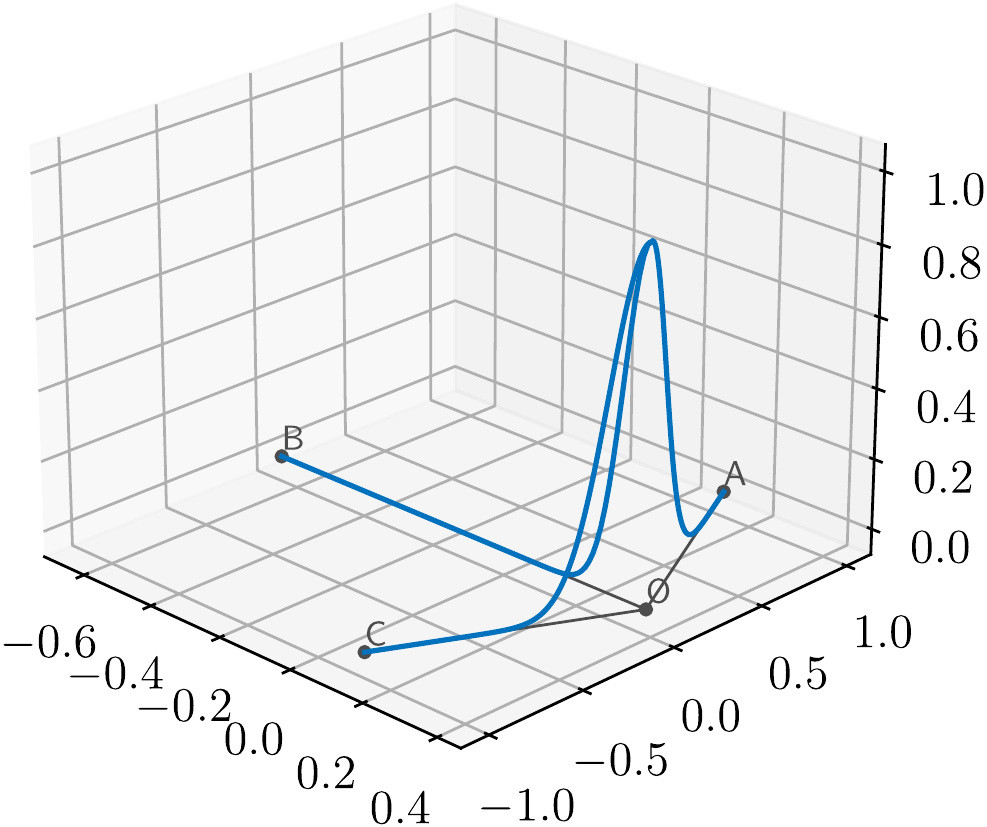}
  \caption{Plot of the function $u$ on graph $\mathcal {G}_{3,sg}$}
  \label{fig:simple_example_01}
\end{figure}

Then, we define the function that we wish to plot through a dictionary where each key corresponds to
an edge. The available keys can be found by the Python instruction \verb+g.Edges.keys()+. Each key
is a tuple made of three strings. The two first are the vertices labels defining the edge and the
third one is an identifier that will be explained later. The values are Python lambda functions with
$x$ belonging to the interval $[0,l_e]$, where $l_e$ is the length of the directed edge $e \in
\mathcal E$. So, $x=0$ corresponds to the initial vertex of $e$ and $x=l_e$ to the last one. We
construct an instance of the WFGraph class with the constructor \verb+WFGraph+ with as arguments the
dictionary \verb+fun+ and the instance of the graph \verb+g+. Since we import the class WFGraph as
\verb+WF+, the instruction may be shorten as it appears in the listing. It remains to use the
\verb+draw+ method of the WFGraph class to plot the function $u$ on $\mathcal {G}_{3,sg}$. The
result is available on Figure~\ref{fig:simple_example_01}. Since the \verb+draw+ function of the
WFGraph class delivers outputs, we use the Python instruction \verb+_ = + to avoid their display.

We use the networkx library to determine the geometric positions of each
vertex on the plane $(Oxy)$. More specifically, the function
\verb+networkx.drawing.layout.kamada_kawai_layout+ is executed on \verb+g_nx+ within
Graph class automatically to compute them. The length of each edge is however not taken
into account (indeed, the networkx library is implemented for non metric graphs). To overcome this issue, we have implemented the method
\verb+Position+ 
in Graph class. This method allows the
user to define by hand the geometric positions of each vertex. Its single argument is one
dictionary where the geometric position is given for each vertex. Finally, we 
draw the graph $\mathcal
{G}_{3,sg}$ with  the method \verb+draw+. For example, the definition of the geometric positions and the representation of the graph is proposed in
Listing~\ref{lst:geompos}. 

\begin{listing}[htbp!]
  \begin{Verbatim}[commandchars=\\\{\},linenos=false,frame=single,framesep=2mm,fontsize=\footnotesize]
\PYG{n}{NewPos}\PYG{o}{=}\PYG{p}{\PYGZob{}}\PYG{l+s+s1}{\PYGZsq{}O\PYGZsq{}}\PYG{p}{:[}\PYG{l+m+mi}{0}\PYG{p}{,}\PYG{l+m+mi}{0}\PYG{p}{],}\PYG{l+s+s1}{\PYGZsq{}A\PYGZsq{}}\PYG{p}{:[}\PYG{o}{\PYGZhy{}}\PYG{l+m+mi}{10}\PYG{p}{,}\PYG{l+m+mi}{0}\PYG{p}{],}\PYG{l+s+s1}{\PYGZsq{}B\PYGZsq{}}\PYG{p}{:[}\PYG{l+m+mi}{10}\PYG{p}{,}\PYG{l+m+mi}{0}\PYG{p}{],}\PYG{l+s+s1}{\PYGZsq{}C\PYGZsq{}}\PYG{p}{:[}\PYG{l+m+mi}{0}\PYG{p}{,}\PYG{l+m+mi}{10}\PYG{p}{]\PYGZcb{}}
\PYG{n}{GR}\PYG{o}{.}\PYG{n}{Position}\PYG{p}{(}\PYG{n}{g}\PYG{p}{,}\PYG{n}{NewPos}\PYG{p}{)}
\PYG{n}{\PYGZus{}} \PYG{o}{=} \PYG{n}{GR}\PYG{o}{.}\PYG{n}{draw}\PYG{p}{(}\PYG{n}{g}\PYG{p}{)}
\end{Verbatim}

  \caption{Definition of geometric positions of vertices}
  \label{lst:geompos}
\end{listing}

\begin{figure}[htpb!]
  \centering
  \begin{tabular}{cc}
    \includegraphics[width=.31\textwidth]{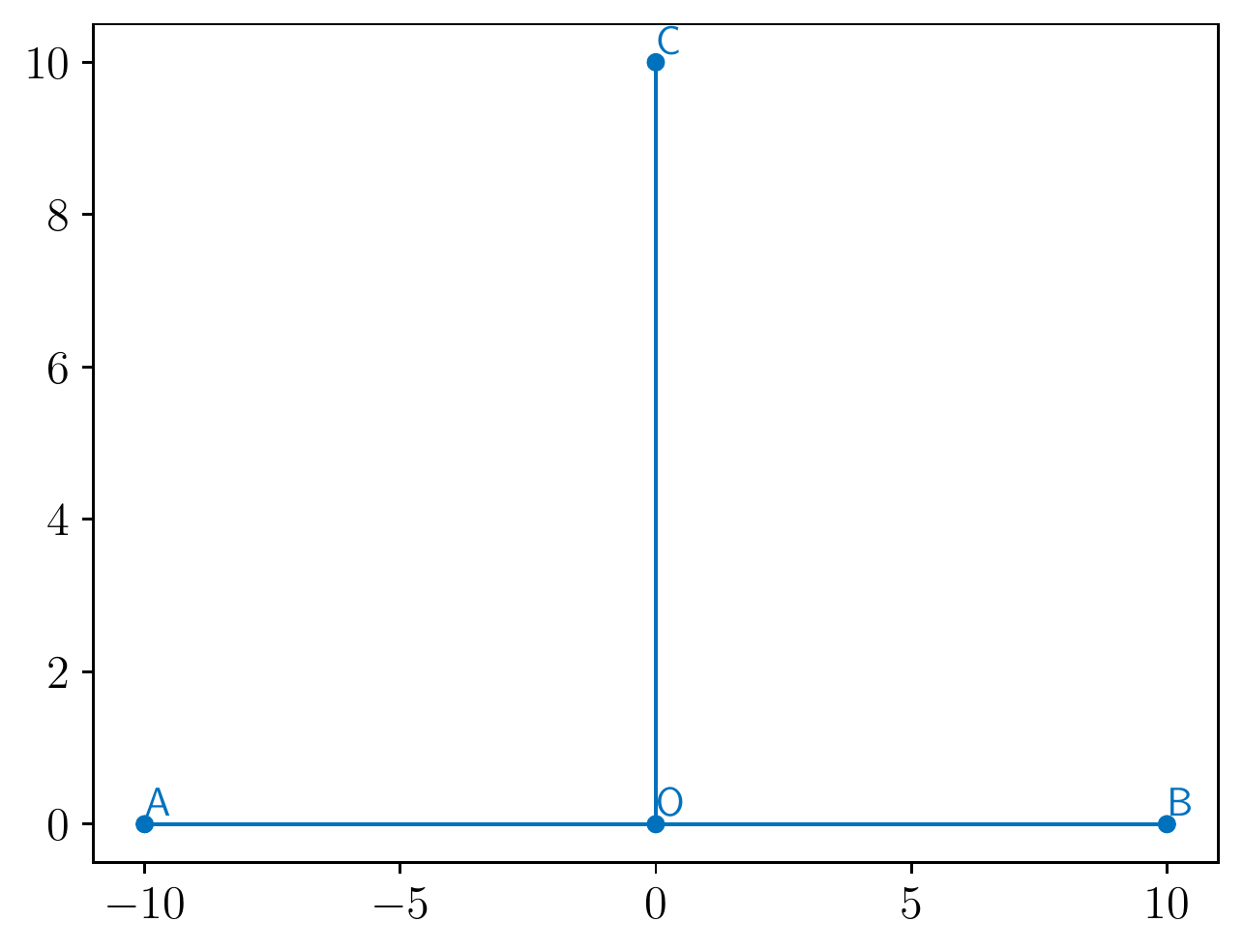} & \includegraphics[width=.31\textwidth]{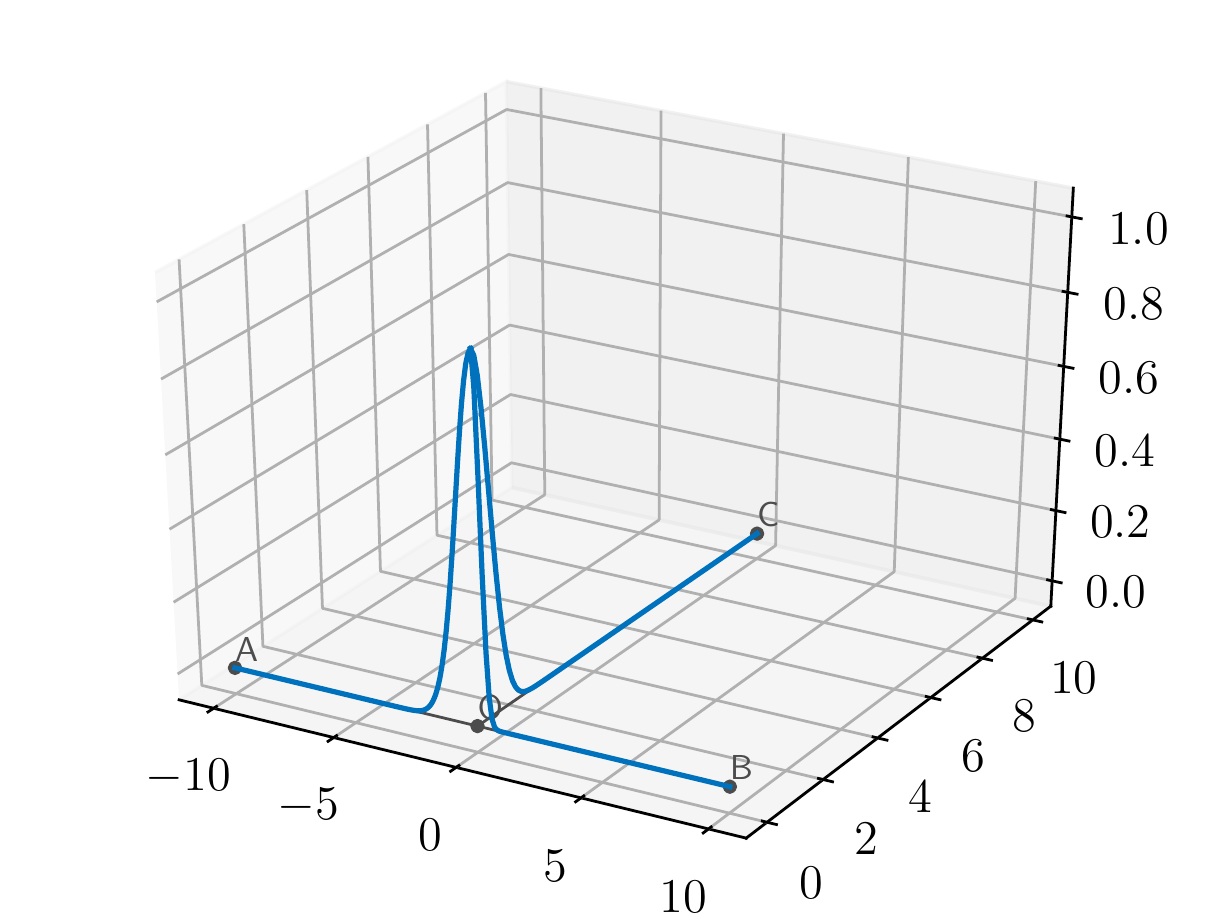}
  \end{tabular}
  \caption{Plot of the graph $\mathcal {G}_{3,sg}$ (left) and of the function $u$ on it (right).}
  \label{fig:simple_example_02}
\end{figure}

The new plot of the function $u$ on $\mathcal {G}_{3,sg}$ and the representation of the
graph are provided in Figure~\ref{fig:simple_example_02}.

\subsection{Basic elements of the Graph class}
The purpose of the Grafidi library is to provide tools to compute numerical solutions of partial differential equations involving the Laplace operator $H$ defined by~\eqref{eq:def-H}-\eqref{eq:bc}.
Actually, the instruction \verb+g = GR(g_nx)+ in Listing~\ref{lst:RD} automatically creates the discretization matrix $[[H]]$ of the operator $H$ following the rules defined in Section~\ref{sec:space_disc}. By default, the standard \textsl{Kirchhoff-Neumann conditions} are considered at each vertex and $N_e=100$ nodes are used to discretize each edge $e\in \mathcal{E}$.  The total number of discretization nodes is $N=\sum_{e\in \mathcal{E}}N_e$.  The matrix is stored in a sparse matrix in Compressed Sparse Column format in \verb+-g.Lap+ (actually,
  \verb+g.Lap+
  is the approximation matrix of $-H$). If needed, the user may declare other boundary conditions at each vertex. The boundary conditions at each vertex  are stored in a Python dictionary, which we call here \verb+bc+. Each key corresponds to a vertex and the values are lists. We provide in the Grafidi library various standard boundary conditions (Kirchhoff-Neumann, Dirichlet, $\delta$, $\delta'$), but more general can be constructed by defining matrices $A$ and $B$ at each vertex as in~\eqref{eq:bc}. We consider again the graph $\mathcal{G}_{3,sg}$ and assume that the space discretization has to be made with $3000$ interior nodes, and that boundary conditions are of homogeneous Dirichlet type at the vertices $A$, $B$ and $C$, and of $\delta$ type with strength $1$ for the vertex $O$. We therefore modify the instruction \verb+g = GR(g_nx)+ of Listing~\ref{lst:RD} to construct a new graph taking into account the new boundary conditions and total number of discretization points (see Listing~\ref{lst:def_bc_N}).
\begin{listing}[htbp!]
  \begin{Verbatim}[commandchars=\\\{\},linenos=false,frame=single,framesep=2mm,fontsize=\footnotesize]
\PYG{n}{bc} \PYG{o}{=} \PYG{p}{\PYGZob{}}\PYG{l+s+s1}{\PYGZsq{}O\PYGZsq{}}\PYG{p}{:[}\PYG{l+s+s1}{\PYGZsq{}Delta\PYGZsq{}}\PYG{p}{,}\PYG{l+m+mi}{1}\PYG{p}{],} \PYG{l+s+s1}{\PYGZsq{}A\PYGZsq{}}\PYG{p}{:[}\PYG{l+s+s1}{\PYGZsq{}Dirichlet\PYGZsq{}}\PYG{p}{],} \PYG{l+s+s1}{\PYGZsq{}B\PYGZsq{}}\PYG{p}{:[}\PYG{l+s+s1}{\PYGZsq{}Dirichlet\PYGZsq{}}\PYG{p}{],} \PYG{l+s+s1}{\PYGZsq{}C\PYGZsq{}}\PYG{p}{:[}\PYG{l+s+s1}{\PYGZsq{}Dirichlet\PYGZsq{}}\PYG{p}{]\PYGZcb{}}
\PYG{n}{N}\PYG{o}{=}\PYG{l+m+mi}{3000}
\PYG{n}{g} \PYG{o}{=} \PYG{n}{GR}\PYG{p}{(}\PYG{n}{g\PYGZus{}nx}\PYG{p}{,}\PYG{n}{N}\PYG{p}{,}\PYG{n}{bc}\PYG{p}{)}
\end{Verbatim}

  \caption{Definition of boundary conditions at vertices and discretization parameter}
  \label{lst:def_bc_N}
\end{listing}
Indeed, the constructor \verb+Graph+ actually takes three arguments: the mandatory instance of the networkx graph
\verb+g_nx+, and two optional arguments, the total number of discretization nodes \verb+N+ and the dictionary \verb+bc+ describing the boundary conditions at each vertex of graph $\mathcal{G}$.

\begin{figure}[htpb!]
  \centering
\begin{tikzpicture}[xscale=0.7]
  \draw[-] (-5,0) -- (5,0);
  \node[label=above:$B$] at (-5,0) {$\bullet$};
  \node[label=above:$A$] at (0,0) {$\bullet$};
  \node[label=above:$C$] at (5,0) {$\bullet$};
  \draw[<->,>=latex] (-5,-0.5) -- (0,-0.5);
  \draw[<->,>=latex] (0,-0.5) -- (5,-0.5);
  \node at (-2.5,-1) {$L$};
  \node at (2.5,-1) {$L$};
\end{tikzpicture}  
  \caption{Simple two edges star graph $\mathcal{G}_{2,sg}$}
  \label{fig:graphesimple}
\end{figure}
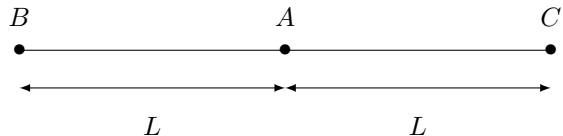

During the creation of the graph \verb+g+, two additional
variables are also automatically created: \verb+g.Edges+ and \verb+g.Nodes+. They allow to store informations
related to the mesh of $\mathcal G$. We describe them on the simple
two-edges star
graph $\mathcal{G}_{2,sg}$ (see Figure~\ref{fig:graphesimple}). It is made of three vertices $A$, $B$, $C$, $A$ being the central node, and 
two edges $[AB]$ and $[AC]$ with identical length $L$.

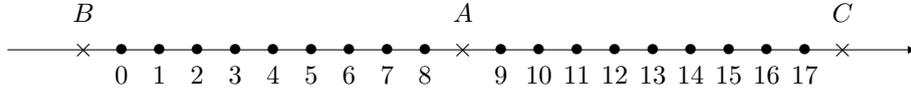
\begin{figure}[htpb!]
  \centering
\begin{tikzpicture}
  \draw[->,>=latex] (-6,0) -- (6,0);
  \node[label=above:$B$] at (-5,0) {$\times$};
  \node[label=above:$A$] at (0,0) {$\times$};
  \node[label=above:$C$] at (5,0) {$\times$};
  \foreach \x in {-4.5,-4.0,...,-0.5} {
     \node[] at (\x,0) {$\bullet$};
     \pgfmathtruncatemacro{\label}{2 * \x + 9}
         \node [below]  at (\x,-0.1) {\label};}
  \foreach \x in {0.5,1,...,4.5}{
     \node[] at (\x,0) {$\bullet$};
     \pgfmathtruncatemacro{\label}{2 * \x + 8}
         \node [below]  at (\x,-0.1) {\label};}
\end{tikzpicture}  
  \caption{Mesh on the simple star graph $\mathcal{G}_{2,sg}$.}
  \label{fig:maillage}
\end{figure}

We describe the mesh on the graph $\mathcal{G}_{2,sg}$. We assume that $L=5$ and we
discretize the graph with $N=18$ interior nodes. Thus,
$\delta x_i=\delta x=1/2$ and each edge is discretized with $N_i=9$, $i=1,2$,
nodes. The associated mesh is drawn on Figure~\ref{fig:maillage}.

The discretization nodes on the edge $[AB]$ are indexed from $0$ to $8$ and
the ones on $[AC]$ are indexed from
$9$ to $17$. All these informations are stored in
the dictionary \verb+g.Edges+. The keys are the edges of the graph made of the
vertices of each edge and a label (two vertices can be linked by many
edges). For the simple two-edges graph, the dictionary is given in Listing~\ref{lst:dict} (more detailed explanations  of the content of the dictionary is provided in the next section).
\begin{listing}[htbp!]
  \begin{Verbatim}[commandchars=\\\{\},linenos=false,frame=single,framesep=2mm,fontsize=\footnotesize]
Edges = \PYGZob{}
   (\PYGZsq{}B\PYGZsq{},\PYGZsq{}A\PYGZsq{},\PYGZsq{}0\PYGZsq{}) : \PYGZob{}\PYGZsq{}N\PYGZsq{}:9, \PYGZsq{}L\PYGZsq{}:5, \PYGZsq{}dx\PYGZsq{}:0.5, \PYGZsq{}Nodes\PYGZsq{}:[\PYGZsq{}B\PYGZsq{},\PYGZsq{}A\PYGZsq{}], \PYGZsq{}TypeC\PYGZsq{}:\PYGZsq{}S\PYGZsq{}, \PYGZsq{}Indexes\PYGZsq{}:[0,8]\PYGZcb{},
   (\PYGZsq{}A\PYGZsq{},\PYGZsq{}C\PYGZsq{},\PYGZsq{}0\PYGZsq{}) : \PYGZob{}\PYGZsq{}N\PYGZsq{}:9, \PYGZsq{}L\PYGZsq{}:5, \PYGZsq{}dx\PYGZsq{}:0.5, \PYGZsq{}Nodes\PYGZsq{}:[\PYGZsq{}A\PYGZsq{},\PYGZsq{}C\PYGZsq{}], \PYGZsq{}TypeC\PYGZsq{}:\PYGZsq{}S\PYGZsq{}, \PYGZsq{}Indexes\PYGZsq{}:[9,17]\PYGZcb{}
\PYGZcb{}
\end{Verbatim}

  \caption{The dictionary \texttt{g.Edges} }
\label{lst:dict}
\end{listing}

The second important variable is the dictionary \verb+g.Nodes+ that contains
various important informations to build the finite differences
approximation of the
operator $H$ on $\mathcal {G}$. The keys of \verb+g.Nodes+ are the identifiers for
each vertex. For the simple 2-star graph, they are \verb+'A'+, \verb+'B'+ and
\verb+'C'+. We associate to each vertex a dictionary with various
keys. We describe below the most relevant keys.
\begin{itemize}
\item \verb+'Degree'+ is an integer containing the degree $d_v$ of the vertex $v$.
\item \verb+'Boundary conditions'+ is a string containing the boundary
  condition set on the
  vertex $v$. The current possibilities are
  \begin{itemize}
  \item \verb+['Dirichlet']+,
  \item \verb+['Kirchhoff']+,
  \item \verb+['Delta', val]+, where \verb+val+ is the characteristic value of the $\delta$ condition,
  \item \verb+['Delta Prime', val]+, where \verb+val+ is the characteristic value of  the $\delta'$ condition,
  \item \verb+['UserDefined', [A_v,B_v]]+, where \verb+[A_v,B_v]+ are matrices used to describe the boundary condition at the vertex $v$.
  \end{itemize}
\item \verb+'Position'+ is a list $[x,y]$ representing the geometric coordinates
  of the vertex $v$.
\end{itemize}

We already met the method \verb+draw+ of Graph class. Some options are available to control figure name, color, width, markersize, textsize of the drawing of the graph (for a complete description, see Appendix). The method \verb+draw+ returns figure and axes matplotlib identifiers. This allows to have a fine control of the figure and its contained elements with matplotlib primitives.

\subsection{A first concrete example: eigenelements of the triple-bridge 
}
We are now able to handle more complex graphs. Since the Grafidi library relies on
the MultiDiGraph - \textsl{Directed
graphs with self loops and parallel edges} - class of networkx library, we can
handle loops and many edges between two single vertices. The declaration of such complex graphs is easy with the Graph class, as we illustrate in the following example.

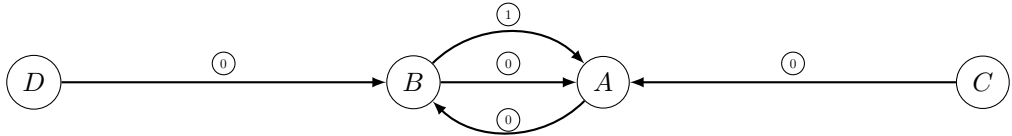
\begin{figure}[htpb!]
  \centering
  \begin{tikzpicture}[scale=.5]
    \node[draw,circle] (A) at (0,0) {$A$};
    \node[draw,circle] (B) at (-5,0) {$B$};
    \node[draw,circle] (D) at (-15,0) {$D$};
    \node[draw,circle] (C) at (10,0) {$C$};
    \draw[thick,->,>=latex] (D) -- (B);
    \draw[thick,->,>=latex] (B) -- (A);
    \draw[thick,->,>=latex] (B) to[bend left=45] (A);
    \draw[thick,->,>=latex] (A) to[bend left=45] (B);
    \draw[thick,->,>=latex] (C) -- (A);
    \draw (-10,0.5) circle [radius=0.3] node[scale=0.5] {$0$};    
    \draw (-2.5,-1.0) circle [radius=0.3] node[scale=0.5] {$0$};    
    \draw (-2.5,0.5) circle [radius=0.3] node[scale=0.5] {$0$};    
    \draw (-2.5,1.8) circle [radius=0.3] node[scale=0.5] {$1$};    
    \draw (5,0.5) circle [radius=0.3] node[scale=0.5] {$0$};    
  \end{tikzpicture}
  \caption{Directed graph $\mathcal{G}_d$}
  \label{fig:dirgraph}
\end{figure}

We want to represent the graph
$\mathcal{G}_d$ defined on Figure~\ref{fig:dirgraph}.

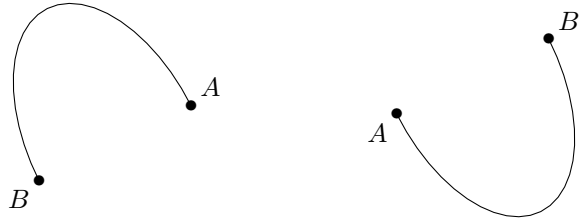
\begin{figure}[htpb!]

  \centering
  \begin{tabular}{ccc}
  \begin{tikzpicture}
  \draw [domain=0:1] plot ({sqrt(5/4)*cos(pi*\x r)*cos(atan(1/2))-2*sin(pi*\x
    r)*sin(atan(1/2))+2},{sqrt(5/4)*cos(pi*\x r)*sin(atan(1/2))+2*sin(pi*\x
    r)*cos(atan(1/2))+3/2});
  \node (A) at (3,2) {$\bullet$};
  \node (B) at (1,1) {$\bullet$};
  \node[above right] at (A) {$A$};
  \node[below left] at (B) {$B$};
\end{tikzpicture}  & \hspace{1cm} &
  \begin{tikzpicture}
  \draw [domain=0:1] plot ({sqrt(5/4)*cos(-pi*\x r)*cos(atan(1/2))-2*sin(-pi*\x
    r)*sin(atan(1/2))+2},{sqrt(5/4)*cos(-pi*\x r)*sin(atan(1/2))+2*sin(-pi*\x
    r)*cos(atan(1/2))+3/2});
  \node (A) at (3,2) {$\bullet$};
  \node (B) at (1,1) {$\bullet$};
  \node[above right] at (A) {$B$};
  \node[below left] at (B) {$A$};
\end{tikzpicture}  
  \end{tabular}
  \caption{The two configurations ``up'' and ``down'' of the oriented curved edge $\protect\overset{\curvearrowright}{AB}$}
  \label{fig:exellip}
\end{figure}

The vertices $A$ and $B$ are connected by three edges:
\begin{itemize}
\item one edge oriented from the vertex $A$ to the vertex $B$,
\item two edges oriented from the vertex $B$ to the vertex $A$.
\end{itemize}
The vertices $A$ and $B$ are also respectively connected to the vertices $C$ and $D$.
We provide in Listing~\ref{lst:compgraph} an example of a standard declaration.
\begin{listing}[htbp!]
  \begin{Verbatim}[commandchars=\\\{\},linenos=false,frame=single,framesep=2mm,fontsize=\footnotesize]
g\PYGZus{}list=[\PYGZdq{}B A \PYGZob{}\PYGZsq{}Length\PYGZsq{}:5\PYGZcb{}\PYGZdq{}, \PYGZdq{}B A \PYGZob{}\PYGZsq{}Length\PYGZsq{}:10\PYGZcb{}\PYGZdq{}, \PYGZdq{}A B \PYGZob{}\PYGZsq{}Length\PYGZsq{}:10\PYGZcb{}\PYGZdq{}, \PYGZdq{}C A \PYGZob{}\PYGZsq{}Length\PYGZsq{}:20\PYGZcb{}\PYGZdq{}, \PYGZdq{}D B \PYGZob{}\PYGZsq{}Length\PYGZsq{}:20\PYGZcb{}\PYGZdq{}]
\end{Verbatim}

  \caption{Declaration of complex graph}
  \label{lst:compgraph}
\end{listing}
The library automatically assigns an \verb+Id+ (such as the ones indicated in Figure~\ref{fig:dirgraph}) and  a type ``segment'' \verb+'S'+ or ``curve'' \verb+'C'+ to each edge. This operation is transparent for the user. If only one edge connects two vertices, the \verb+Id+ is set to \verb+'0'+ and the chosen type is \verb+'S'+. On the contrary, the algorithm chooses between \verb+'S'+ and \verb+'C'+ and the \verb+Id+ is incrementally increased starting from \verb+'0'+ when multiple edges connect the same two vertices. If a selected  edge is of type \verb+'C'+, it will be represented as curved line $\overset{\curvearrowright}{AB}$ (actually an half-ellipsis of length \verb+Length+) going from $A$ to $B$ counterclockwise (as a consequence, the edge will be ``up'' or ``down'' depending on the position of the
  vertices, see Figure~\ref{fig:exellip}).
  
The user can also explicitly provide the edge type and \verb+Id+ informations as in Listing~\ref{lst:txtfile}.
\begin{listing}[htbp!]
  \begin{Verbatim}[commandchars=\\\{\},linenos=false,frame=single,framesep=2mm,fontsize=\footnotesize]
g\PYGZus{}list=[\PYGZdq{}B A \PYGZob{}\PYGZsq{}Length\PYGZsq{}: 5,\PYGZsq{}Line\PYGZsq{}:\PYGZsq{}S\PYGZsq{},\PYGZsq{}Id\PYGZsq{}:\PYGZsq{}0\PYGZsq{}\PYGZcb{}\PYGZdq{},\PYGZbs{}
        \PYGZdq{}B A \PYGZob{}\PYGZsq{}Length\PYGZsq{}:10,\PYGZsq{}Line\PYGZsq{}:\PYGZsq{}C\PYGZsq{},\PYGZsq{}Id\PYGZsq{}:\PYGZsq{}1\PYGZsq{}\PYGZcb{}\PYGZdq{},\PYGZbs{}
        \PYGZdq{}A B \PYGZob{}\PYGZsq{}Length\PYGZsq{}:10,\PYGZsq{}Line\PYGZsq{}:\PYGZsq{}C\PYGZsq{},\PYGZsq{}Id\PYGZsq{}:\PYGZsq{}0\PYGZsq{}\PYGZcb{}\PYGZdq{},\PYGZbs{}
        \PYGZdq{}C A \PYGZob{}\PYGZsq{}Length\PYGZsq{}:20,\PYGZsq{}Line\PYGZsq{}:\PYGZsq{}S\PYGZsq{},\PYGZsq{}Id\PYGZsq{}:\PYGZsq{}0\PYGZsq{}\PYGZcb{}\PYGZdq{},\PYGZbs{}
        \PYGZdq{}D B \PYGZob{}\PYGZsq{}Length\PYGZsq{}:20,\PYGZsq{}Line\PYGZsq{}:\PYGZsq{}S\PYGZsq{},\PYGZsq{}Id\PYGZsq{}:\PYGZsq{}0\PYGZsq{}\PYGZcb{}\PYGZdq{}]
\end{Verbatim}

  \caption{User defined description of the graph $\mathcal{G}_d$.}
\label{lst:txtfile}
\end{listing}

\begin{figure}[htpb!]
  \centering
  \includegraphics[width=.6\textwidth]{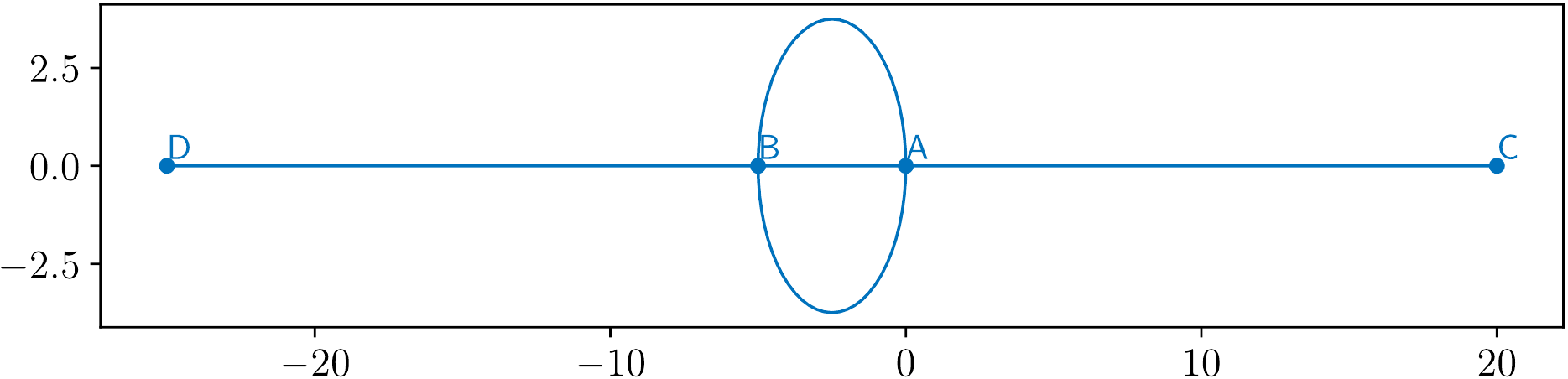}
  \caption{Plot of the graph $\mathcal{G}_d$ with Grafidi library.}
  \label{fig:graph_classe}
\end{figure}

The plot with the Grafidi library of the graph $\mathcal{G}_d$ with positions adjusted is presented on Figure~\ref{fig:graph_classe}.

As an illustration, we now present the computations of some eigenelements of the operator $H$ on
the graph $\mathcal{G}_d$ with Kirchhoff boundary conditions at the vertices $A$ and
$B$ and Dirichlet ones at the vertices $C$ and $D$. Since the approximation matrix of $\partial_{xx}$ is
automatically {generated} and stored in \verb+g.Lap+, we can compute the
eigenelements of $[[{H}]]=-$\verb+g.Lap+. We present in Listing~\ref{lst:eigenelements} the easiest way to compute the first four
eigenvalues/eigenvectors and to draw the eigenvectors on $\mathcal{G}_d$. It is understood that all libraries appearing in Listing~\ref{lst:RD} are already imported.
\begin{listing}[htbp!]



  \begin{Verbatim}[commandchars=\\\{\},linenos=false,frame=single,framesep=2mm,fontsize=\footnotesize]
\PYG{k+kn}{import} \PYG{n+nn}{scipy.sparse} \PYG{k}{as} \PYG{n+nn}{scs}

\PYG{n}{g\PYGZus{}list}\PYG{o}{=}\PYG{p}{[}\PYG{l+s+s2}{\PYGZdq{}B A \PYGZob{}\PYGZsq{}Length\PYGZsq{}:5\PYGZcb{}\PYGZdq{}}\PYG{p}{,} \PYG{l+s+s2}{\PYGZdq{}B A \PYGZob{}\PYGZsq{}Length\PYGZsq{}:10\PYGZcb{}\PYGZdq{}}\PYG{p}{,} \PYG{l+s+s2}{\PYGZdq{}A B \PYGZob{}\PYGZsq{}Length\PYGZsq{}:10\PYGZcb{}\PYGZdq{}}\PYG{p}{,} \PYG{l+s+s2}{\PYGZdq{}C A \PYGZob{}\PYGZsq{}Length\PYGZsq{}:20\PYGZcb{}\PYGZdq{}}\PYG{p}{,} \PYG{l+s+s2}{\PYGZdq{}D B \PYGZob{}\PYGZsq{}Length\PYGZsq{}:20\PYGZcb{}\PYGZdq{}}\PYG{p}{]}
\PYG{n}{g\PYGZus{}nx} \PYG{o}{=} \PYG{n}{nx}\PYG{o}{.}\PYG{n}{parse\PYGZus{}edgelist}\PYG{p}{(}\PYG{n}{g\PYGZus{}list}\PYG{p}{,}\PYG{n}{create\PYGZus{}using}\PYG{o}{=}\PYG{n}{nx}\PYG{o}{.}\PYG{n}{MultiDiGraph}\PYG{p}{())}
\PYG{n}{g} \PYG{o}{=} \PYG{n}{GR}\PYG{p}{(}\PYG{n}{g\PYGZus{}nx}\PYG{p}{)}

\PYG{n}{bc} \PYG{o}{=} \PYG{p}{\PYGZob{}}\PYG{l+s+s1}{\PYGZsq{}A\PYGZsq{}}\PYG{p}{:[}\PYG{l+s+s1}{\PYGZsq{}Kirchhoff\PYGZsq{}}\PYG{p}{],} \PYG{l+s+s1}{\PYGZsq{}B\PYGZsq{}}\PYG{p}{:[}\PYG{l+s+s1}{\PYGZsq{}Kirchhoff\PYGZsq{}}\PYG{p}{],} \PYG{l+s+s1}{\PYGZsq{}C\PYGZsq{}}\PYG{p}{:[}\PYG{l+s+s1}{\PYGZsq{}Dirichlet\PYGZsq{}}\PYG{p}{],} \PYG{l+s+s1}{\PYGZsq{}D\PYGZsq{}}\PYG{p}{:[}\PYG{l+s+s1}{\PYGZsq{}Dirichlet\PYGZsq{}}\PYG{p}{]\PYGZcb{}}
\PYG{n}{N}\PYG{o}{=}\PYG{l+m+mi}{3000}
\PYG{n}{g} \PYG{o}{=} \PYG{n}{GR}\PYG{p}{(}\PYG{n}{g\PYGZus{}nx}\PYG{p}{,}\PYG{n}{N}\PYG{p}{,}\PYG{n}{bc}\PYG{p}{)}
\PYG{n}{NewPos}\PYG{o}{=}\PYG{p}{\PYGZob{}}\PYG{l+s+s1}{\PYGZsq{}A\PYGZsq{}}\PYG{p}{:[}\PYG{l+m+mi}{0}\PYG{p}{,}\PYG{l+m+mi}{0}\PYG{p}{],}\PYG{l+s+s1}{\PYGZsq{}B\PYGZsq{}}\PYG{p}{:[}\PYG{o}{\PYGZhy{}}\PYG{l+m+mi}{5}\PYG{p}{,}\PYG{l+m+mi}{0}\PYG{p}{],}\PYG{l+s+s1}{\PYGZsq{}C\PYGZsq{}}\PYG{p}{:[}\PYG{l+m+mi}{20}\PYG{p}{,}\PYG{l+m+mi}{0}\PYG{p}{],}\PYG{l+s+s1}{\PYGZsq{}D\PYGZsq{}}\PYG{p}{:[}\PYG{o}{\PYGZhy{}}\PYG{l+m+mi}{25}\PYG{p}{,}\PYG{l+m+mi}{0}\PYG{p}{]\PYGZcb{}}
\PYG{n}{GR}\PYG{o}{.}\PYG{n}{Position}\PYG{p}{(}\PYG{n}{g}\PYG{p}{,}\PYG{n}{NewPos}\PYG{p}{)}

\PYG{p}{[}\PYG{n}{EigVals}\PYG{p}{,} \PYG{n}{EigVecs}\PYG{p}{]} \PYG{o}{=} \PYG{n}{scs}\PYG{o}{.}\PYG{n}{linalg}\PYG{o}{.}\PYG{n}{eigs}\PYG{p}{(}\PYG{o}{\PYGZhy{}}\PYG{n}{g}\PYG{o}{.}\PYG{n}{Lap}\PYG{p}{,}\PYG{n}{k}\PYG{o}{=}\PYG{l+m+mi}{4}\PYG{p}{,}\PYG{n}{sigma}\PYG{o}{=}\PYG{l+m+mi}{0}\PYG{p}{)}
\PYG{n}{Fig}\PYG{o}{=}\PYG{n}{plt}\PYG{o}{.}\PYG{n}{figure}\PYG{p}{(}\PYG{n}{figsize}\PYG{o}{=}\PYG{p}{[}\PYG{l+m+mi}{9}\PYG{p}{,}\PYG{l+m+mi}{6}\PYG{p}{])}
\PYG{k}{for} \PYG{n}{k} \PYG{o+ow}{in} \PYG{n+nb}{range}\PYG{p}{(}\PYG{n}{EigVals}\PYG{o}{.}\PYG{n}{size}\PYG{p}{):}
    \PYG{n}{ax}\PYG{o}{=}\PYG{n}{Fig}\PYG{o}{.}\PYG{n}{add\PYGZus{}subplot}\PYG{p}{(}\PYG{l+m+mi}{2}\PYG{p}{,}\PYG{l+m+mi}{2}\PYG{p}{,}\PYG{n}{k}\PYG{o}{+}\PYG{l+m+mi}{1}\PYG{p}{,}\PYG{n}{projection}\PYG{o}{=}\PYG{l+s+s1}{\PYGZsq{}3d\PYGZsq{}}\PYG{p}{)}
    \PYG{n}{EigVec} \PYG{o}{=} \PYG{n}{WF}\PYG{p}{(}\PYG{n}{np}\PYG{o}{.}\PYG{n}{real}\PYG{p}{(}\PYG{n}{EigVecs}\PYG{p}{[:,}\PYG{n}{k}\PYG{p}{]),}\PYG{n}{g}\PYG{p}{)}
    \PYG{n}{EigVec} \PYG{o}{=} \PYG{n}{EigVec}\PYG{o}{/}\PYG{n}{WF}\PYG{o}{.}\PYG{n}{norm}\PYG{p}{(}\PYG{n}{EigVec}\PYG{p}{,}\PYG{l+m+mi}{2}\PYG{p}{)}
    \PYG{n}{\PYGZus{}}\PYG{o}{=}\PYG{n}{WF}\PYG{o}{.}\PYG{n}{draw}\PYG{p}{(}\PYG{n}{EigVec}\PYG{p}{,}\PYG{n}{AxId}\PYG{o}{=}\PYG{n}{ax}\PYG{p}{)}
    \PYG{n}{ax}\PYG{o}{.}\PYG{n}{set\PYGZus{}title}\PYG{p}{(}\PYG{l+s+sa}{r}\PYG{l+s+s1}{\PYGZsq{}\PYGZdl{}\PYGZbs{}lambda\PYGZus{}}\PYG{l+s+si}{\PYGZob{}\PYGZcb{}}\PYG{l+s+s1}{=\PYGZdl{}\PYGZsq{}}\PYG{o}{.}\PYG{n}{format}\PYG{p}{(}\PYG{n}{k}\PYG{p}{)}\PYG{o}{+}\PYG{l+s+sa}{f}\PYG{l+s+s1}{\PYGZsq{}}\PYG{l+s+si}{\PYGZob{}}\PYG{n}{np}\PYG{o}{.}\PYG{n}{real}\PYG{p}{(}\PYG{n}{EigVals}\PYG{p}{[}\PYG{n}{k}\PYG{p}{])}\PYG{l+s+si}{:}\PYG{l+s+s1}{f}\PYG{l+s+si}{\PYGZcb{}}\PYG{l+s+s1}{\PYGZsq{}}\PYG{p}{)}
\end{Verbatim}

\caption{Computation of some eigenelements of $[[{H}]]$ on $\mathcal{G}_d$}
\label{lst:eigenelements}
\end{listing}

\begin{figure}[!htbp]
  \centering
  \includegraphics[width=\textwidth]{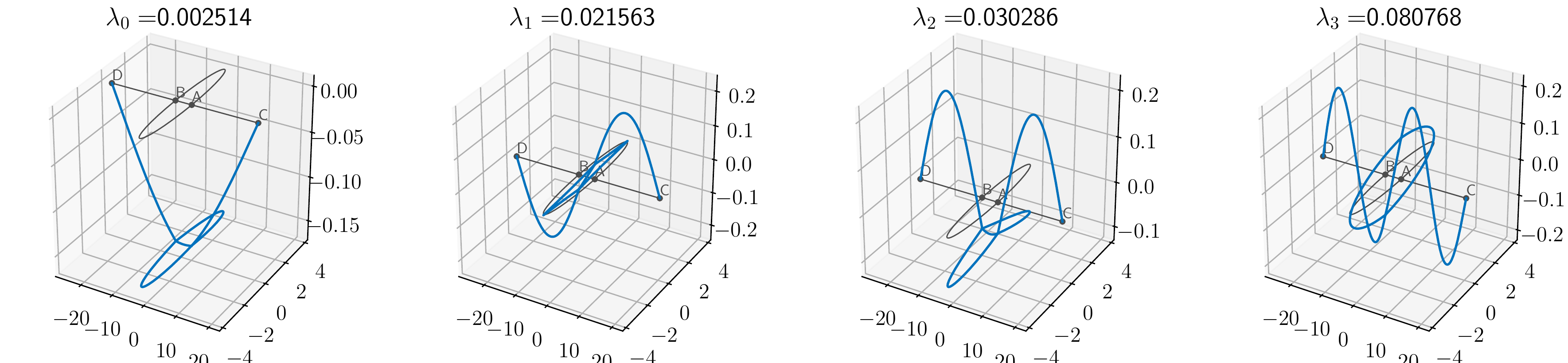}
  \caption{The first four eigenvectors of $[[{H}]]$ on graph $\mathcal{G}_d$. 
  }
  \label{fig:eigenelem}
\end{figure}

Listing~\ref{lst:eigenelements} works as follows. To compute the eigenelements of
$[[{H}]]$, we use the function \verb+linalg.eigs+ of the library \verb+scipy.sparse+. We transform
each eigenfunction (stored in the matrix \verb+EigVecs+) as an instance of the WFGraph
class by the instruction \verb+EigVec = WF(np.real(EigVecs[:,k]),g)+, where \verb+g+ is the graph
instance of Graph class representing $\mathcal{G}_d$. Next, we normalize the
eigenfunction. One notices that the $L^2$ norm of an instance of WFGraph can be simply computed with
the instruction \verb+WF.norm(EigVec,2)+. We are also able to divide a WFGraph entity by a scalar
(\verb+EigVec/WF.norm(EigVec,2)+). Each eigenvector is finally plotted with the command
\verb+WF.draw(EigVec,AxId=ax)+. The option \verb+AxId+ allows to plot the eigenvector on
the matplotlib axes \verb+ax+. The fours eigenvectors with their associated eigenvalues $\lambda_j$
are represented in Figure~\ref{fig:eigenelem}.

\section{Numerical methods for stationary and time dependent Schr\"odinger equations}
\label{sec:numerical-methods}

In this section, we discuss the implementation with the Grafidi library of various methods to compute grounds states or dynamical solutions of time-dependent Schrödinger equations on nonlinear quantum graphs.

\subsection{Computation of ground states on quantum graphs}

We begin with the computation of ground states. For a given second order differential operator $H$ on a quantum graph $\mathcal{G}$, a ground state is a minimizer of the Schr\"odinger energy $E$ at fixed
mass $M$, where
\[
E(u)=\frac12\dual{Hu}{u}-\frac12\int_{\mathcal G}G(|u|^2)dx,\quad G'=g,\quad M(u)=\norm{u}_{L^2(\mathcal G)}^2,
  \]
  where $g$ is the nonlinearity. In the following, we consider the case of a power-type nonlinearity
  \[
    g(u) = |u|^{p-1}u, \quad p>1.
  \]
To compute ground states, the most common methods are gradient methods. Here, we will cover two popular gradient methods: the Continuous Normalized Gradient Flow (CNGF), which we have analyzed in the context of quantum graphs in~\cite{BeDuLe20}, and a nonlinear (preconditionned) conjugate gradient flow (see~\cite{antoine2017efficient,danaila2017computation}), which we implement in the particular context of graphs without further theoretical analysis.

\subsubsection{The continuous normalized gradient flow}

We start with the CNGF method. We fix $\delta t>0$ a certain gradient step and $m>0$
the mass of the ground state. Let $\rho=\sqrt{m}$ be the $L^2$-norm of the ground
state. The method is divided into two steps: first a semi-implicit gradient descent step then a projection on the constraint manifold (here the $L^2$-sphere of radius $\rho$). In practice, we construct a sequence $\{u^n\}_{n\geq0}$ (which will converge to the ground state) given by
\begin{equation*}
\left\{\begin{array}{ll}
u_*^{n+1} = u^n - \delta t\left( Hu_*^{n+1} - |u^n|^{p-1}u_*^{n+1}\right),
\\ u^{n+1} = \rho u_*^{n+1}/\|u_*^{n+1}\|_{L^2(\mathcal{G})},
\end{array}\right.
\end{equation*}
where the initial data $u^0\in L^{2}(\mathcal{G})$ is chosen such that $\|u^0\|_{L^2(\mathcal{G})} = \rho$.
The implementation is described in Algorithm~\ref{alg:CNGF}, where we have chosen a stopping criterion corresponding  to the stagnation of the sequence of vectors $[u^n]$ in the $\ell^2$-norm. The gradient step requires to solve a linear system whose matrix is
\begin{equation*}
[[M_n]] = [[\mathrm{Id}]] + \delta t \left([[H]] - [[|u^{n}|^{p-1}]]\right).
\end{equation*}
Here, the matrix $[[|u^{n}|^{p-1}]]$ is a diagonal matrix constructed from the vector $[|u^{n}|^{p-1}]$.

\begin{algorithm}
\caption{CNGF algorithm \label{alg:CNGF}}
\begin{algorithmic} 
\REQUIRE $[u_0] \in \ell^2(\mathcal{G})$ with $\|[u_0]\|_{\ell^2} = \rho$,
$\varepsilon>0$, $\texttt{Stop\_Crit} = \texttt{True}$, $n = 0$ and $\texttt{Iter\_max}=1000$
\WHILE{$\texttt{Stop\_Crit}$ \textbf{and} $n\leq \texttt{Iter\_max}$}
\STATE{$\mathbf{Solve}\quad\left([[\mathrm{Id}]] + \delta t [[H]] - \delta t [[|u^{n}|^{p-1}]]\right)[u_*^{n+1}] = [u^n]$ }
\STATE{$[u^{n+1}] \gets \rho [u^{n+1}_*]/\|[u^{n+1}_*]\|_{\ell^2}$}
\STATE{$\texttt{Stop\_Crit} \gets \|[u^{n+1}] - [u^n]\|_{\ell^2} / \|[u^n]\|_{\ell^2}>\varepsilon$}
\STATE{$n\gets n+1$}
\ENDWHILE
\end{algorithmic}
\end{algorithm}

We now proceed to translate Algorithm~\ref{alg:CNGF} (with $p = 3$) into a
Python script  using the Grafidi library. First, we need to construct a quantum
graph. We choose to use the same graph as in Listing~\ref{lst:eigenelements}. Our code can be seen in Listing~\ref{lst:CNFG} (we avoid repetition in the listings, and consider that Listing~\ref{lst:eigenelements} is executed prior to Listing~\ref{lst:CNFG}).

\begin{listing}[htbp!]




  \begin{Verbatim}[commandchars=\\\{\},linenos=false,frame=single,framesep=2mm,fontsize=\footnotesize]
\PYG{n}{fun} \PYG{o}{=} \PYG{p}{\PYGZob{}\PYGZcb{}}
\PYG{n}{fun}\PYG{p}{[(}\PYG{l+s+s1}{\PYGZsq{}D\PYGZsq{}}\PYG{p}{,} \PYG{l+s+s1}{\PYGZsq{}B\PYGZsq{}}\PYG{p}{,} \PYG{l+s+s1}{\PYGZsq{}0\PYGZsq{}}\PYG{p}{)]}\PYG{o}{=}\PYG{k}{lambda} \PYG{n}{x}\PYG{p}{:} \PYG{n}{np}\PYG{o}{.}\PYG{n}{exp}\PYG{p}{(}\PYG{o}{\PYGZhy{}}\PYG{l+m+mf}{10e\PYGZhy{}2}\PYG{o}{*}\PYG{p}{(}\PYG{n}{x}\PYG{o}{\PYGZhy{}}\PYG{l+m+mi}{20}\PYG{p}{)}\PYG{o}{**}\PYG{l+m+mi}{2}\PYG{p}{)}
\PYG{n}{fun}\PYG{p}{[(}\PYG{l+s+s1}{\PYGZsq{}C\PYGZsq{}}\PYG{p}{,} \PYG{l+s+s1}{\PYGZsq{}A\PYGZsq{}}\PYG{p}{,} \PYG{l+s+s1}{\PYGZsq{}0\PYGZsq{}}\PYG{p}{)]}\PYG{o}{=}\PYG{k}{lambda} \PYG{n}{x}\PYG{p}{:} \PYG{n}{np}\PYG{o}{.}\PYG{n}{exp}\PYG{p}{(}\PYG{o}{\PYGZhy{}}\PYG{l+m+mf}{10e\PYGZhy{}2}\PYG{o}{*}\PYG{p}{(}\PYG{n}{x}\PYG{o}{\PYGZhy{}}\PYG{l+m+mi}{20}\PYG{p}{)}\PYG{o}{**}\PYG{l+m+mi}{2}\PYG{p}{)}
\PYG{n}{fun}\PYG{p}{[(}\PYG{l+s+s1}{\PYGZsq{}A\PYGZsq{}}\PYG{p}{,} \PYG{l+s+s1}{\PYGZsq{}B\PYGZsq{}}\PYG{p}{,} \PYG{l+s+s1}{\PYGZsq{}0\PYGZsq{}}\PYG{p}{)]}\PYG{o}{=}\PYG{k}{lambda} \PYG{n}{x}\PYG{p}{:} \PYG{l+m+mi}{1}\PYG{o}{\PYGZhy{}}\PYG{p}{(}\PYG{n}{x}\PYG{o}{\PYGZhy{}}\PYG{l+m+mi}{10}\PYG{p}{)}\PYG{o}{*}\PYG{n}{x}\PYG{o}{/}\PYG{l+m+mi}{50}
\PYG{n}{fun}\PYG{p}{[(}\PYG{l+s+s1}{\PYGZsq{}B\PYGZsq{}}\PYG{p}{,} \PYG{l+s+s1}{\PYGZsq{}A\PYGZsq{}}\PYG{p}{,} \PYG{l+s+s1}{\PYGZsq{}0\PYGZsq{}}\PYG{p}{)]}\PYG{o}{=}\PYG{k}{lambda} \PYG{n}{x}\PYG{p}{:} \PYG{l+m+mi}{1}\PYG{o}{+}\PYG{p}{(}\PYG{n}{x}\PYG{o}{\PYGZhy{}}\PYG{l+m+mi}{5}\PYG{p}{)}\PYG{o}{*}\PYG{n}{x}\PYG{o}{/}\PYG{l+m+mi}{20}
\PYG{n}{fun}\PYG{p}{[(}\PYG{l+s+s1}{\PYGZsq{}B\PYGZsq{}}\PYG{p}{,} \PYG{l+s+s1}{\PYGZsq{}A\PYGZsq{}}\PYG{p}{,} \PYG{l+s+s1}{\PYGZsq{}1\PYGZsq{}}\PYG{p}{)]}\PYG{o}{=}\PYG{k}{lambda} \PYG{n}{x}\PYG{p}{:} \PYG{l+m+mi}{1}\PYG{o}{+}\PYG{p}{(}\PYG{n}{x}\PYG{o}{\PYGZhy{}}\PYG{l+m+mi}{10}\PYG{p}{)}\PYG{o}{*}\PYG{n}{x}\PYG{o}{/}\PYG{l+m+mi}{30}

\PYG{n}{u}   \PYG{o}{=} \PYG{n}{WF}\PYG{p}{(}\PYG{n}{fun}\PYG{p}{,}\PYG{n}{g}\PYG{p}{)}
\PYG{n}{rho} \PYG{o}{=} \PYG{l+m+mi}{1}
\PYG{n}{u}   \PYG{o}{=} \PYG{n}{rho}\PYG{o}{*}\PYG{n}{u}\PYG{o}{/}\PYG{n}{WF}\PYG{o}{.}\PYG{n}{norm}\PYG{p}{(}\PYG{n}{u}\PYG{p}{,}\PYG{l+m+mi}{2}\PYG{p}{)}

\PYG{k}{def} \PYG{n+nf}{E}\PYG{p}{(}\PYG{n}{u}\PYG{p}{):}
    \PYG{k}{return} \PYG{o}{\PYGZhy{}}\PYG{l+m+mf}{0.5}\PYG{o}{*}\PYG{n}{WF}\PYG{o}{.}\PYG{n}{Lap}\PYG{p}{(}\PYG{n}{u}\PYG{p}{)}\PYG{o}{.}\PYG{n}{dot}\PYG{p}{(}\PYG{n}{u}\PYG{p}{)} \PYG{o}{\PYGZhy{}} \PYG{l+m+mf}{0.25}\PYG{o}{*}\PYG{n}{WF}\PYG{o}{.}\PYG{n}{norm}\PYG{p}{(}\PYG{n}{u}\PYG{p}{,}\PYG{l+m+mi}{4}\PYG{p}{)}\PYG{o}{**}\PYG{l+m+mi}{4}
\PYG{n}{En0} \PYG{o}{=} \PYG{n}{E}\PYG{p}{(}\PYG{n}{u}\PYG{p}{)}

\PYG{n}{delta\PYGZus{}t} \PYG{o}{=} \PYG{l+m+mf}{10e\PYGZhy{}1}
\PYG{n}{Epsilon} \PYG{o}{=} \PYG{l+m+mf}{10e\PYGZhy{}8}
\PYG{n}{M\PYGZus{}1} \PYG{o}{=} \PYG{n}{g}\PYG{o}{.}\PYG{n}{Id} \PYG{o}{\PYGZhy{}} \PYG{n}{delta\PYGZus{}t}\PYG{o}{*}\PYG{n}{g}\PYG{o}{.}\PYG{n}{Lap}
\PYG{k}{for} \PYG{n}{n} \PYG{o+ow}{in} \PYG{n+nb}{range}\PYG{p}{(}\PYG{l+m+mi}{1000}\PYG{p}{):}
    \PYG{n}{u\PYGZus{}old} \PYG{o}{=} \PYG{n}{u}
    \PYG{n}{M} \PYG{o}{=} \PYG{n}{M\PYGZus{}1} \PYG{o}{\PYGZhy{}} \PYG{n}{delta\PYGZus{}t}\PYG{o}{*}\PYG{n}{GR}\PYG{o}{.}\PYG{n}{Diag}\PYG{p}{(}\PYG{n}{g}\PYG{p}{,}\PYG{n+nb}{abs}\PYG{p}{(}\PYG{n}{u}\PYG{p}{)}\PYG{o}{**}\PYG{l+m+mi}{2}\PYG{p}{)}
    \PYG{n}{u} \PYG{o}{=} \PYG{n}{WF}\PYG{o}{.}\PYG{n}{Solve}\PYG{p}{(}\PYG{n}{M}\PYG{p}{,}\PYG{n}{u}\PYG{p}{)}
    \PYG{n}{u} \PYG{o}{=} \PYG{n}{rho}\PYG{o}{*}\PYG{n}{WF}\PYG{o}{.}\PYG{n}{abs}\PYG{p}{(}\PYG{n}{u}\PYG{p}{)}\PYG{o}{/}\PYG{n}{WF}\PYG{o}{.}\PYG{n}{norm}\PYG{p}{(}\PYG{n}{u}\PYG{p}{,}\PYG{l+m+mi}{2}\PYG{p}{)}
    \PYG{n}{En} \PYG{o}{=} \PYG{n}{E}\PYG{p}{(}\PYG{n}{u}\PYG{p}{)}
    \PYG{n+nb}{print}\PYG{p}{(}\PYG{l+s+sa}{f}\PYG{l+s+s2}{\PYGZdq{}Energy evolution: }\PYG{l+s+si}{\PYGZob{}}\PYG{n}{En}\PYG{o}{\PYGZhy{}}\PYG{n}{En0} \PYG{l+s+si}{:}\PYG{l+s+s2}{ 12.8e}\PYG{l+s+si}{\PYGZcb{}}\PYG{l+s+s2}{\PYGZdq{}}\PYG{p}{,}\PYG{n}{end}\PYG{o}{=}\PYG{l+s+s1}{\PYGZsq{}}\PYG{l+s+se}{\PYGZbs{}r}\PYG{l+s+s1}{\PYGZsq{}}\PYG{p}{)}
    \PYG{n}{En0} \PYG{o}{=} \PYG{n}{En}
    \PYG{n}{Stop\PYGZus{}crit} \PYG{o}{=} \PYG{n}{WF}\PYG{o}{.}\PYG{n}{norm}\PYG{p}{(}\PYG{n}{u}\PYG{o}{\PYGZhy{}}\PYG{n}{u\PYGZus{}old}\PYG{p}{,}\PYG{l+m+mi}{2}\PYG{p}{)}\PYG{o}{/}\PYG{n}{WF}\PYG{o}{.}\PYG{n}{norm}\PYG{p}{(}\PYG{n}{u\PYGZus{}old}\PYG{p}{,}\PYG{l+m+mi}{2}\PYG{p}{)}\PYG{o}{\PYGZlt{}}\PYG{n}{Epsilon}
    \PYG{k}{if} \PYG{n}{Stop\PYGZus{}crit}\PYG{p}{:}
        \PYG{k}{break}

\PYG{n}{\PYGZus{}}\PYG{o}{=}\PYG{n}{WF}\PYG{o}{.}\PYG{n}{draw}\PYG{p}{(}\PYG{n}{u}\PYG{p}{)}
\PYG{n+nb}{print}\PYG{p}{()}
\end{Verbatim}

\caption{Computation of a ground state using the CNGF method.}
\label{lst:CNFG}
\end{listing}






A few comments are in order. The initial data $u^0$ is set as a function that is quadratic on the edges connecting $A$ and $B$, increasing from $D$ to $B$ as
$\exp(-0.01(x-20)^2)$ (where $20$ is the length of $[DB]$) and increasing from
$C$ to $A$ as $\exp(-0.01(x-20)^2)$ (where $20$ is the length of $[CA]$). An
instance $u$ of {WFGraph} on \texttt{g} which corresponds to $u^0$ is
constructed accordingly.
The variable
$\rho = 1$
 corresponding to the $L^2$-norm is set
and the variable $u$ is normalized by using the function \texttt{norm} of
{WFGraph}. A function \texttt{E} is defined that corresponds to the energy and
we can see that we have used the \texttt{Lap} function of {WFGraph} to apply the
operator $[[H]]$ to $u$ as well as the function \texttt{dot} to compute the
scalar product. We set the variables $\delta t = 10^{-1}$ and $\varepsilon=
10^{-8}$. The part of the matrix $[[M_n]]$ that is independent of $n$ is built
in the variable \texttt{M\_1} which is the sum of $[[\mathrm{Id}]]$ (given by
the variable \texttt{g.Id} from Graph class) and $-\delta t[[H]]$ (where $-[[H]]$ is given by the
variable \texttt{g.Lap} from Graph class) and we also note that the matrix is
sparse. When entering the loop (with at most $1000$ iterations), we make a copy of $u$,
then construct the matrix $[[M_n]]$ by adding the diagonal matrix from the
nonlinearity (given through the function \texttt{GR.Diag} from {Graph}). The
linear system whose matrix is $[[M_n]]$ and right-hand-side $[[u^n]]$ is solved
thanks to the function \texttt{Solve} from {WFGraph}. Then, the variable $u$ is
normalized, the evolution of the energy is printed, the stopping criterion is
computed through the boolean variable \texttt{Stop\_Crit} and, finally, we
verify if the stopping criterion is attained (in which case we exit the loop and
draw $u$). 

\begin{figure}[htpb!]
\includegraphics[width=0.8\textwidth]{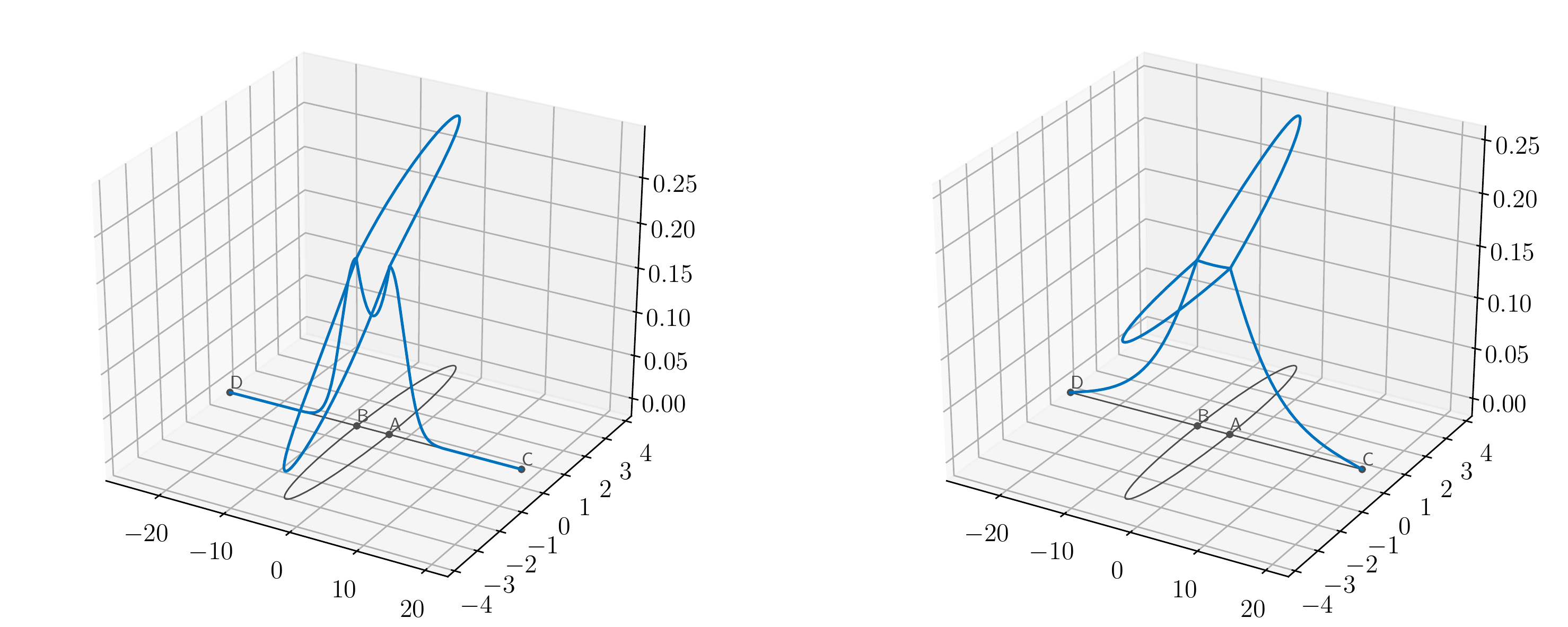}

\caption{Initial data (left) and ground state (right) obtained with the code from Listing~\ref{lst:CNFG}.}
  \label{fig:sol_CNFG}
\end{figure}

In the end, we obtain a ground state depicted in Figure~\ref{fig:sol_CNFG} which is computed in
$665$
iterations.

\subsubsection{The nonlinear conjugate gradient flow}

We now turn to the more sophisticated nonlinear conjugate gradient method. It is an extension of the conjugate gradient method that is used to solve linear systems. Here, we choose to use in the context of quantum graphs the method described for full spaces in~\cite[Algorithm 2]{antoine2017efficient}, which uses a preconditionner providing robustness. The method consists in the construction of a sequence $\{u^n\}_{n\geq0}$ (converging to the ground state), which is recursively defined by 
\begin{equation*}
\left\{\begin{array}{ll}
r^n = \mathcal{P}_{T,u^n}(Hu^n - |u^n|^{p-1}u^n)
\\ \beta^n = \max\left(0,\langle r^n - r^{n-1},Pr^{n}\rangle / \langle r^{n-1},Pr^{n-1} \rangle\right)
\\ d^n = - Pr^{n} + \beta_n p^{n-1}
\\ p^n = \mathcal{P}_{T,u^n} d^n
\\ \theta^n = \min_{\theta\in[-\pi,\pi]}E(\cos(\theta) u^n + \sin(\theta) \mathcal{P}_{S} p^n)
\\ u^{n+1} = \cos(\theta^n) u^n + \sin(\theta^n) \mathcal{P}_{S} p^n
\end{array}\right.
\end{equation*}
where $\mathcal{P}_{T,u}$ is the orthogonal projection on the tangent manifold of the sphere $\mathbb{S}_\rho = \{u\in L^2(\mathcal{G}): \|u\|_{L^2} = \rho\}$ at $u$ given by
\begin{equation*}
\mathcal{P}_{T,u}v = v - \frac{\langle v,u \rangle}{\|u\|_{L^2}^2} u,
\end{equation*}
$\mathcal{P}_S$ is the orthogonal projection on $\mathbb{S}_\rho$ given by
\begin{equation*}
\mathcal{P}_S v = \frac{\rho v}{\|v\|_{L^2}},
\end{equation*}
and $P = (1 + H)^{-1}$ is the preconditionner. The implementation of the method is described in Algorithm~\ref{alg:NCG}, where we have added a preliminary gradient descent step to initialize the iterative procedure.

\begin{algorithm}[!h]
\caption{Nonlinear Conjugate Gradient algorithm \label{alg:NCG}}
\begin{algorithmic} 
\REQUIRE $[u^{-1}] \in \ell^2(\mathcal{G})$ with $\|[u^{-1}]\|_{\ell^2} = \rho$,
$\varepsilon>0$, $\texttt{Stop\_Crit} = \texttt{True}$, $n = 0$ and $\texttt{Iter\_max}=500$
\STATE{$\lambda^{-1} \gets \langle [[H]] [u^{-1}] - [[|u^{-1}|^{p-1}]][u^{-1}],[u^{-1}]\rangle/\|[u^{-1}]\|_{\ell^2}$}
\STATE{$[r^{-1}] \gets [[H]] [u^{-1}] - [[|u^{-1}|^{p-1}]][u^{-1}] - \lambda^{-1} [u^{-1}]$}
\STATE{$\mathbf{Solve}\quad (\alpha[[\mathrm{Id}]]  + [[H]])[v^{-1}] = [r^{-1}]$}
\STATE{$[p^{-1}] \gets [v^{-1}] - \langle[v^{-1}],[u^{-1}]\rangle/\|[u^{-1}]\|_{\ell^2} [u^{-1}]$}
\STATE{$[\ell^{-1}] \gets \rho[p^{-1}]/\|[p^{-1}]\|_{\ell^2}$}
\STATE{$\mathbf{Minimize}\quad f(\theta^{-1}) = E(\cos(\theta^{-1})[u^{-1}] + \sin(\theta^{-1}) [\ell^{-1}]),\quad\theta^{-1}\in[-\pi,\pi]$}
\STATE{$[u^{0}] \gets \cos(\theta^{-1})[u^{-1}] + \sin(\theta^{-1}) [\ell^{-1}]$}

\WHILE{$\texttt{Stop\_Crit}$ \textbf{and} $n\leq \texttt{Iter\_max}$}
\STATE{$\lambda^n \gets \langle [[H]] [u^n] - [[|u^n|^{p-1}]][u^n],[u^n]\rangle/\|[u^n]\|_{\ell^2}$}
\STATE{$[r^n] \gets [[H]] [u^n] - [[|u^n|^{p-1}]][u^n] - \lambda^n [u^n]$}
\STATE{$\mathbf{Solve}\quad (\alpha[[\mathrm{Id}]]+[[H]])[v^n] = [r^n]$}
\STATE{$\beta^n \gets \max\left(0,\langle [r^n] - [r^{n-1}],[v^{n}]\rangle / \langle [r^{n-1}],[v^{n-1}] \rangle\right)$}
\STATE{$[d^n]\gets -[v^{n}] + \beta^n [p^{n-1}]$}
\STATE{$[p^n] \gets [d^n] - \langle [d^n],[u^n]\rangle/\|[u^n]\|_{\ell^2}[u^n]$}
\STATE{$[\ell^n] \gets \rho[p^n]/\|[p^n]\|_{\ell^2}$}
\STATE{$\mathbf{Minimize}\quad f(\theta^n) = E(\cos(\theta^n)[u^n] + \sin(\theta^n) [\ell^n]),\quad\theta^n\in[-\pi,\pi]$}
\STATE{$[u^{n+1}] \gets \cos(\theta^n)[u^n] + \sin(\theta^n) [\ell^n]$}
\STATE{$\texttt{Stop\_Crit} \gets \|[u^{n+1}] - [u^n]\|_{\ell^2} / \|[u^n]\|_{\ell^2}>\varepsilon$}
\STATE{$n\gets n+1$}
\ENDWHILE
\end{algorithmic}
\end{algorithm}

The corresponding code in Python, with the help of the Grafidi library, is given in Listings~\ref{lst:NCG}. We use the same quantum graph as in Section~\ref{sec:first-steps}, with in particular a $\delta$ condition with parameter $1$ at $O$. The initial function will be the function $u$ defined in~\eqref{eq:def-u}, normalized to verify the mass constraint. Listings~\ref{lst:RD}-\ref{lst:def_bc_N} are assumed to have been executed prior to Listings~\ref{lst:NCG}.




\begin{code}





    \begin{Verbatim}[commandchars=\\\{\},linenos=false,frame=single,framesep=2mm,fontsize=\footnotesize]
\PYG{k+kn}{import} \PYG{n+nn}{scipy.optimize} \PYG{k}{as} \PYG{n+nn}{sco}

\PYG{n}{rho} \PYG{o}{=} \PYG{l+m+mi}{2}
\PYG{n}{Epsilon} \PYG{o}{=} \PYG{l+m+mf}{10e\PYGZhy{}8}

\PYG{n}{u} \PYG{o}{=} \PYG{n}{rho}\PYG{o}{*}\PYG{n}{u}\PYG{o}{/}\PYG{n}{WF}\PYG{o}{.}\PYG{n}{norm}\PYG{p}{(}\PYG{n}{u}\PYG{p}{,}\PYG{l+m+mi}{2}\PYG{p}{)}

\PYG{k}{def} \PYG{n+nf}{E}\PYG{p}{(}\PYG{n}{u}\PYG{p}{):}
    \PYG{k}{return} \PYG{o}{\PYGZhy{}}\PYG{l+m+mf}{0.5}\PYG{o}{*}\PYG{n}{WF}\PYG{o}{.}\PYG{n}{Lap}\PYG{p}{(}\PYG{n}{u}\PYG{p}{)}\PYG{o}{.}\PYG{n}{dot}\PYG{p}{(}\PYG{n}{u}\PYG{p}{)} \PYG{o}{\PYGZhy{}} \PYG{l+m+mf}{0.25}\PYG{o}{*}\PYG{n}{WF}\PYG{o}{.}\PYG{n}{norm}\PYG{p}{(}\PYG{n}{u}\PYG{p}{,}\PYG{l+m+mi}{4}\PYG{p}{)}\PYG{o}{**}\PYG{l+m+mi}{4}
\PYG{k}{def} \PYG{n+nf}{P\PYGZus{}S}\PYG{p}{(}\PYG{n}{u}\PYG{p}{):}
    \PYG{k}{return} \PYG{n}{rho}\PYG{o}{*}\PYG{n}{u}\PYG{o}{/}\PYG{n}{WF}\PYG{o}{.}\PYG{n}{norm}\PYG{p}{(}\PYG{n}{u}\PYG{p}{,}\PYG{l+m+mi}{2}\PYG{p}{)}
\PYG{k}{def} \PYG{n+nf}{P\PYGZus{}T}\PYG{p}{(}\PYG{n}{u}\PYG{p}{,}\PYG{n}{v}\PYG{p}{):}
    \PYG{k}{return} \PYG{n}{v} \PYG{o}{\PYGZhy{}} \PYG{n}{v}\PYG{o}{.}\PYG{n}{dot}\PYG{p}{(}\PYG{n}{u}\PYG{p}{)}\PYG{o}{/}\PYG{p}{(}\PYG{n}{WF}\PYG{o}{.}\PYG{n}{norm}\PYG{p}{(}\PYG{n}{u}\PYG{p}{,}\PYG{l+m+mi}{2}\PYG{p}{)}\PYG{o}{**}\PYG{l+m+mi}{2}\PYG{p}{)}\PYG{o}{*}\PYG{n}{u}
\PYG{k}{def} \PYG{n+nf}{GradE}\PYG{p}{(}\PYG{n}{u}\PYG{p}{):}
    \PYG{k}{return} \PYG{o}{\PYGZhy{}}\PYG{n}{WF}\PYG{o}{.}\PYG{n}{Lap}\PYG{p}{(}\PYG{n}{u}\PYG{p}{)}\PYG{o}{\PYGZhy{}}\PYG{n}{WF}\PYG{o}{.}\PYG{n}{abs}\PYG{p}{(}\PYG{n}{u}\PYG{p}{)}\PYG{o}{**}\PYG{l+m+mi}{2}\PYG{o}{*}\PYG{n}{u}
\PYG{k}{def} \PYG{n+nf}{Pr}\PYG{p}{(}\PYG{n}{u}\PYG{p}{):}
    \PYG{k}{return} \PYG{n}{WF}\PYG{o}{.}\PYG{n}{Solve}\PYG{p}{(}\PYG{l+m+mf}{0.5}\PYG{o}{*}\PYG{n}{g}\PYG{o}{.}\PYG{n}{Id}\PYG{o}{\PYGZhy{}}\PYG{n}{g}\PYG{o}{.}\PYG{n}{Lap}\PYG{p}{,}\PYG{n}{u}\PYG{p}{)}
\PYG{k}{def} \PYG{n+nf}{E\PYGZus{}proj}\PYG{p}{(}\PYG{n}{theta}\PYG{p}{,}\PYG{n}{u}\PYG{p}{,}\PYG{n}{v}\PYG{p}{):}
    \PYG{k}{return} \PYG{n}{E}\PYG{p}{(}\PYG{n}{np}\PYG{o}{.}\PYG{n}{cos}\PYG{p}{(}\PYG{n}{theta}\PYG{p}{)}\PYG{o}{*}\PYG{n}{u}\PYG{o}{+}\PYG{n}{np}\PYG{o}{.}\PYG{n}{sin}\PYG{p}{(}\PYG{n}{theta}\PYG{p}{)}\PYG{o}{*}\PYG{n}{v}\PYG{p}{)}
\PYG{k}{def} \PYG{n+nf}{argmin\PYGZus{}E}\PYG{p}{(}\PYG{n}{u}\PYG{p}{,}\PYG{n}{v}\PYG{p}{):}
    \PYG{n}{theta} \PYG{o}{=} \PYG{n}{sco}\PYG{o}{.}\PYG{n}{fminbound}\PYG{p}{(}\PYG{n}{E\PYGZus{}proj}\PYG{p}{,}\PYG{o}{\PYGZhy{}}\PYG{n}{np}\PYG{o}{.}\PYG{n}{pi}\PYG{p}{,}\PYG{n}{np}\PYG{o}{.}\PYG{n}{pi}\PYG{p}{,(}\PYG{n}{u}\PYG{p}{,}\PYG{n}{v}\PYG{p}{),}\PYG{n}{xtol} \PYG{o}{=} \PYG{l+m+mf}{1e\PYGZhy{}14}\PYG{p}{,}\PYG{n}{maxfun} \PYG{o}{=} \PYG{l+m+mi}{1000}\PYG{p}{)}
    \PYG{k}{return} \PYG{n}{np}\PYG{o}{.}\PYG{n}{cos}\PYG{p}{(}\PYG{n}{theta}\PYG{p}{)}\PYG{o}{*}\PYG{n}{u}\PYG{o}{+}\PYG{n}{np}\PYG{o}{.}\PYG{n}{sin}\PYG{p}{(}\PYG{n}{theta}\PYG{p}{)}\PYG{o}{*}\PYG{n}{v}
\PYG{n}{En} \PYG{o}{=} \PYG{n}{E}\PYG{p}{(}\PYG{n}{u}\PYG{p}{)}

\PYG{n}{rm1} \PYG{o}{=} \PYG{n}{P\PYGZus{}T}\PYG{p}{(}\PYG{n}{u}\PYG{p}{,}\PYG{o}{\PYGZhy{}}\PYG{n}{GradE}\PYG{p}{(}\PYG{n}{u}\PYG{p}{))}
\PYG{n}{vm1} \PYG{o}{=} \PYG{n}{Pr}\PYG{p}{(}\PYG{n}{rm1}\PYG{p}{)}
\PYG{n}{pnm1} \PYG{o}{=} \PYG{n}{P\PYGZus{}T}\PYG{p}{(}\PYG{n}{u}\PYG{p}{,}\PYG{n}{Pr}\PYG{p}{(}\PYG{n}{rm1}\PYG{p}{))}
\PYG{n}{lm1} \PYG{o}{=} \PYG{n}{P\PYGZus{}S}\PYG{p}{(}\PYG{n}{pnm1}\PYG{p}{)}
\PYG{n}{u} \PYG{o}{=} \PYG{n}{argmin\PYGZus{}E}\PYG{p}{(}\PYG{n}{u}\PYG{p}{,}\PYG{n}{lm1}\PYG{p}{)}

\PYG{k}{for} \PYG{n}{n} \PYG{o+ow}{in} \PYG{n+nb}{range}\PYG{p}{(}\PYG{l+m+mi}{500}\PYG{p}{):}
    \PYG{n}{r} \PYG{o}{=} \PYG{n}{P\PYGZus{}T}\PYG{p}{(}\PYG{n}{u}\PYG{p}{,}\PYG{o}{\PYGZhy{}}\PYG{n}{GradE}\PYG{p}{(}\PYG{n}{u}\PYG{p}{))}
    \PYG{n}{v} \PYG{o}{=} \PYG{n}{Pr}\PYG{p}{(}\PYG{n}{r}\PYG{p}{)}
    \PYG{n}{beta} \PYG{o}{=} \PYG{n+nb}{max}\PYG{p}{(}\PYG{l+m+mi}{0}\PYG{p}{,(}\PYG{n}{r}\PYG{o}{\PYGZhy{}}\PYG{n}{rm1}\PYG{p}{)}\PYG{o}{.}\PYG{n}{dot}\PYG{p}{(}\PYG{n}{v}\PYG{p}{)}\PYG{o}{/}\PYG{n}{rm1}\PYG{o}{.}\PYG{n}{dot}\PYG{p}{(}\PYG{n}{vm1}\PYG{p}{))}
    \PYG{n}{rm1} \PYG{o}{=} \PYG{n}{r}
    \PYG{n}{vm1} \PYG{o}{=} \PYG{n}{v}
    \PYG{n}{d} \PYG{o}{=} \PYG{o}{\PYGZhy{}}\PYG{n}{v} \PYG{o}{+} \PYG{n}{beta}\PYG{o}{*}\PYG{n}{pnm1}
    \PYG{n}{p} \PYG{o}{=} \PYG{n}{P\PYGZus{}T}\PYG{p}{(}\PYG{n}{u}\PYG{p}{,}\PYG{n}{d}\PYG{p}{)}
    \PYG{n}{pm1} \PYG{o}{=} \PYG{n}{p}
    \PYG{n}{l} \PYG{o}{=} \PYG{n}{P\PYGZus{}S}\PYG{p}{(}\PYG{n}{p}\PYG{p}{)}
    \PYG{n}{um1} \PYG{o}{=} \PYG{n}{u}
    \PYG{n}{u} \PYG{o}{=} \PYG{n}{argmin\PYGZus{}E}\PYG{p}{(}\PYG{n}{u}\PYG{p}{,}\PYG{n}{l}\PYG{p}{)}
    \PYG{n}{En0} \PYG{o}{=} \PYG{n}{En}
    \PYG{n}{En} \PYG{o}{=} \PYG{n}{E}\PYG{p}{(}\PYG{n}{u}\PYG{p}{)}
    \PYG{n+nb}{print}\PYG{p}{(}\PYG{l+s+sa}{f}\PYG{l+s+s2}{\PYGZdq{}Energy evolution: }\PYG{l+s+si}{\PYGZob{}}\PYG{n}{En}\PYG{o}{\PYGZhy{}}\PYG{n}{En0} \PYG{l+s+si}{:}\PYG{l+s+s2}{ 12.8e}\PYG{l+s+si}{\PYGZcb{}}\PYG{l+s+s2}{\PYGZdq{}}\PYG{p}{,}\PYG{n}{end}\PYG{o}{=}\PYG{l+s+s1}{\PYGZsq{}}\PYG{l+s+se}{\PYGZbs{}r}\PYG{l+s+s1}{\PYGZsq{}}\PYG{p}{)}
    \PYG{n}{Stop\PYGZus{}crit} \PYG{o}{=} \PYG{n}{WF}\PYG{o}{.}\PYG{n}{norm}\PYG{p}{(}\PYG{n}{u}\PYG{o}{\PYGZhy{}}\PYG{n}{um1}\PYG{p}{,}\PYG{l+m+mi}{2}\PYG{p}{)}\PYG{o}{/}\PYG{n}{WF}\PYG{o}{.}\PYG{n}{norm}\PYG{p}{(}\PYG{n}{um1}\PYG{p}{,}\PYG{l+m+mi}{2}\PYG{p}{)}\PYG{o}{\PYGZlt{}}\PYG{n}{Epsilon}
    \PYG{k}{if} \PYG{n}{Stop\PYGZus{}crit}\PYG{p}{:}
        \PYG{k}{break}
\PYG{n}{\PYGZus{}}\PYG{o}{=}\PYG{n}{WF}\PYG{o}{.}\PYG{n}{draw}\PYG{p}{(}\PYG{n}{u}\PYG{p}{)}
\PYG{n+nb}{print}\PYG{p}{()}
\end{Verbatim}

  
  \caption{Computation of a ground state using the nonlinear conjugate gradient method.}
\label{lst:NCG}
\end{code}







\vspace*{4mm}

As for Listing~\ref{lst:CNFG}, a few comments are in order. The initial data $u^0$ is set as a function that is decreasing from $O$ to $A$, $O$ to $B$ and $O$ to $C$ as $\exp(-x^2)$. A function \texttt{P\_S} is defined that corresponds to $\mathcal{P}_S$, another one \texttt{P\_T} corresponds to $\mathcal{P}_{T,u}$ and another one \texttt{GradE} corresponds to the gradient of the energy. Furthermore, a function \texttt{Pr} computes the application of the preconditionner $[[P]]$ to an instance of {WFGraph} and returns the result as an instance of {WFGraph}. The function \texttt{E\_proj} computes the energy $E(\cos(\theta) u + \sin(\theta) v)$ with variables $u,v$ as instance of {WFGraph} and $\theta$ a scalar. The function \texttt{argmin\_E} is defined to return $w = \cos(\theta) u + \sin(\theta) v$ where $\theta$ is the solution of the minimum of $f(\theta) = E(\cos(\theta) u + \sin(\theta) v)$ and, moreover, it uses the function \texttt{fminbound} from Scipy (the maximum of iterations is fixed to $500$ and the tolerance to $10^{-8}$).

\begin{figure}[htpb!]
\includegraphics[width=.31\textwidth]{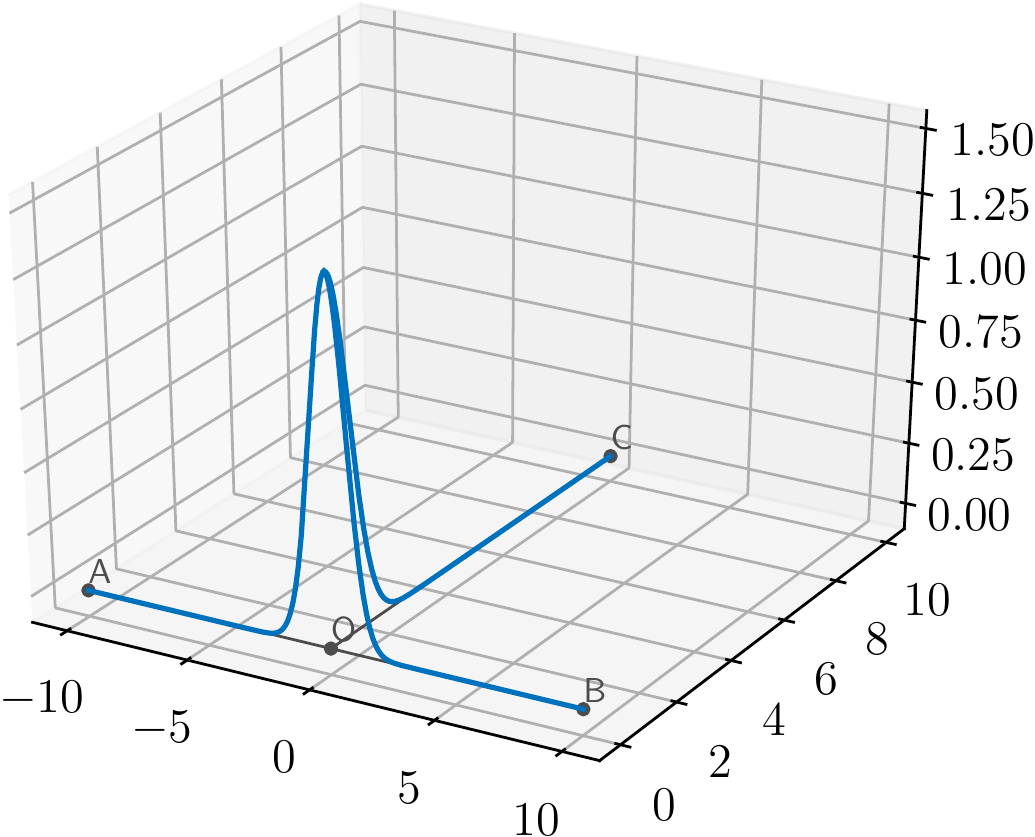}\hspace*{1.5cm}
\includegraphics[width=.31\textwidth]{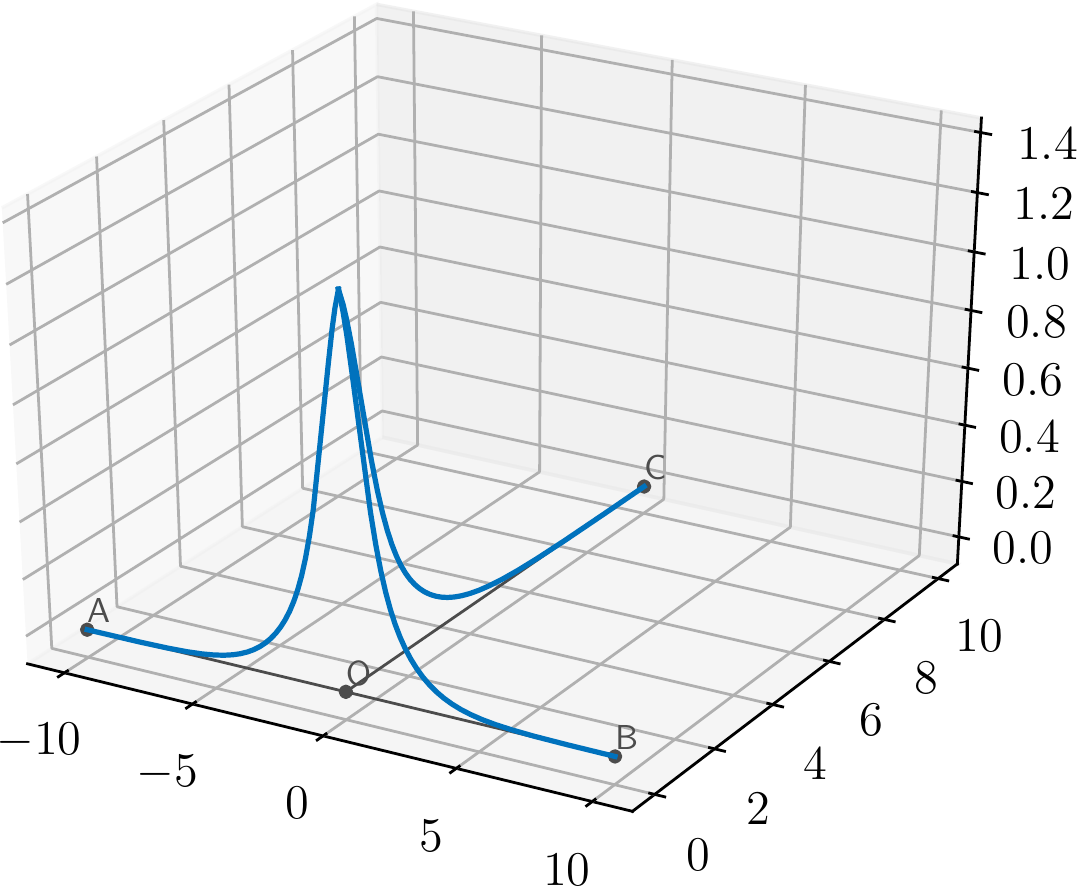}
\caption{Initial data (left) and ground state (right) obtained with the code from Listing~\ref{lst:NCG}.}
  \label{fig:sol_NCG}
\end{figure}

In the end, we obtain a ground state depicted in Figure~\ref{fig:sol_NCG} which is computed in $9$ iterations.

\subsection{Simulation of solutions of time-dependent nonlinear Schr\"odinger equations on graphs}
\label{sec:simulations}

In this section, we discuss the dynamical simulations of nonlinear Schr\"odinger equations on quantum graphs. To be more specific, we wish to simulate the solution on $\mathbb{R}^+\times\mathcal{G}$ of the following time-dependent equation
\begin{equation}\label{eq:dyn}
\left\{\begin{array}{ll}
i\partial_t \psi= H \psi - |\psi|^{2}\psi,
\\ \psi(t=0) = \psi_0\in L^2(\mathcal{G}).
\end{array}\right.
\end{equation}
We have chosen to present our results for the cubic power nonlinearity, but extension to other types of nonlinearity is straightforward.

\subsubsection{The Crank-Nicolson relaxation scheme}
One popular scheme to discretize~\eqref{eq:dyn} in time is the Crank-Nicolson scheme~\cite{delfour1981finite} which is of second order.  Since the equation is nonlinear,
the main drawback of the Crank-Nicolson scheme is the need to use a fixed-point method at each time step, which can be quite costly. To avoid this issue, we use the relaxation scheme proposed in~\cite{besse2004relaxation} which is semi-implicit and of second order. Let $\delta t>0$ be the time  step. The relaxation scheme applied to~\eqref{eq:dyn} is given by
\begin{equation}\label{eq:relax}
\left\{\begin{array}{ll}
\displaystyle \frac{\phi^{n+\frac12} + \phi^{n-\frac12}}{2} = - |\psi^{n}|^2\vspace{0.5em}
\\\displaystyle i\left(\frac{\psi^{n+1} - \psi^n}{\delta t}\right) = H\left(\frac{\psi^{n+1} + \psi^n}{2}\right) + \phi^{n+\frac12}\left(\frac{\psi^{n+1} + \psi^n}{2}\right),\quad\forall n\geq 0,
\vspace{0.5em}\\ \phi^{-\frac12} = -|\psi^{0}|^2 \quad\textrm{and}\quad \psi^0 = \psi_0\in L^2(\mathcal{G}),
\end{array}\right.
\end{equation}
where $\psi^{n}$ is an approximation of the solution $\psi$ of ~\eqref{eq:dyn} at time $n\delta t$. By introducing the intermediate variable $\psi^{n+\frac12} = (\psi^{n+1} + \psi^n)/2$, we deduce Algorithm~\ref{alg:relax}.

\begin{algorithm}
\caption{Relaxation scheme \label{alg:relax}}
\begin{algorithmic} 
\REQUIRE $[\psi^0] \in \ell^2(\mathcal{G})$, $\delta t>0$, $T>0$ and $N = \lceil T/\delta t \rceil$
\STATE{$[\phi^{-\frac12}] = - |[\psi^0]|^2 $}
\FOR{$n = 1,\ldots,N$}
\STATE{$[\phi^{n+\frac12}] = - 2 |[\psi^n]|^2 - [\phi^{n-\frac12}]$}
\STATE{$\mathbf{Solve}\quad\left([[\mathrm{Id}]] + i \delta t/2 [[H]] + i\delta t/2 [[\phi^{n+\frac12}]]\right)[\psi^{n+\frac12}] = [\psi^n]$ }
\STATE{$[\psi^{n+1}] \gets 2 [\psi^{n+\frac12}] - [\psi^n]$}
\ENDFOR
\end{algorithmic}
\end{algorithm}

\begin{figure}[htpb!]
\includegraphics[width=.31\textwidth]{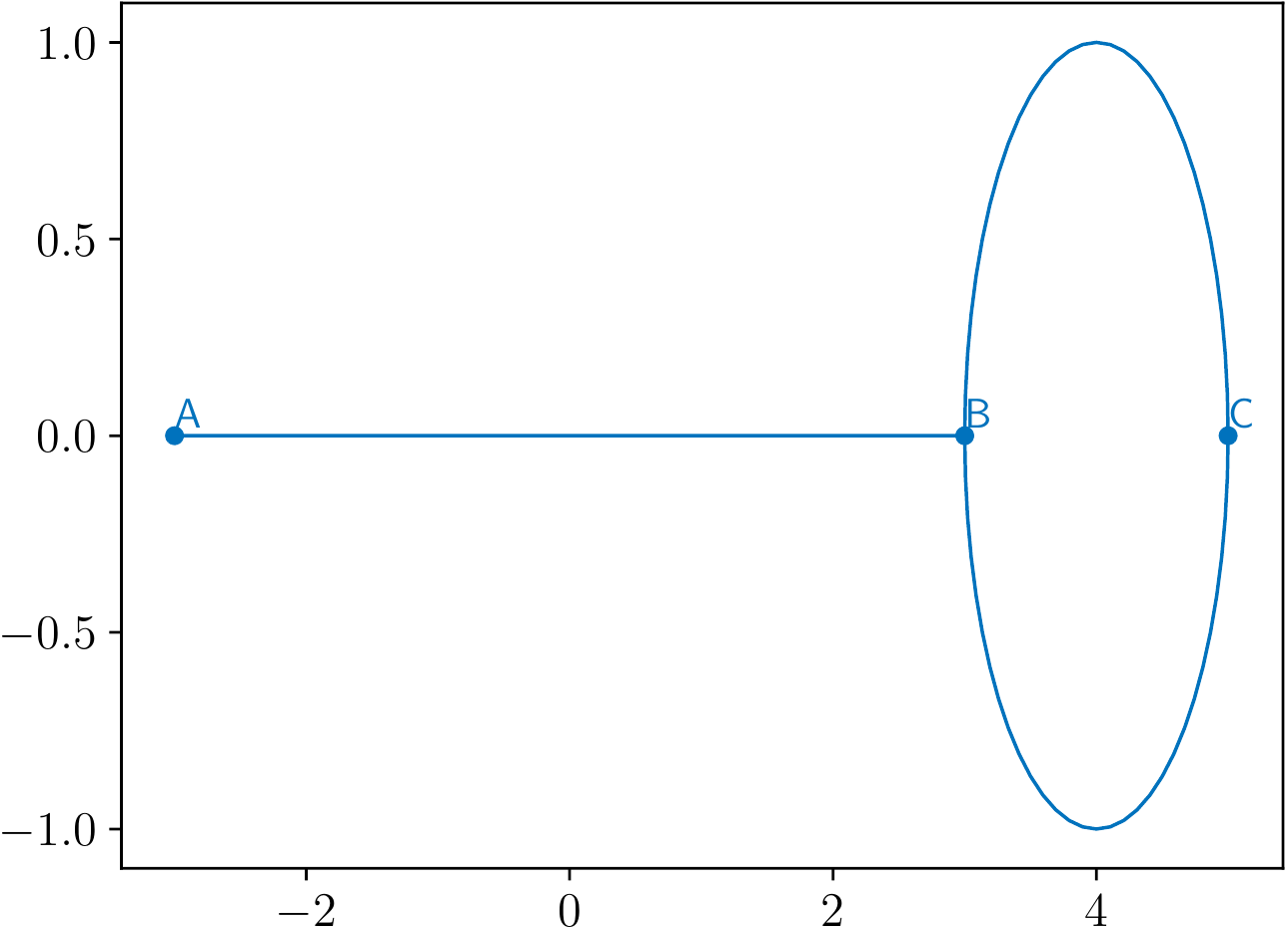}
\caption{Graph for the simulation with the relaxation scheme.\label{fig:signpost_graph}}
\end{figure}

We now wish to perform a simulation on the tadpole graph depicted in Figure~\ref{fig:signpost_graph} with the following lengths: $|AB| = 6$ and $|BC| = |CB| =\pi$. Observe here that we have to introduce an auxiliary vertex $C$ with Kirchhoff condition, and the loop is constructed as two half-loop edges connecting $B$ and $C$. This is of no consequence for the behavior of wave functions on the loop, as it was observed in~\cite[Remark 1.3.3]{BeKu13} that a vertex with Kirchhoff conditions with only two incident edges can always be removed. 

The boundary conditions that we take for the operator $H$ are Dirichlet at $A$ and Kirchhoff at $B$ and $C$. The initial data is taken as a bright soliton in the middle of the segment $[AB]$ with an initial velocity, \textit{i.e.}
\begin{equation*}
\psi^0(x) = \left\{\begin{array}{ll}
\displaystyle\frac{m}{2\sqrt{2}}\sech\left(\frac{m(x-3)}{4}\right)e^{i c x},&\text{for }x\in[AB],
\\ 0,& \text{for }x\in [BC]\cup[CB],
\end{array}\right.
\end{equation*}
with $m = 20$ and $c = 3$. Our simulation ends at time $T=1$ with a step
time of $\delta t = 10^{-3}$. We observe that the numerical scheme~\eqref{eq:relax} involves
complex valued functions $\psi^{n+1}$,  $\psi^n$ and $\phi^{n+1/2}$. We therefore have to explicitly
declare \texttt{WFGraph} instances with complex type. To this aim, we set
\verb+WF(fun,g,Dtype='complex')+. The argument \texttt{Dtype} is by default set to
\verb+'float'+. The choice \verb+Dtype='complex'+ allows to use \texttt{numpy.complex128} arrays
for linear algebra operations. This leads us to Listing~\ref{lst:relax} where we
implemented the relaxation scheme with the Grafidi library. This listing gives us the opportunity to discuss the outputs of \verb+draw+ function of WFGraph and their
usage. The return of \verb+draw+ is a three components tuple \verb+K,fig,ax+:
\begin{itemize}
\item \verb+fig+ is the matplotlib figure identifier where the plots are made,
\item \verb+ax+ is the matplotlib axes included in \verb+fig+,
\item \verb+K+ collection of elements actually drawn in \verb+ax+.
\end{itemize}
When we call \verb+draw+ with \verb+K+ as second argument, it automatically
updates the collection of elements in \verb+K+ into the figure \verb+fig+
without completely redrawing it, which is more efficient. In order to apply this
modification, we need to use both \verb+fig.canvas.draw()+ and \verb+plt.pause(0.01)+.

\begin{listing}[!htbp]




  \begin{Verbatim}[commandchars=\\\{\},linenos=false,frame=single,framesep=2mm,fontsize=\footnotesize]
\PYG{n}{g\PYGZus{}list} \PYG{o}{=} \PYG{p}{[}\PYG{l+s+s2}{\PYGZdq{}A B \PYGZob{}\PYGZsq{}Length\PYGZsq{}: 6\PYGZcb{}\PYGZdq{}}\PYG{p}{,} \PYG{l+s+s2}{\PYGZdq{}B C \PYGZob{}\PYGZsq{}Length\PYGZsq{}:3.14159\PYGZcb{}\PYGZdq{}}\PYG{p}{,} \PYG{l+s+s2}{\PYGZdq{}C B \PYGZob{}\PYGZsq{}Length\PYGZsq{}:3.14159\PYGZcb{}\PYGZdq{}}\PYG{p}{]}
\PYG{n}{g\PYGZus{}nx} \PYG{o}{=} \PYG{n}{nx}\PYG{o}{.}\PYG{n}{parse\PYGZus{}edgelist}\PYG{p}{(}\PYG{n}{g\PYGZus{}list}\PYG{p}{,}\PYG{n}{create\PYGZus{}using}\PYG{o}{=}\PYG{n}{nx}\PYG{o}{.}\PYG{n}{MultiDiGraph}\PYG{p}{())}
\PYG{n}{bc} \PYG{o}{=} \PYG{p}{\PYGZob{}}\PYG{l+s+s1}{\PYGZsq{}A\PYGZsq{}}\PYG{p}{:[}\PYG{l+s+s1}{\PYGZsq{}Dirichlet\PYGZsq{}}\PYG{p}{],} \PYG{l+s+s1}{\PYGZsq{}B\PYGZsq{}}\PYG{p}{:[}\PYG{l+s+s1}{\PYGZsq{}Kirchhoff\PYGZsq{}}\PYG{p}{],} \PYG{l+s+s1}{\PYGZsq{}C\PYGZsq{}}\PYG{p}{:[}\PYG{l+s+s1}{\PYGZsq{}Kirchhoff\PYGZsq{}}\PYG{p}{]\PYGZcb{}}
\PYG{n}{N}\PYG{o}{=}\PYG{l+m+mi}{3000}
\PYG{n}{g} \PYG{o}{=} \PYG{n}{GR}\PYG{p}{(}\PYG{n}{g\PYGZus{}nx}\PYG{p}{,}\PYG{n}{N}\PYG{p}{,}\PYG{n}{bc}\PYG{p}{)}
\PYG{n}{NewPos}\PYG{o}{=}\PYG{p}{\PYGZob{}}\PYG{l+s+s1}{\PYGZsq{}A\PYGZsq{}}\PYG{p}{:[}\PYG{o}{\PYGZhy{}}\PYG{l+m+mi}{3}\PYG{p}{,}\PYG{l+m+mi}{0}\PYG{p}{],}\PYG{l+s+s1}{\PYGZsq{}B\PYGZsq{}}\PYG{p}{:[}\PYG{l+m+mi}{3}\PYG{p}{,}\PYG{l+m+mi}{0}\PYG{p}{],}\PYG{l+s+s1}{\PYGZsq{}C\PYGZsq{}}\PYG{p}{:[}\PYG{l+m+mi}{5}\PYG{p}{,}\PYG{l+m+mi}{0}\PYG{p}{]\PYGZcb{}}
\PYG{n}{GR}\PYG{o}{.}\PYG{n}{Position}\PYG{p}{(}\PYG{n}{g}\PYG{p}{,}\PYG{n}{NewPos}\PYG{p}{)}

\PYG{n}{m} \PYG{o}{=} \PYG{l+m+mi}{20}
\PYG{n}{c} \PYG{o}{=} \PYG{l+m+mi}{3}
\PYG{n}{fun} \PYG{o}{=} \PYG{p}{\PYGZob{}\PYGZcb{}}
\PYG{n}{fun}\PYG{p}{[(}\PYG{l+s+s1}{\PYGZsq{}A\PYGZsq{}}\PYG{p}{,} \PYG{l+s+s1}{\PYGZsq{}B\PYGZsq{}}\PYG{p}{,} \PYG{l+s+s1}{\PYGZsq{}0\PYGZsq{}}\PYG{p}{)]}\PYG{o}{=}\PYG{k}{lambda} \PYG{n}{x}\PYG{p}{:} \PYG{n}{m}\PYG{o}{/}\PYG{l+m+mi}{2}\PYG{o}{/}\PYG{n}{np}\PYG{o}{.}\PYG{n}{sqrt}\PYG{p}{(}\PYG{l+m+mi}{2}\PYG{p}{)}\PYG{o}{/}\PYG{n}{np}\PYG{o}{.}\PYG{n}{cosh}\PYG{p}{(}\PYG{n}{m}\PYG{o}{*}\PYG{p}{(}\PYG{n}{x}\PYG{o}{\PYGZhy{}}\PYG{l+m+mi}{3}\PYG{p}{)}\PYG{o}{/}\PYG{l+m+mi}{4}\PYG{p}{)}\PYG{o}{*}\PYG{n}{np}\PYG{o}{.}\PYG{n}{exp}\PYG{p}{(}\PYG{l+m+mi}{1}\PYG{n}{j}\PYG{o}{*}\PYG{n}{c}\PYG{o}{*}\PYG{n}{x}\PYG{p}{)}
\PYG{n}{psi} \PYG{o}{=} \PYG{n}{WF}\PYG{p}{(}\PYG{n}{fun}\PYG{p}{,}\PYG{n}{g}\PYG{p}{,}\PYG{n}{Dtype}\PYG{o}{=}\PYG{l+s+s1}{\PYGZsq{}complex\PYGZsq{}}\PYG{p}{)}

\PYG{n}{K}\PYG{p}{,}\PYG{n}{fig}\PYG{p}{,}\PYG{n}{ax}\PYG{o}{=}\PYG{n}{WF}\PYG{o}{.}\PYG{n}{draw}\PYG{p}{(}\PYG{n}{WF}\PYG{o}{.}\PYG{n}{abs}\PYG{p}{(}\PYG{n}{psi}\PYG{p}{))}

\PYG{n}{T} \PYG{o}{=} \PYG{l+m+mi}{1}
\PYG{n}{delta\PYGZus{}t} \PYG{o}{=} \PYG{l+m+mf}{1e\PYGZhy{}3}
\PYG{n}{phi} \PYG{o}{=} \PYG{o}{\PYGZhy{}}\PYG{n}{WF}\PYG{o}{.}\PYG{n}{abs}\PYG{p}{(}\PYG{n}{psi}\PYG{p}{)}\PYG{o}{**}\PYG{l+m+mi}{2}
\PYG{n}{M\PYGZus{}1} \PYG{o}{=} \PYG{n}{g}\PYG{o}{.}\PYG{n}{Id} \PYG{o}{\PYGZhy{}} \PYG{l+m+mi}{1}\PYG{n}{j}\PYG{o}{*}\PYG{n}{delta\PYGZus{}t}\PYG{o}{/}\PYG{l+m+mi}{2}\PYG{o}{*}\PYG{n}{g}\PYG{o}{.}\PYG{n}{Lap}
\PYG{k}{for} \PYG{n}{n} \PYG{o+ow}{in} \PYG{n+nb}{range}\PYG{p}{(}\PYG{n+nb}{int}\PYG{p}{(}\PYG{n}{T}\PYG{o}{/}\PYG{n}{delta\PYGZus{}t}\PYG{p}{)}\PYG{o}{+}\PYG{l+m+mi}{1}\PYG{p}{):}
    \PYG{n}{phi} \PYG{o}{=} \PYG{o}{\PYGZhy{}}\PYG{l+m+mi}{2}\PYG{o}{*}\PYG{n}{WF}\PYG{o}{.}\PYG{n}{abs}\PYG{p}{(}\PYG{n}{psi}\PYG{p}{)}\PYG{o}{**}\PYG{l+m+mi}{2} \PYG{o}{\PYGZhy{}} \PYG{n}{phi}
    \PYG{n}{M} \PYG{o}{=} \PYG{n}{M\PYGZus{}1} \PYG{o}{+} \PYG{l+m+mi}{1}\PYG{n}{j}\PYG{o}{*}\PYG{n}{delta\PYGZus{}t}\PYG{o}{/}\PYG{l+m+mi}{2}\PYG{o}{*}\PYG{n}{GR}\PYG{o}{.}\PYG{n}{Diag}\PYG{p}{(}\PYG{n}{g}\PYG{p}{,}\PYG{n}{phi}\PYG{p}{)}
    \PYG{n}{varphi} \PYG{o}{=} \PYG{n}{WF}\PYG{o}{.}\PYG{n}{Solve}\PYG{p}{(}\PYG{n}{M}\PYG{p}{,}\PYG{n}{psi}\PYG{p}{)}
    \PYG{n}{psi} \PYG{o}{=} \PYG{l+m+mi}{2}\PYG{o}{*}\PYG{n}{varphi} \PYG{o}{\PYGZhy{}} \PYG{n}{psi}
    \PYG{k}{if} \PYG{n}{n}\PYG{o}{\PYGZpc{}}\PYG{l+m+mi}{100}\PYG{o}{==}\PYG{l+m+mi}{0}\PYG{p}{:}
        \PYG{n}{\PYGZus{}}\PYG{o}{=}\PYG{n}{WF}\PYG{o}{.}\PYG{n}{draw}\PYG{p}{(}\PYG{n}{WF}\PYG{o}{.}\PYG{n}{abs}\PYG{p}{(}\PYG{n}{psi}\PYG{p}{),}\PYG{n}{K}\PYG{p}{)}
        \PYG{n}{fig}\PYG{o}{.}\PYG{n}{canvas}\PYG{o}{.}\PYG{n}{draw}\PYG{p}{()}
        \PYG{n}{plt}\PYG{o}{.}\PYG{n}{pause}\PYG{p}{(}\PYG{l+m+mf}{0.01}\PYG{p}{)}

\PYG{n}{\PYGZus{}}\PYG{o}{=}\PYG{n}{WF}\PYG{o}{.}\PYG{n}{draw}\PYG{p}{(}\PYG{n}{WF}\PYG{o}{.}\PYG{n}{abs}\PYG{p}{(}\PYG{n}{psi}\PYG{p}{),}\PYG{n}{K}\PYG{p}{)}
\end{Verbatim}

\caption{Simulation of a soliton traveling in a quantum graph with a relaxation scheme.}
\label{lst:relax}
\end{listing}

\begin{figure}[!htbp]
  \centering
  \begin{tabular}{cccc}
    \includegraphics[height=.14\textheight]{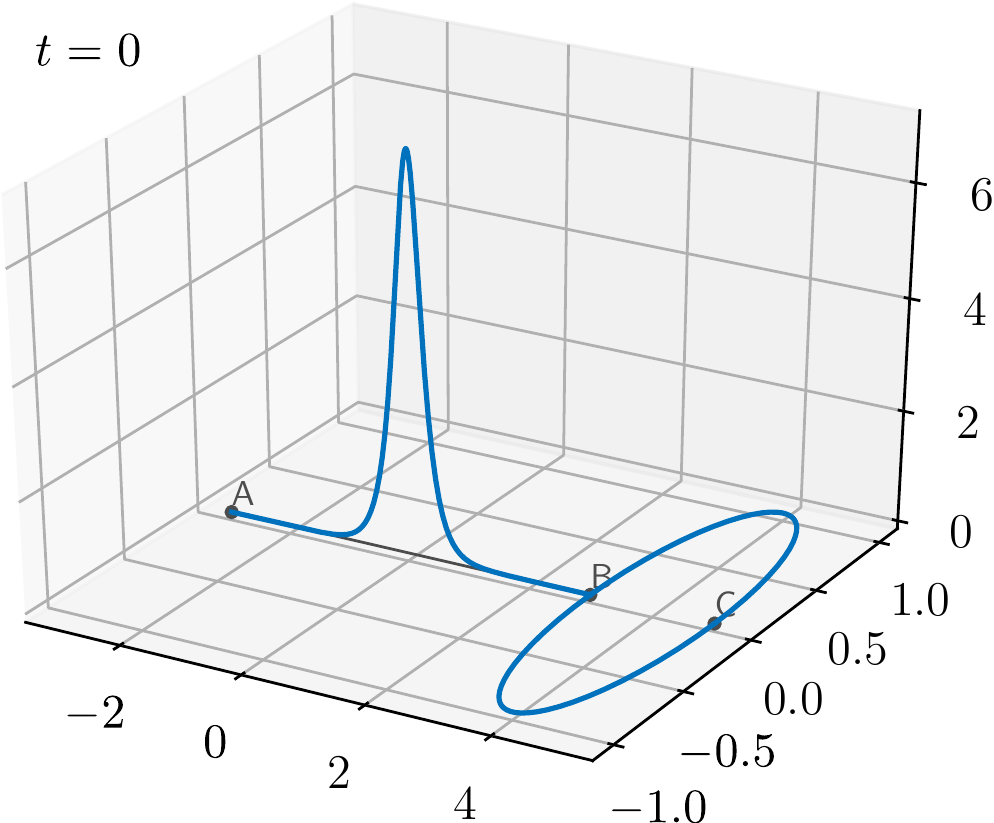}& \includegraphics[height=.14\textheight]{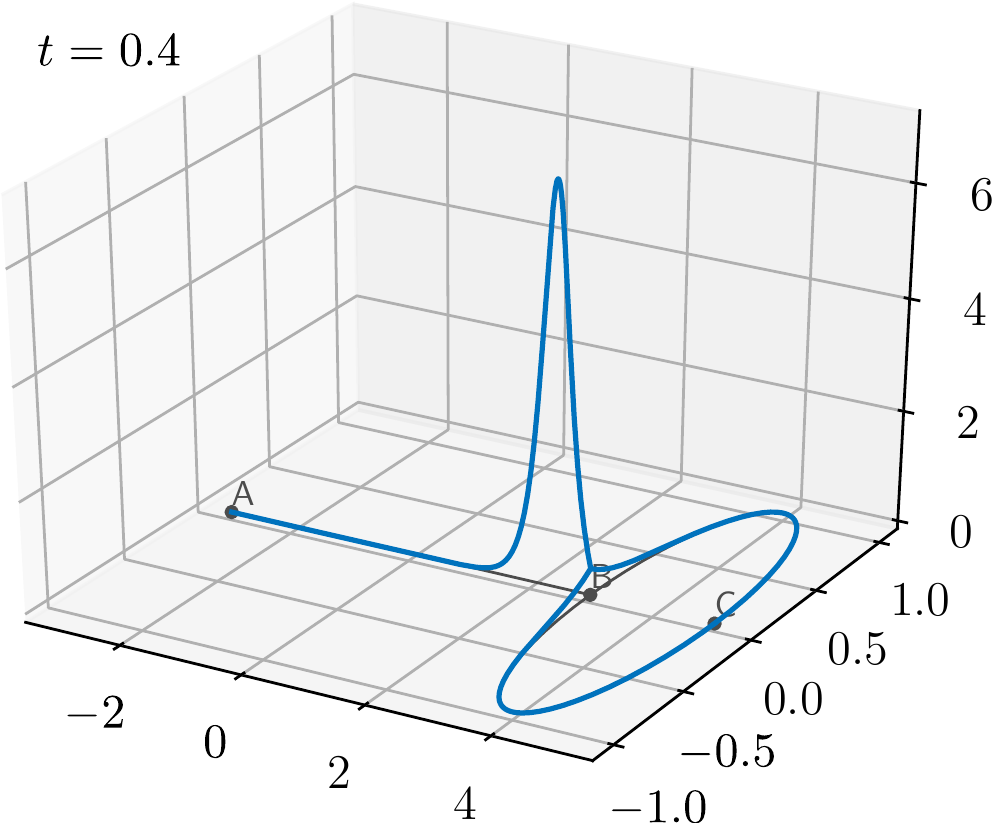}&
    \includegraphics[height=.14\textheight]{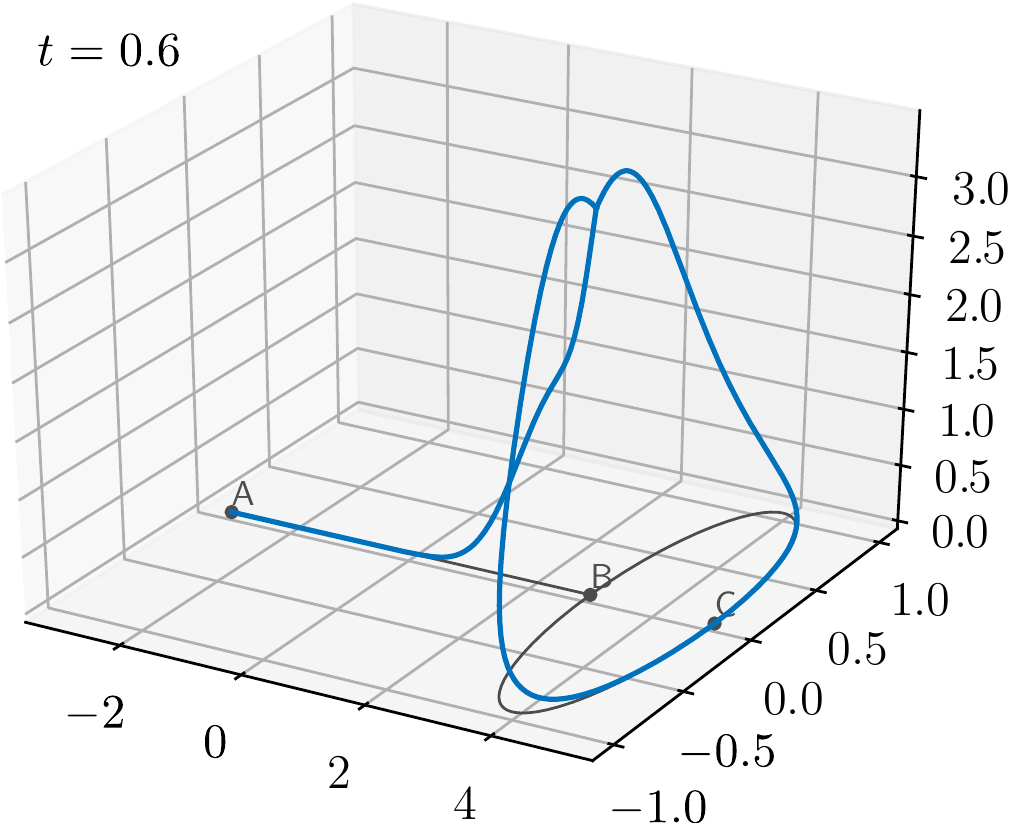}& \includegraphics[height=.14\textheight]{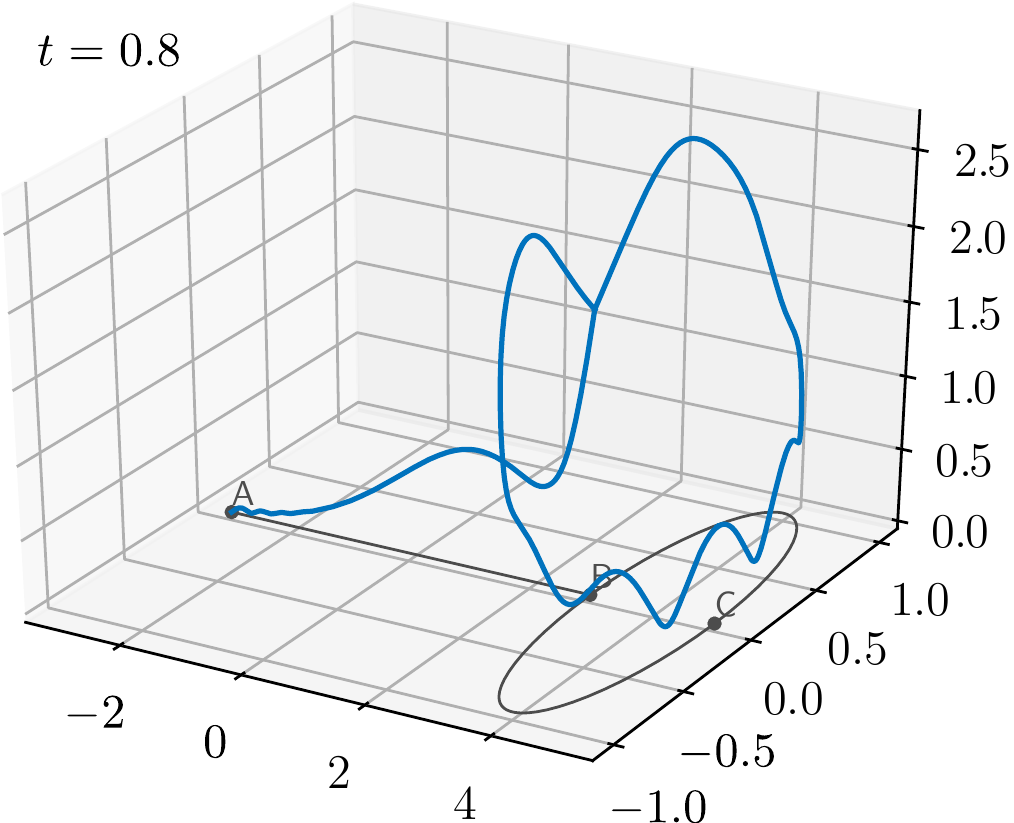}
  \end{tabular}
  \caption{Simulation of~\eqref{eq:dyn} at times $T = 0,0.4,0.6,0.8$ with the relaxation scheme.}
  \label{fig:relaxsim}
\end{figure}

The result of the simulation can be seen in Figure~\ref{fig:relaxsim} where the absolute value of $\psi$ at different times is given.


\subsubsection{The Strang splitting scheme}
Another popular approach for the simulation of nonlinear Schr\"odinger evolution is the so-called splitting method~\cite{weideman1986split}. As is well-known, the idea behind splitting methods is to ``split'' the full evolution equation into several (simpler) dynamical equations which are solved successively at each time step. In the case of~\eqref{eq:dyn}, we split the equation into a linear part and a nonlinear part. The equation corresponding to the linear part is
\begin{equation}\label{eq:SplitLin}
i\partial_t \psi = H \psi,
\end{equation}
and the equation associated to the nonlinear part is
\begin{equation*}
i\partial_t \psi = - |\psi|^2\psi.
\end{equation*}
This is motivated by the fact that the solution for the nonlinear part can be obtained explicitly. We use a Strang splitting scheme of second order~\cite{strang1968construction}.  For a given time step $\delta t>0$, we obtain the following method, for any $n\geq 0$,
\begin{equation*}
\left\{\begin{array}{ll}
\psi^{n+\frac13} = e^{i \delta t/2 |\psi^{n}|^2} \psi^{n},\vspace{0.5em}
\\ \frac{\psi^{n+\frac23}-\psi^{n+\frac13}}{\delta t} = H \left(\frac{\psi^{n+\frac23}+\psi^{n+\frac13}}{2}\right),\vspace{0.5em}
\\ \psi^{n+1} = e^{i \delta t/2 |\psi^{n+\frac23}|^2} \psi^{n+\frac23},
\end{array}\right.
\end{equation*}
where we have used a Crank-Nicolson scheme to discretize in time Equation~\eqref{eq:SplitLin}. Through the introduction of an intermediate variable $\psi^{n+\frac12} = (\psi^{n+\frac23}+\psi^{n+\frac13})/2$, we deduce Algorithm~\ref{alg:splitting}.

\begin{algorithm}
\caption{Strang splitting scheme \label{alg:splitting}}
\begin{algorithmic} 
\REQUIRE $[\psi^0] \in \ell^2(\mathcal{G})$, $\delta t>0$, $T>0$ and $N = \lceil T/\delta t \rceil$
\FOR{$n = 1,\ldots,N$}
\STATE{$[\psi^{n+\frac13}] \gets \texttt{exp}(i \delta t/2|[\psi^n]|^2)[\psi^n]$}
\STATE{$\mathbf{Solve}\quad\left([[\mathrm{Id}]] + i \delta t/2 [[H]]\right)[\psi^{n+\frac12}] = [\psi^{n+\frac13}]$ }
\STATE{$[\psi^{n+\frac23}] \gets 2 [\psi^{n+\frac12}] - [\psi^{n+\frac13}]$}
\STATE{$[\psi^{n+1}] \gets \texttt{exp}(i \delta t/2|[\psi^{n+\frac23}]|^2)[\psi^{n+\frac23}]$}
\ENDFOR
\end{algorithmic}
\end{algorithm}

\begin{figure}[htpb!]
\includegraphics[width=.31\textwidth]{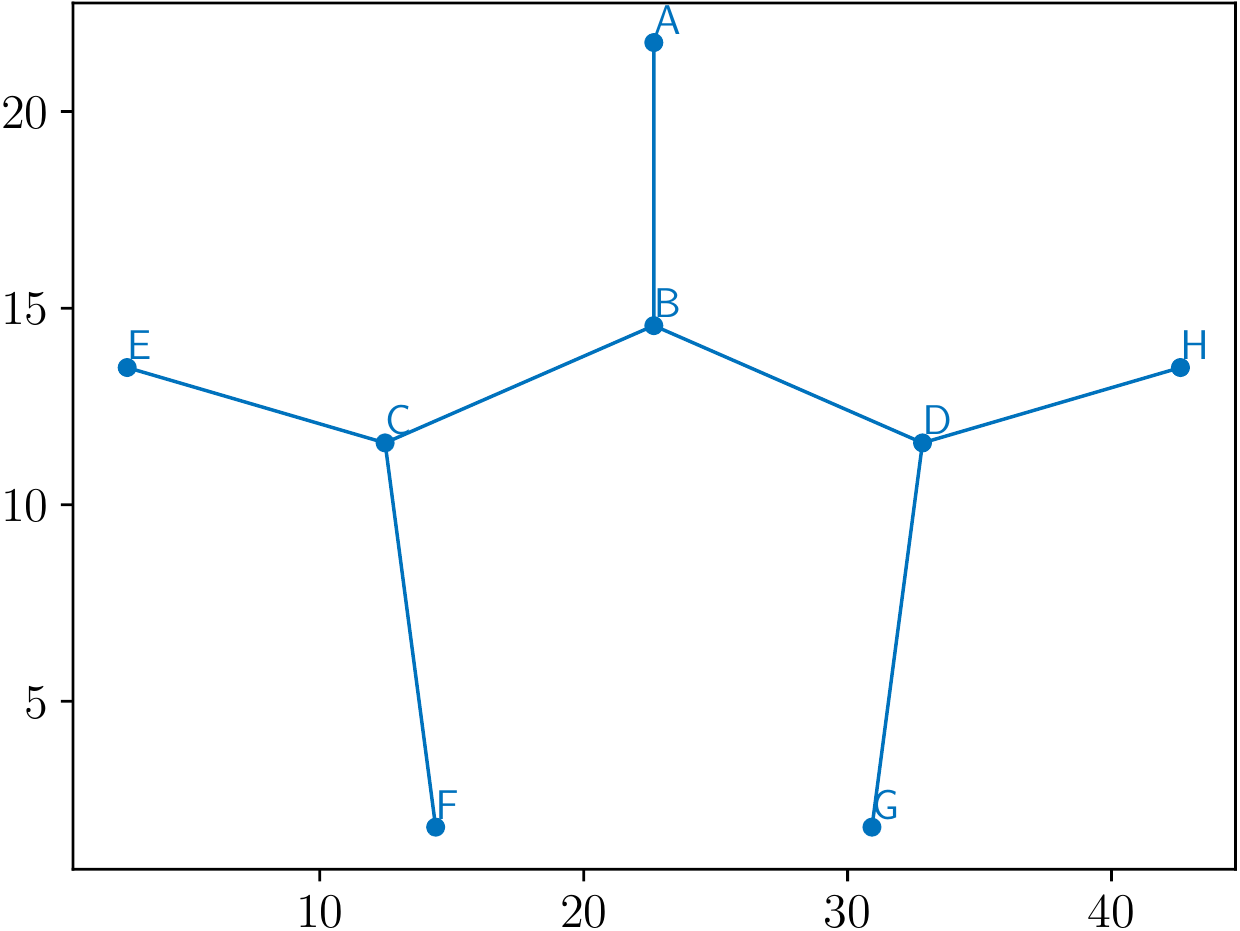}
\caption{Graph of the simulation for the splitting scheme.\label{fig:tree_graph}}
\end{figure}

We now wish to perform a simulation on the graph depicted in Figure~\ref{fig:tree_graph} with the following lengths: $|AB| = 6$, $|BC| = |BD| = 10.61$ and $|CE| = |CF| = |DG| = |DH| = 9.96$.

The boundary conditions that we would like for the operator $H$ are Kirchhoff at $B$, $C$ and $D$ and Dirichlet for all the others. The initial data is a bright soliton in the middle of segment $[AB]$ with an initial velocity $c = 3$ and  $m = 15$. Our simulation ends at time $T=2$ with a step time of $\delta t = 10^{-3}$. This leads us to Listing~\ref{lst:splitting} where we implemented the Strang splitting scheme with the Grafidi library.

\begin{listing}[!htbp]






  \begin{Verbatim}[commandchars=\\\{\},linenos=false,frame=single,framesep=2mm,fontsize=\footnotesize]
\PYG{n}{g\PYGZus{}list}\PYG{o}{=}\PYG{p}{[}\PYG{l+s+s2}{\PYGZdq{}A B \PYGZob{}\PYGZsq{}Length\PYGZsq{}:7.20\PYGZcb{}\PYGZdq{}}\PYG{p}{,} \PYG{l+s+s2}{\PYGZdq{}B C \PYGZob{}\PYGZsq{}Length\PYGZsq{}:10.61\PYGZcb{}\PYGZdq{}}\PYG{p}{,} \PYG{l+s+s2}{\PYGZdq{}B D \PYGZob{}\PYGZsq{}Length\PYGZsq{}:10.61\PYGZcb{}\PYGZdq{}}\PYG{p}{,}\PYGZbs{}
    \PYG{l+s+s2}{\PYGZdq{}C E \PYGZob{}\PYGZsq{}Length\PYGZsq{}:9.96\PYGZcb{}\PYGZdq{}}\PYG{p}{,} \PYG{l+s+s2}{\PYGZdq{}C F \PYGZob{}\PYGZsq{}Length\PYGZsq{}:9.96\PYGZcb{}\PYGZdq{}}\PYG{p}{,} \PYG{l+s+s2}{\PYGZdq{}D G \PYGZob{}\PYGZsq{}Length\PYGZsq{}:9.96\PYGZcb{}\PYGZdq{}}\PYG{p}{,} \PYGZbs{}
    \PYG{l+s+s2}{\PYGZdq{}D H \PYGZob{}\PYGZsq{}Length\PYGZsq{}:9.96\PYGZcb{}\PYGZdq{}}\PYG{p}{]}
\PYG{n}{g\PYGZus{}nx} \PYG{o}{=} \PYG{n}{nx}\PYG{o}{.}\PYG{n}{parse\PYGZus{}edgelist}\PYG{p}{(}\PYG{n}{g\PYGZus{}list}\PYG{p}{,}\PYG{n}{create\PYGZus{}using}\PYG{o}{=}\PYG{n}{nx}\PYG{o}{.}\PYG{n}{MultiDiGraph}\PYG{p}{())}
\PYG{n}{bc} \PYG{o}{=} \PYG{p}{\PYGZob{}}\PYG{l+s+s1}{\PYGZsq{}A\PYGZsq{}}\PYG{p}{:[}\PYG{l+s+s1}{\PYGZsq{}Dirichlet\PYGZsq{}}\PYG{p}{],}\PYG{l+s+s1}{\PYGZsq{}B\PYGZsq{}}\PYG{p}{:[}\PYG{l+s+s1}{\PYGZsq{}Kirchhoff\PYGZsq{}}\PYG{p}{],}\PYG{l+s+s1}{\PYGZsq{}C\PYGZsq{}}\PYG{p}{:[}\PYG{l+s+s1}{\PYGZsq{}Kirchhoff\PYGZsq{}}\PYG{p}{],}\PYGZbs{}
    \PYG{l+s+s1}{\PYGZsq{}D\PYGZsq{}}\PYG{p}{:[}\PYG{l+s+s1}{\PYGZsq{}Kirchhoff\PYGZsq{}}\PYG{p}{],}\PYG{l+s+s1}{\PYGZsq{}E\PYGZsq{}}\PYG{p}{:[}\PYG{l+s+s1}{\PYGZsq{}Dirichlet\PYGZsq{}}\PYG{p}{],}\PYG{l+s+s1}{\PYGZsq{}F\PYGZsq{}}\PYG{p}{:[}\PYG{l+s+s1}{\PYGZsq{}Dirichlet\PYGZsq{}}\PYG{p}{],}\PYGZbs{}
    \PYG{l+s+s1}{\PYGZsq{}G\PYGZsq{}}\PYG{p}{:[}\PYG{l+s+s1}{\PYGZsq{}Dirichlet\PYGZsq{}}\PYG{p}{],}\PYG{l+s+s1}{\PYGZsq{}H\PYGZsq{}}\PYG{p}{:[}\PYG{l+s+s1}{\PYGZsq{}Dirichlet\PYGZsq{}}\PYG{p}{]\PYGZcb{}}
\PYG{n}{N}\PYG{o}{=}\PYG{l+m+mi}{3000}
\PYG{n}{g} \PYG{o}{=} \PYG{n}{GR}\PYG{p}{(}\PYG{n}{g\PYGZus{}nx}\PYG{p}{,}\PYG{n}{N}\PYG{p}{,}\PYG{n}{bc}\PYG{p}{)}

\PYG{n}{NewPos} \PYG{o}{=} \PYG{p}{\PYGZob{}}
    \PYG{l+s+s1}{\PYGZsq{}A\PYGZsq{}}\PYG{p}{:} \PYG{p}{[}\PYG{l+m+mf}{22.656}\PYG{p}{,} \PYG{l+m+mf}{21.756}\PYG{p}{],} \PYG{l+s+s1}{\PYGZsq{}B\PYGZsq{}}\PYG{p}{:} \PYG{p}{[}\PYG{l+m+mf}{22.656}\PYG{p}{,} \PYG{l+m+mf}{14.556}\PYG{p}{],} \PYG{l+s+s1}{\PYGZsq{}C\PYGZsq{}}\PYG{p}{:} \PYG{p}{[}\PYG{l+m+mf}{12.473}\PYG{p}{,} \PYG{l+m+mf}{11.573}\PYG{p}{],}\PYGZbs{}
    \PYG{l+s+s1}{\PYGZsq{}D\PYGZsq{}}\PYG{p}{:} \PYG{p}{[}\PYG{l+m+mf}{32.838}\PYG{p}{,} \PYG{l+m+mf}{11.573}\PYG{p}{],} \PYG{l+s+s1}{\PYGZsq{}E\PYGZsq{}}\PYG{p}{:} \PYG{p}{[}\PYG{l+m+mf}{2.7}\PYG{p}{,} \PYG{l+m+mf}{13.49}\PYG{p}{],}     \PYG{l+s+s1}{\PYGZsq{}F\PYGZsq{}}\PYG{p}{:} \PYG{p}{[}\PYG{l+m+mf}{14.39}\PYG{p}{,} \PYG{l+m+mf}{1.8}\PYG{p}{],}\PYGZbs{}
    \PYG{l+s+s1}{\PYGZsq{}G\PYGZsq{}}\PYG{p}{:} \PYG{p}{[}\PYG{l+m+mf}{30.922}\PYG{p}{,} \PYG{l+m+mf}{1.8}\PYG{p}{],}    \PYG{l+s+s1}{\PYGZsq{}H\PYGZsq{}}\PYG{p}{:} \PYG{p}{[}\PYG{l+m+mf}{42.612}\PYG{p}{,} \PYG{l+m+mf}{13.49}\PYG{p}{]\PYGZcb{}}
\PYG{n}{GR}\PYG{o}{.}\PYG{n}{Position}\PYG{p}{(}\PYG{n}{g}\PYG{p}{,}\PYG{n}{NewPos}\PYG{p}{)}

\PYG{n}{m} \PYG{o}{=} \PYG{l+m+mi}{15}
\PYG{n}{c} \PYG{o}{=} \PYG{l+m+mi}{3}
\PYG{n}{x0} \PYG{o}{=} \PYG{l+m+mf}{7.2}\PYG{o}{/}\PYG{l+m+mi}{2}
\PYG{n}{fun} \PYG{o}{=} \PYG{p}{\PYGZob{}\PYGZcb{}}
\PYG{n}{fun}\PYG{p}{[(}\PYG{l+s+s1}{\PYGZsq{}A\PYGZsq{}}\PYG{p}{,} \PYG{l+s+s1}{\PYGZsq{}B\PYGZsq{}}\PYG{p}{,} \PYG{l+s+s1}{\PYGZsq{}0\PYGZsq{}}\PYG{p}{)]}\PYG{o}{=}\PYG{k}{lambda} \PYG{n}{x}\PYG{p}{:} \PYG{n}{m}\PYG{o}{/}\PYG{l+m+mi}{2}\PYG{o}{/}\PYG{n}{np}\PYG{o}{.}\PYG{n}{sqrt}\PYG{p}{(}\PYG{l+m+mi}{2}\PYG{p}{)}\PYG{o}{/}\PYG{n}{np}\PYG{o}{.}\PYG{n}{cosh}\PYG{p}{(}\PYG{n}{m}\PYG{o}{*}\PYG{p}{(}\PYG{n}{x}\PYG{o}{\PYGZhy{}}\PYG{n}{x0}\PYG{p}{)}\PYG{o}{/}\PYG{l+m+mi}{4}\PYG{p}{)}\PYG{o}{*}\PYG{n}{np}\PYG{o}{.}\PYG{n}{exp}\PYG{p}{(}\PYG{l+m+mi}{1}\PYG{n}{j}\PYG{o}{*}\PYG{n}{c}\PYG{o}{*}\PYG{n}{x}\PYG{p}{)}
\PYG{n}{psi} \PYG{o}{=} \PYG{n}{WF}\PYG{p}{(}\PYG{n}{fun}\PYG{p}{,}\PYG{n}{g}\PYG{p}{,}\PYG{n}{Dtype}\PYG{o}{=}\PYG{l+s+s1}{\PYGZsq{}complex\PYGZsq{}}\PYG{p}{)}

\PYG{n}{K}\PYG{p}{,}\PYG{n}{fig}\PYG{p}{,}\PYG{n}{ax}\PYG{o}{=}\PYG{n}{WF}\PYG{o}{.}\PYG{n}{draw}\PYG{p}{(}\PYG{n}{WF}\PYG{o}{.}\PYG{n}{abs}\PYG{p}{(}\PYG{n}{psi}\PYG{p}{))}

\PYG{n}{T} \PYG{o}{=} \PYG{l+m+mi}{2}
\PYG{n}{delta\PYGZus{}t} \PYG{o}{=} \PYG{l+m+mf}{1e\PYGZhy{}3}
\PYG{n}{M} \PYG{o}{=} \PYG{n}{g}\PYG{o}{.}\PYG{n}{Id} \PYG{o}{\PYGZhy{}} \PYG{l+m+mi}{1}\PYG{n}{j}\PYG{o}{*}\PYG{n}{delta\PYGZus{}t}\PYG{o}{*}\PYG{n}{g}\PYG{o}{.}\PYG{n}{Lap}\PYG{o}{/}\PYG{l+m+mi}{2}

\PYG{k}{for} \PYG{n}{n} \PYG{o+ow}{in} \PYG{n+nb}{range}\PYG{p}{(}\PYG{n+nb}{int}\PYG{p}{(}\PYG{n}{T}\PYG{o}{/}\PYG{n}{delta\PYGZus{}t}\PYG{p}{)):}
    \PYG{n}{psi} \PYG{o}{=} \PYG{n}{psi}\PYG{o}{*}\PYG{n}{WF}\PYG{o}{.}\PYG{n}{exp}\PYG{p}{(}\PYG{l+m+mi}{1}\PYG{n}{j}\PYG{o}{*}\PYG{n}{delta\PYGZus{}t}\PYG{o}{/}\PYG{l+m+mi}{2}\PYG{o}{*}\PYG{n}{WF}\PYG{o}{.}\PYG{n}{abs}\PYG{p}{(}\PYG{n}{psi}\PYG{p}{)}\PYG{o}{**}\PYG{l+m+mi}{2}\PYG{p}{)}
    \PYG{n}{varphi} \PYG{o}{=} \PYG{n}{WF}\PYG{o}{.}\PYG{n}{Solve}\PYG{p}{(}\PYG{n}{M}\PYG{p}{,}\PYG{n}{psi}\PYG{p}{)}
    \PYG{n}{psi} \PYG{o}{=} \PYG{l+m+mi}{2}\PYG{o}{*}\PYG{n}{varphi} \PYG{o}{\PYGZhy{}} \PYG{n}{psi}
    \PYG{n}{psi} \PYG{o}{=} \PYG{n}{psi}\PYG{o}{*}\PYG{n}{WF}\PYG{o}{.}\PYG{n}{exp}\PYG{p}{(}\PYG{l+m+mi}{1}\PYG{n}{j}\PYG{o}{*}\PYG{n}{delta\PYGZus{}t}\PYG{o}{/}\PYG{l+m+mi}{2}\PYG{o}{*}\PYG{n}{WF}\PYG{o}{.}\PYG{n}{abs}\PYG{p}{(}\PYG{n}{psi}\PYG{p}{)}\PYG{o}{**}\PYG{l+m+mi}{2}\PYG{p}{)}
    \PYG{k}{if} \PYG{n}{n}\PYG{o}{\PYGZpc{}}\PYG{l+m+mi}{100}\PYG{o}{==}\PYG{l+m+mi}{0}\PYG{p}{:}
        \PYG{n}{\PYGZus{}}\PYG{o}{=}\PYG{n}{WF}\PYG{o}{.}\PYG{n}{draw}\PYG{p}{(}\PYG{n}{WF}\PYG{o}{.}\PYG{n}{abs}\PYG{p}{(}\PYG{n}{psi}\PYG{p}{),}\PYG{n}{K}\PYG{p}{)}
        \PYG{n}{fig}\PYG{o}{.}\PYG{n}{canvas}\PYG{o}{.}\PYG{n}{draw}\PYG{p}{()}
        \PYG{n}{plt}\PYG{o}{.}\PYG{n}{pause}\PYG{p}{(}\PYG{l+m+mf}{0.01}\PYG{p}{)}

\PYG{n}{\PYGZus{}}\PYG{o}{=}\PYG{n}{WF}\PYG{o}{.}\PYG{n}{draw}\PYG{p}{(}\PYG{n}{WF}\PYG{o}{.}\PYG{n}{abs}\PYG{p}{(}\PYG{n}{psi}\PYG{p}{),}\PYG{n}{K}\PYG{p}{)}
\end{Verbatim}

\caption{Simulation of a soliton traveling in a tree-shaped quantum graph with a splitting scheme.}
\label{lst:splitting}
\end{listing}

\begin{figure}[htpb!]
  \centering
  \begin{tabular}{cccc}
    \includegraphics[height=.14\textheight]{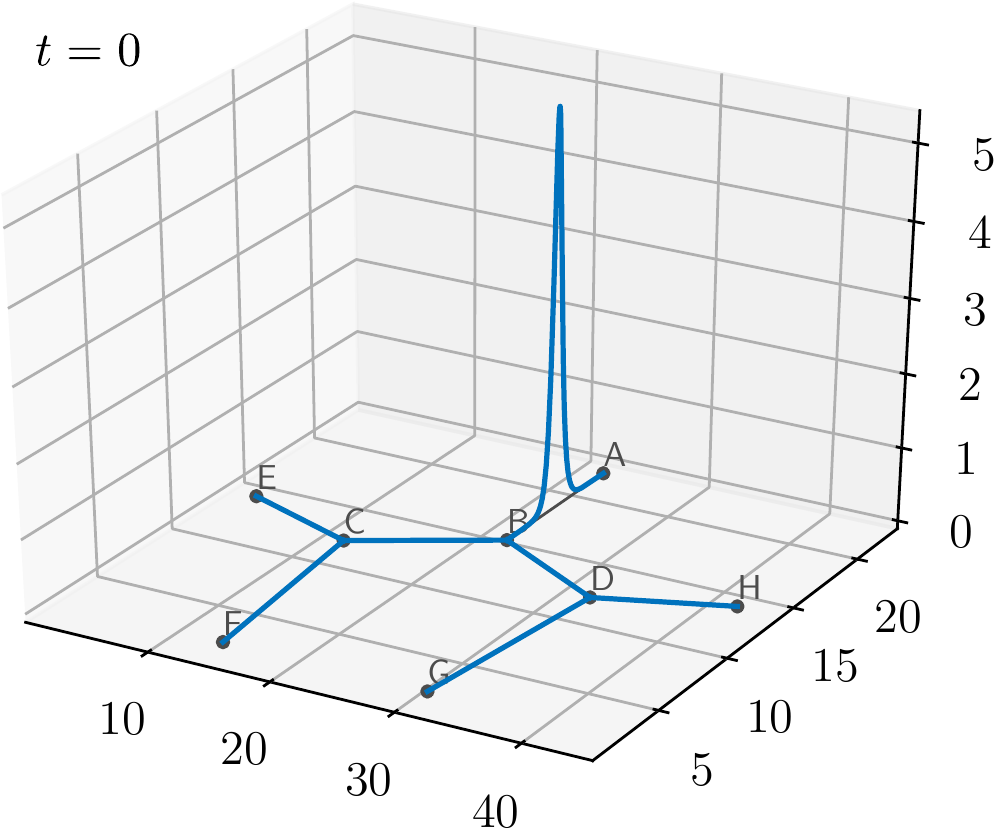}& \includegraphics[height=.14\textheight]{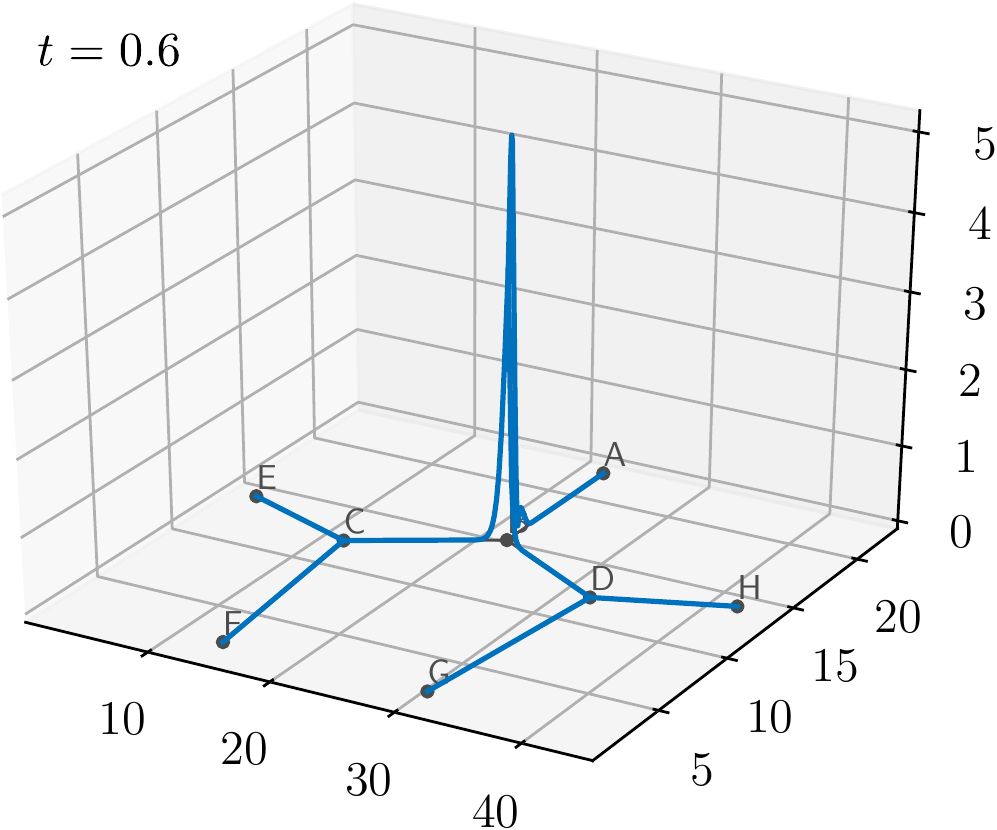}&
    \includegraphics[height=.14\textheight]{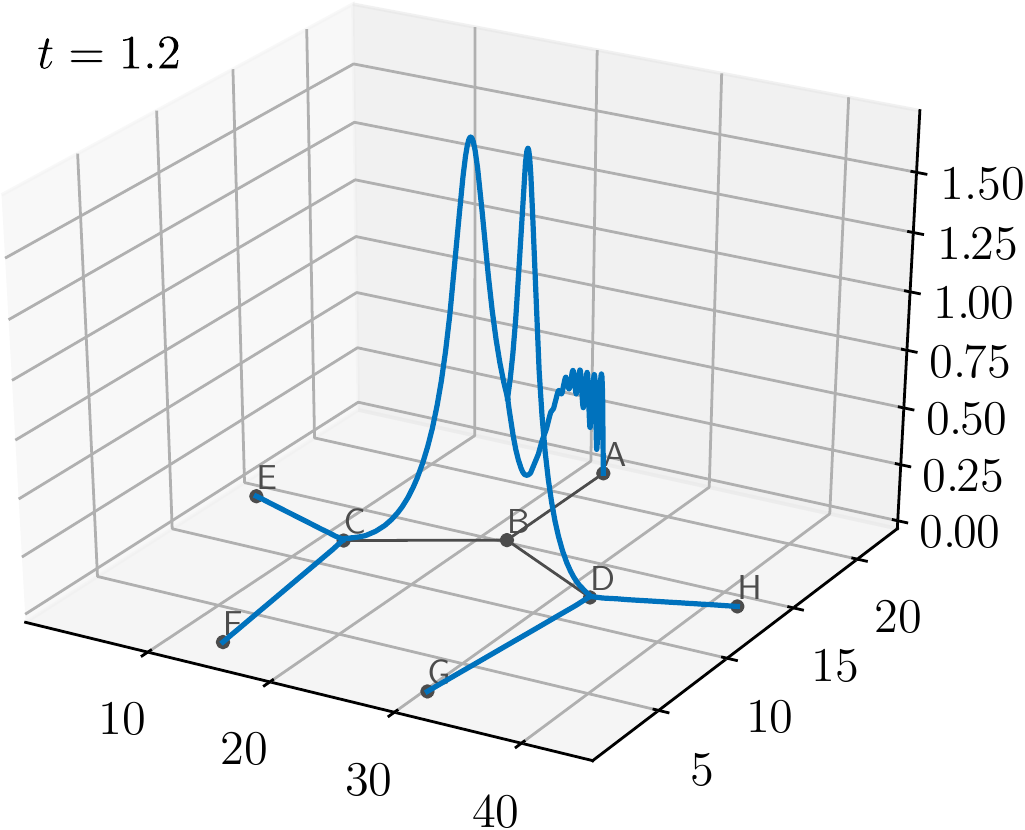}& \includegraphics[height=.14\textheight]{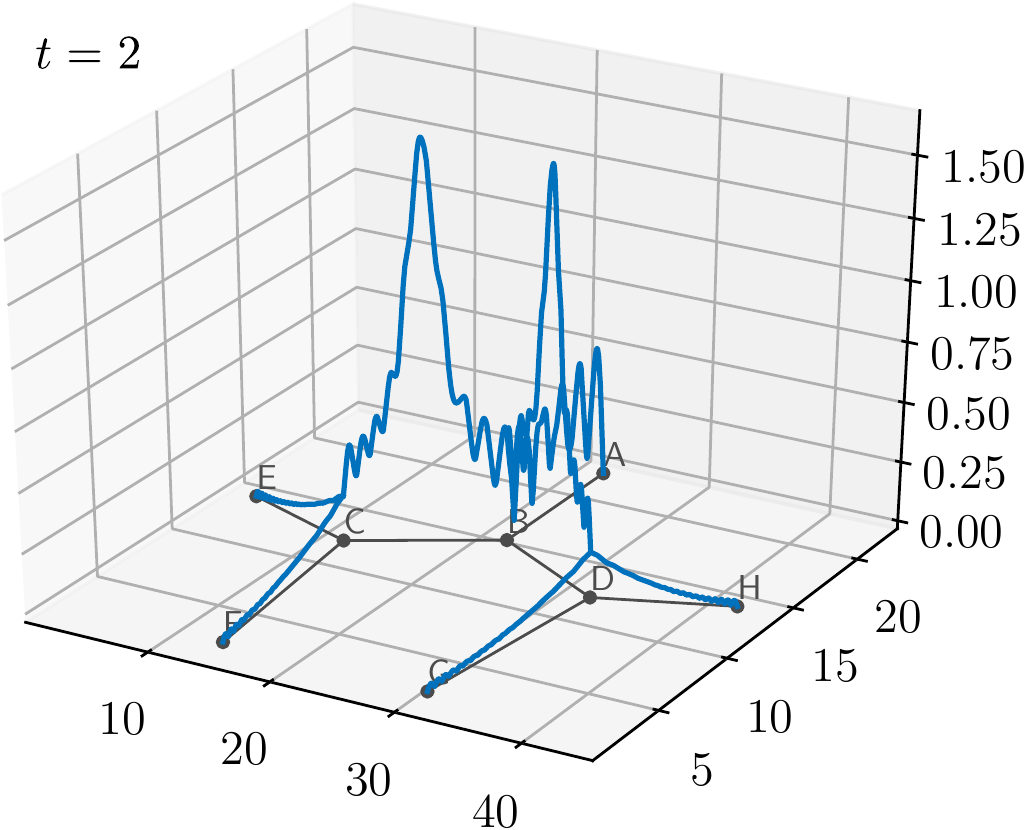}
  \end{tabular}
  \caption{Simulation of~\eqref{eq:dyn} at times $t = 0,0.6,1.2,2$ with a splitting scheme. }
  \label{fig:splittingsim}
\end{figure}

The result of the simulation can be seen in Figure~\ref{fig:splittingsim} where the absolute value of $\psi$ is plotted at different times. The ripples are expected to appear, as when reaching a vertex the solution will split between waves going through the vertex and a reflected wave, which will itself interact with the rest of the incident wave. 


\section{Ground states: numerical experiments and theoretical validation}
\label{sec:ground-states}

In this section, we present various numerical computations of ground states. In many cases,
explicit exact solutions are available. We show the efficiency of the CNGF scheme  for all these cases. Even though
the CNGF method was built for a general nonlinearity, we focus in this section
on the computations of the ground states of the focusing cubic nonlinear
Schr\"odinger (NLS) equation on
a graph $\mathcal{G}$, that reads
\begin{equation}
  \label{eq:nls2}
  i\psi_t = H \psi - \lambda |\psi|^{p-1} \psi,\quad x\in \mathcal{G},
\end{equation}
with $\lambda>0$.

Unless otherwise specified, we assume that $\lambda=1$ and $p=3$.

In what follows, we discuss only the power case nonlinearity and we focus on the results involving the obtention of ground states by minimization of the energy under a fixed mass constraint. It is in general not an easy task to prove that the standing wave profiles obtained by other techniques (e.g. bifurcation) are (or are not) minimizers at fixed mass (even locally).

Recall that for the classical nonlinear Schr\"odinger equation, the equation is said to be $L^2$-subcritical if (in one dimension) $1<p<5$. This is also the range of exponents for which standing waves are stables, and for which they can be obtained as minimizers of the energy at fixed mass. 
Metric graphs being based on one dimensional structures (segments and half-lines), the interesting range of exponents for the nonlinearity is $1<p\leq 5$, with the expectation of additional difficulties in the analysis at the critical case $p=5$. The global dimension of the graph might induce further restriction on the set of possible exponents, e.g. for the $2$-d grid $\mathbb Z^2$, which is locally $1$-d but globally $2$-d (we will comment on that later on).

For the cubic nonlinear Schr\"odinger equation on a finite (bounded or unbounded) graph, and at sufficiently large mass,  Berkolaiko, Marzuola and Pelinovsky~\cite{BeMaPe19}  established the following results. For any edge of the graph, there exists a bound state located  on the graph, i.e. it is positive, achieves its maximum on the edge, and the mass of the bound state is concentrated on the edge up to an exponentially small error (see~\cite[Theorem 1.1]{BeMaPe19} for a precise statement). Moreover, comparing the energies of these bound states, the authors are able to find the one with the smallest energy at fixed mass. Note, however, that the bound state with the smallest energy has not been proven yet to be the ground state. Heuristic arguments in favor of this hypothesis are given in~\cite[Section 4.4]{BeMaPe19}. 
The results of~\cite{BeMaPe19} have to be put in perspective with the results established by Adami, Serra and Tilli~\cite{AdSeTi19} for generic sub-critical power type nonlinearities. Indeed, by very elegant purely variational techniques,  Adami, Serra and Tilli~\cite{AdSeTi19}  established for non-compact graphs the existence of positive bound states achieving their maximum on any chosen finite edge. These bound states are obtained by purely variational techniques: it is proved that they are global minimizers of the energy among the class of functions with fixed mass, \emph{and the additional constraint that the functions should achieve their maximum on the given edge}. It turns out unexpectedly that the minimizer so obtained lies in fact in the interior of the constraint, hence it may also be characterized as a \emph{local} (but obviously not necessarily global) minimizer of the energy at fixed mass.
In the same vein,  the existence of local minimizers of the energy for fixed mass has also been established by Pierotti, Soave and Verzini~\cite{PiSoVe20} in cases where no ground state exists.
As the estimates~\cite[(4.6) and (4.7)]{BeMaPe19} indicate, a pendant edge is clearly preferable to a non-pendant one. However, for non-pendant edges, the differences between energies are quite small.

From the preceding discussion, we infer that extra-care is required when performing numerical experiments, as the outcome of the algorithm may very be only a local minimizer and not a global one.

We have divided this section into four parts, depending on the kind of graphs considered: compact graphs, graphs with finitely many edges, one of which is semi-infinite, periodic graphs and,finally, trees. If the vertices conditions are not specified, it means that Kirchhoff conditions are assumed.

\subsection{Compact graphs}

Compact graphs are made of a finite number of edges, all of which are of finite length. 
On compact graphs, the existence of minimizers in the subcritical case $1<p<5$ is granted by Gagliardo-Nirenberg inequality and the compactness of the injection of $H^1(\mathcal G)$ into $L^p(\mathcal G)$ for $1\leq p\leq \infty$. Hence the main question becomes to identify (or, in the absence of suitable candidates, to describe) the minimizer. Several works have been recently devoted to general compact graphs :~\cite{BeMaPe19,CaDoSe18,Do18,DoGhMiPi20,KuSh20}.
For the simplest compact graphs like the line segment or the ring, the minimizer is (usually) known and this offers us good test cases for our algorithm. Results applying to general compact graphs are not always easy to test numerically (e.g.  in~\cite{Do18}, Dovetta proved for any compact graph, for any $1<p<5$ and for any mass the existence of a sequence of bound state whose energy goes to infinity, but capturing this sequence at the numerical level would require the development of new specific tools). 
However, it was established in~\cite{CaDoSe18} that constant solutions on compact graphs are the ground state (for sub-critical nonlinearities) for sufficiently small mass, a feature which is easy to observe numerically.

The simplest of compact graphs are the segment (two vertices connected by an edge) and the ring (one vertex and an edge connecting the vertex to himself). As the ring case was considered in detail from a variational point of view in~\cite{GuLeTs17}, we chose to conduct experiments in this case and compare the numerical outcomes with the theoretical results of~\cite{GuLeTs17}. Beside the elementary cases of the segment and the ring, many compact graphs are of interest. We will present some experiments performed in the case of the dumbbell graph, for which several recent solid theoretical studies exist (see e.g.~\cite{Go19,MaPe16}).

\subsubsection{The ring}

\begin{figure}[htbp!]
  \centering
  \begin{tabular}{cc}
    \includegraphics[width=.31\textwidth]{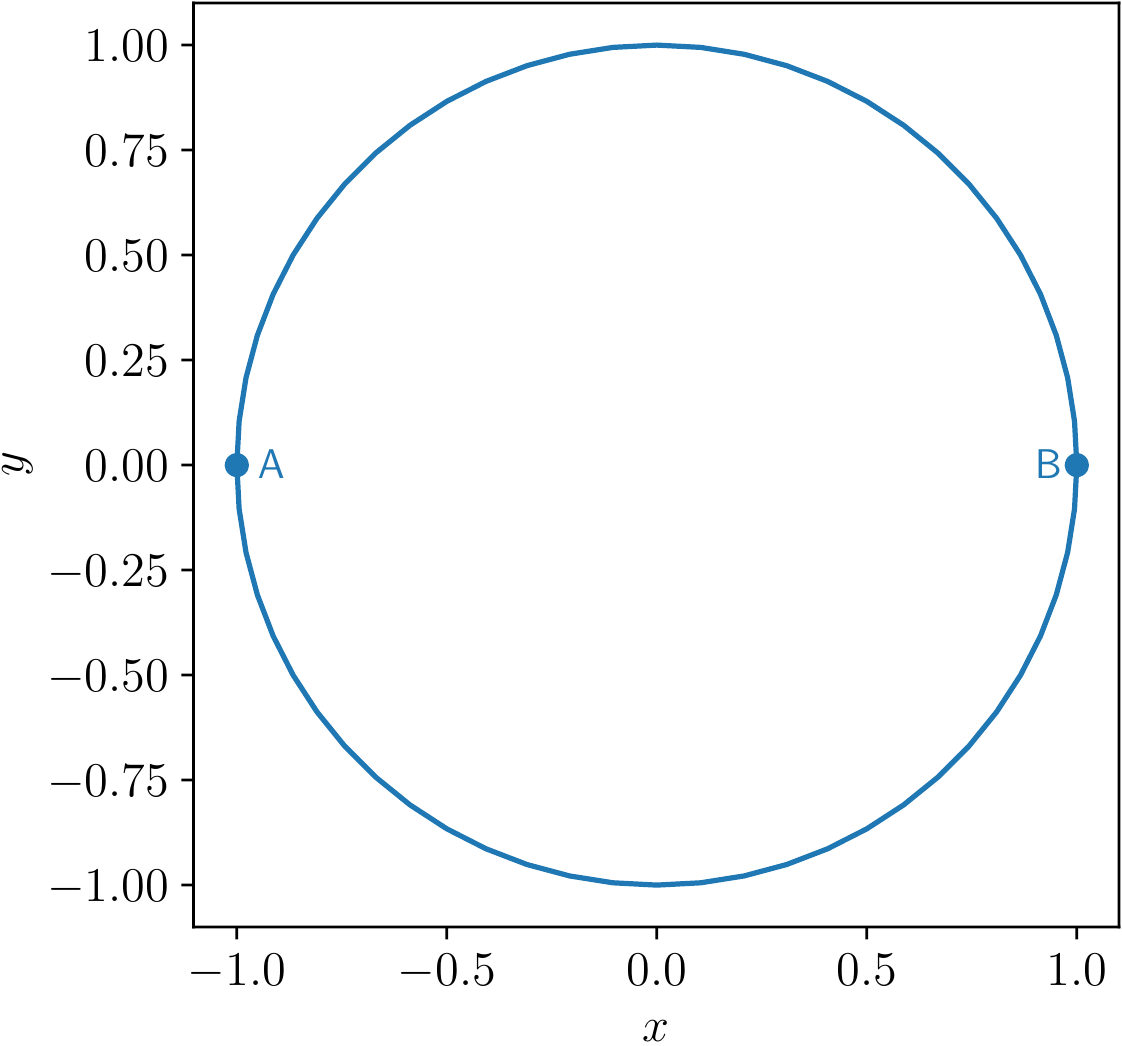} &
    \includegraphics[width=.31\textwidth]{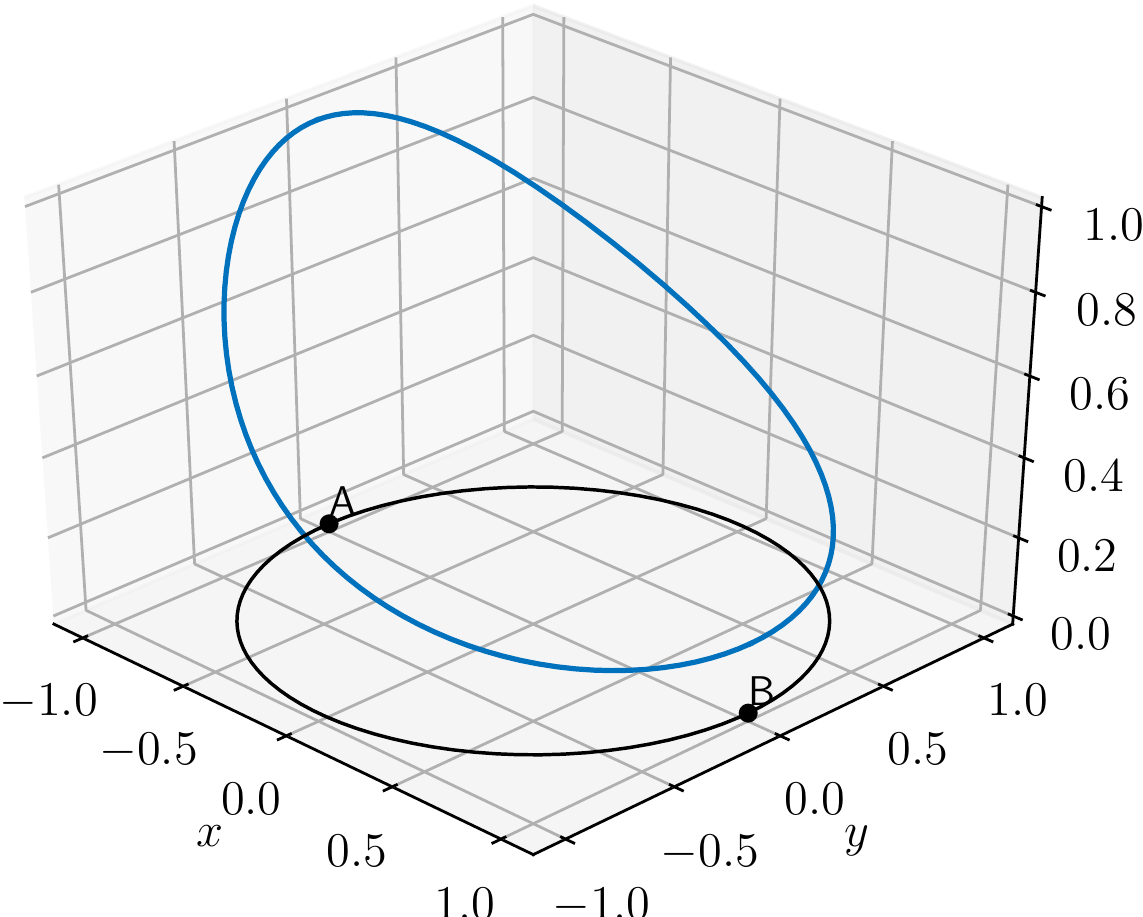}
  \end{tabular}
  \caption{Ring of radius $1$ (left) and the numerical ground state
    (right) when $\lambda=2$.}
  \label{fig:sol_one_loop}
\end{figure}

\begin{figure}[htbp!]
  \centering
  \begin{tabular}{cc}
    \includegraphics[width=.31\textwidth]{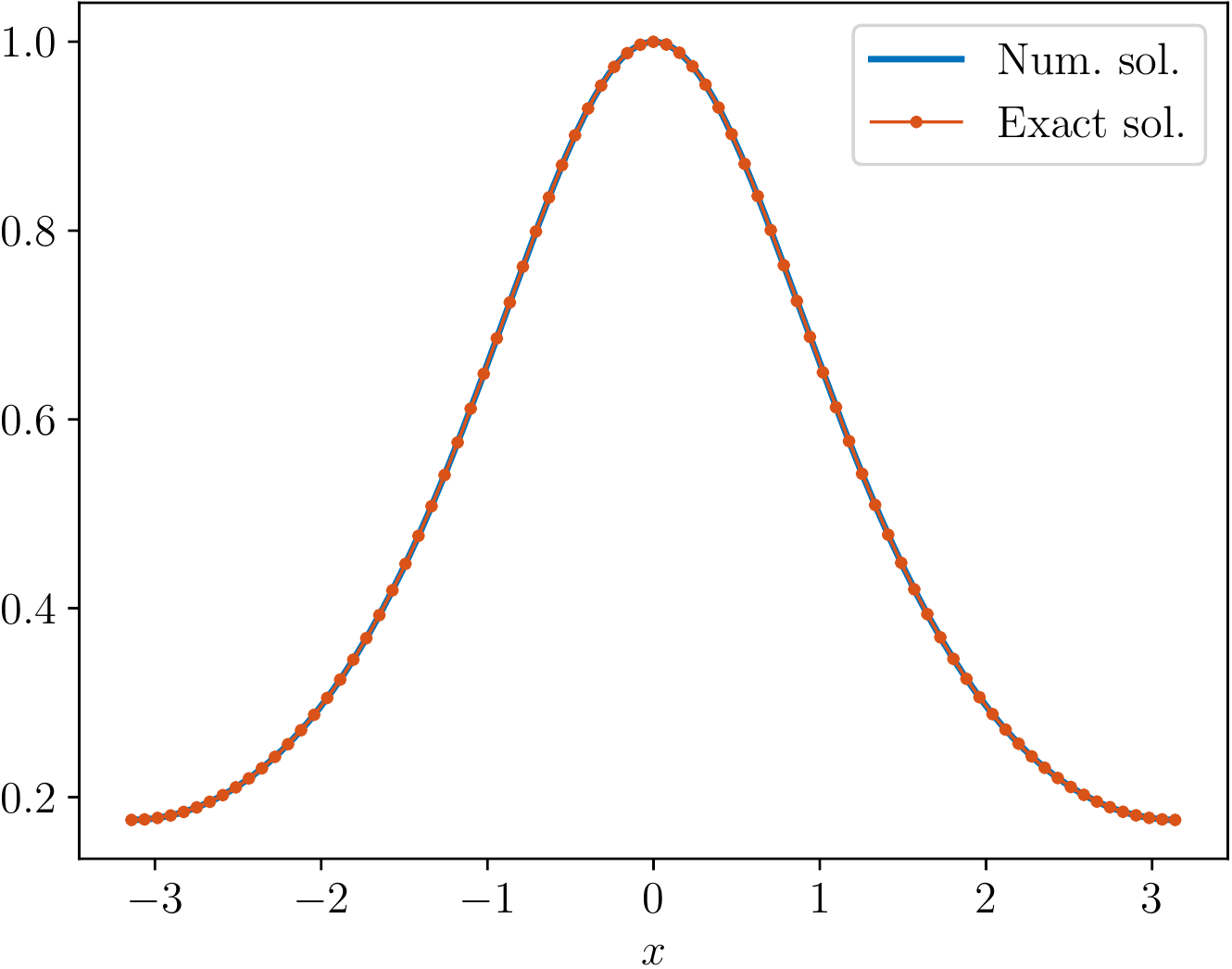} &
    \includegraphics[width=.31\textwidth]{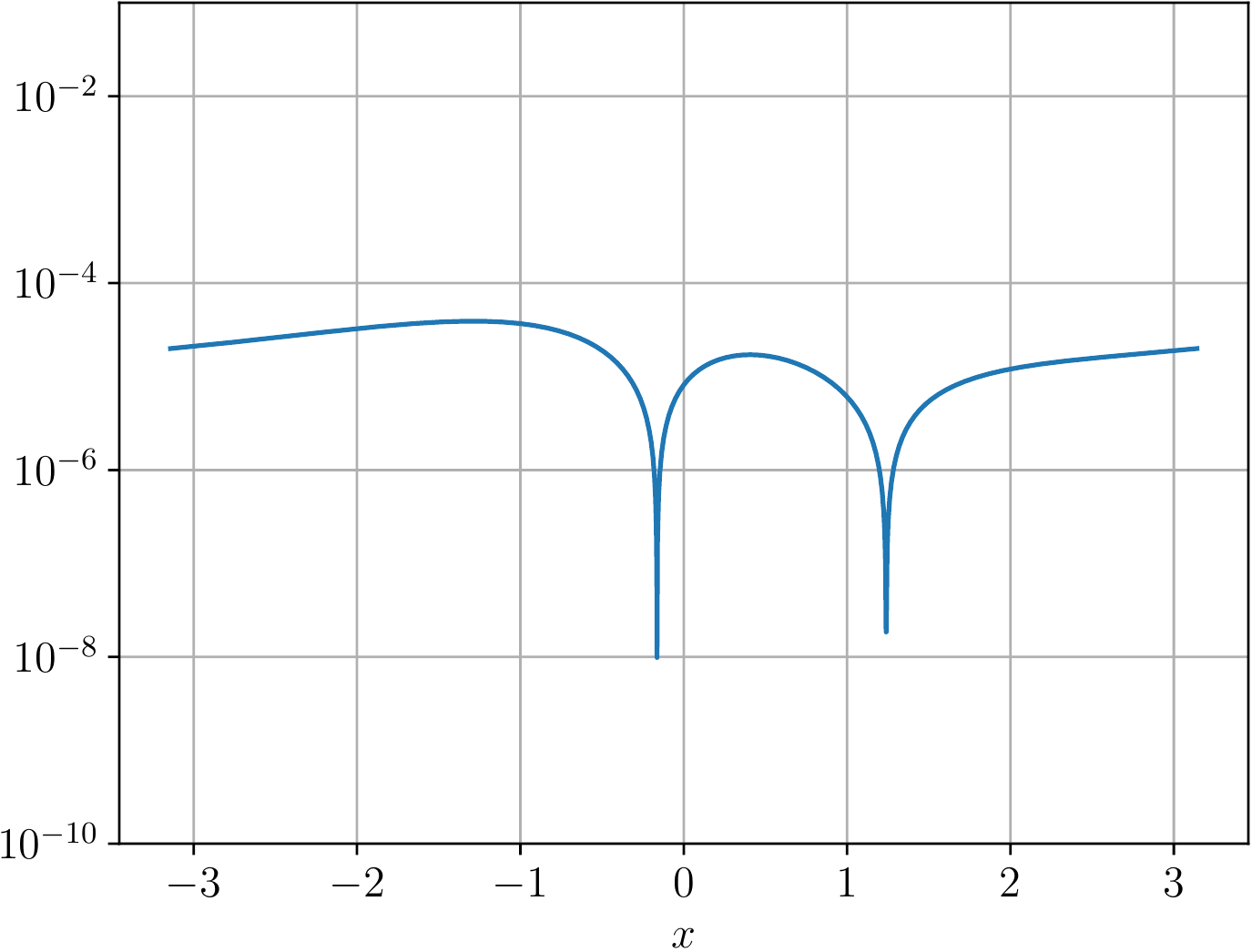}
  \end{tabular}
  \caption{Comparison between exact and numerical solution for
     the ring (where the $\log10$ of the difference is depicted on the right).}
  \label{fig:comp_sol_one_loop}
\end{figure}

From a numerical point of view, we obtain a ring (i.e. a one loop graph) by gluing together two half circles with Kirchhoff
conditions at the vertices (as already explained in Section~\ref{sec:simulations}, it is innocuous for the functions on the graphs).
Considering the loop graph with an edge of length $T$ is equivalent to work on the line $\R$ with $T$-periodic functions, i.e to work in the functional setting:
\[
H^1_{\text{loc}}\cap P_T,\quad P_T=\{f\in
L^2_{\text{loc}}(\mathbb{R}):\ f(x+T)=f(x),\ \forall x \in \mathbb{R}\}. 
  \]
Minimizers in $H^1_{\text{loc}}\cap P_T$ of the Schr\"odinger energy 
  \begin{equation}
  E_{\text{ring}}(\psi)=\frac12 \int_0^T|\psi'(x)|^2 \, dx - \frac14 \int_0^T
  |\psi(x)|^4\, dx
  \label{eq:minprob1}
\end{equation}
at fixed mass $m$ were described explicitly in~\cite{GuLeTs17} in terms of Jacobi elliptic functions.
 Recall that the
  function $\dn$ is the Jacobi elliptic function defined by
  \begin{equation}
    \label{eq:dn}
    \dn(x;k)=\sqrt{1-k^2 \sin^2(\phi)},\quad k\in(0,1),    
  \end{equation}
  where $\phi$ is defined through the inverse of the incomplete elliptic
  integral of the first kind
  \[
    x=F(\phi,k)=\int_0^\phi \frac{d\theta}{\sqrt{1-k^2\sin^2(\theta)}}.
  \]
 The snoidal and cnoidal functions are given 
  by
  \begin{equation}
    \label{eq:sn_cn}
    \text{sn}(x;k)=\sin(\phi),\qquad \text{cn}(x;k)=\cos(\phi).
  \end{equation}
  Recall also that the complete elliptic integrals of first and second kind are given by
  $K(k)=F(\pi/2;k)$ and $E(k)=E(\pi/2;k)$, where
  \[
    E(\phi;k)=\int_0^\phi \sqrt{1-k^2\sin^2(\theta)}\, d\theta.
  \]
The
solutions of the minimization problem~\eqref{eq:minprob1} are given as follows. 
\begin{enumerate}
\item For all $0<m<\frac{2\pi^2}{\lambda T}$, the unique minimizer (up to a
  phase shift) is the
  constant function
  \[
    \psi_{\text{ring}}=\sqrt{\frac{m}{T}}.
    \]
\item For all $  \frac{2\pi^2}{\lambda T}<m<\infty$, the unique minimizer (up to 
  phase shift and translation) is the rescaled dnoidal function
  \[
    \dn_{\alpha,\beta,k}=\frac{1}{\alpha} \dn\left(\frac{\cdot}{\beta};k\right)
    \]
  where the parameters $\alpha$, $\beta$ and $k$ are uniquely determined.
\item\label{item:3} If $\lambda=2$, given $k\in (0,1)$, $T=2K(k)$,  and $m=2E(k)$, the unique minimizer (up to phase shift and translation) 
  is
  \[
    \dn=\dn(\cdot,k).
  \]
\end{enumerate}
We place ourselves in the case of item (\ref{item:3}) and compute the ground state on the one loop graph
with perimeter $2\pi$ and $\lambda=2$ . The
parameter $k$ is therefore such that $k^2=0.9691073732421548$ and we fix the mass to $2E(k)$. We discretize each
half circle with $N_e=1000$ grid nodes. The gradient step is $\delta
t=10^{-2}$. Our experiment gives a remarkable agreement between the theoretical minimizer and the numerically computed minimizer, as shown in Figures~\ref{fig:sol_one_loop} and~\ref{fig:comp_sol_one_loop}.

\subsubsection{The dumbbell
}

\begin{figure}[htbp!]
  \centering
  \includegraphics[width=.5\textwidth]{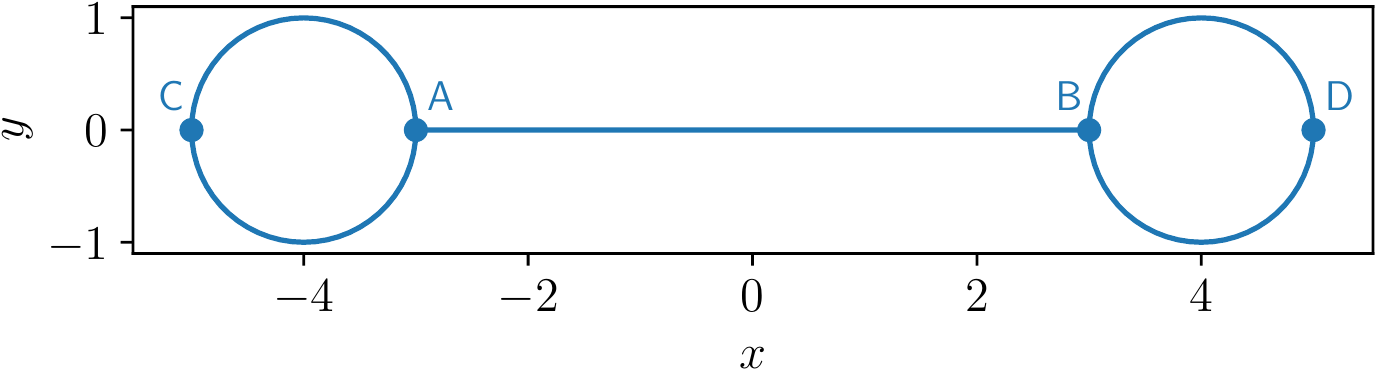}
  \caption{The dumbbell graph}
  \label{fig:dumbbell_graph}
\end{figure}

The dumbbell graph is a structure made of two rings attached to a central line
segment subject to Kirchhoff conditions at the junctions
(see Figure~\ref{fig:dumbbell_graph}).
Each ring
can be assimilated as a loop as in the previous section and  is
therefore numerically made of two glued half circles. The central line segment has a length
$2L$ and the perimeter of each loop is $2\pi$. We set $\lambda=2$, $L=3$ and consider the the
minimizers of the energy
\[
  E_{\text{dumbbell}}(\psi)=\frac12 \int_{\mathcal{G}}|\psi'(x)|^2 -\frac12
  |\psi(x)|^4 \, dx
\]
with fixed mass $M_{\text{dumbbell}}$. 
According to~\cite{MaPe16}, there exist $m^*$ and $m^{**}$ (explicitly known) such that $0<m^*<m^{**}<\infty$ and the following behavior for standing wave profiles on the dumbbell graph holds. For $0 < M_{\text{dumbbell}}<m^*$, the ground state is given by the constant
solution $\psi(x)=p$, where $p$ is implicitly given by
  \begin{equation*}
    M_{\text{dumbbell}} = 2(L+2\pi) p^2.
\end{equation*}
This constant solution undertakes a symmetry breaking bifurcation at $m^*$ and
a symmetry preserving bifurcation at $m^{**}$, which result in the appearance
of new positive non-constant solutions. The asymmetric standing wave is a ground
state for $M_{\text{dumbbell}} \gtrsim m^*$, and the symmetric standing wave is
not a ground state for $M_{\text{dumbbell}} \gtrsim m^{**}$. In our case, the values for $m^*$ and $m^{**}$ are
\[
  m^*=0.18646428284896863, \qquad m^{**}=1.2334076715778846.
\]
Observe that the three profiles described above are expected to be local minimizers of the energy at fixed mass, hence we should be able to find them with our numerical algorithm, provided the initial data is suitably chosen. We have found that the three following initial data were leading to the various desired behaviors (in the following, $\nu$ is a normalization constant adjusted in such a way that the mass constraint is verified):
\begin{itemize}
\item the constant initial data : $\psi_1\equiv \nu$,
\item a gaussian centered on the left loop : $\psi_2(x)_{|CA}=\nu e^{-10x^2}$ and $0$ elsewhere,
  \item a gaussian centered at $x=2$ on the central edge : $\psi_3(x)_{|AB}=\nu e^{-10(x-2)^2}$ and $0$ elsewhere.
\end{itemize}

\begin{figure}[htbp!]
  \centering
  \begin{subfigure}[t]{0.31\textwidth}
    \centering
    \includegraphics[width=\textwidth]{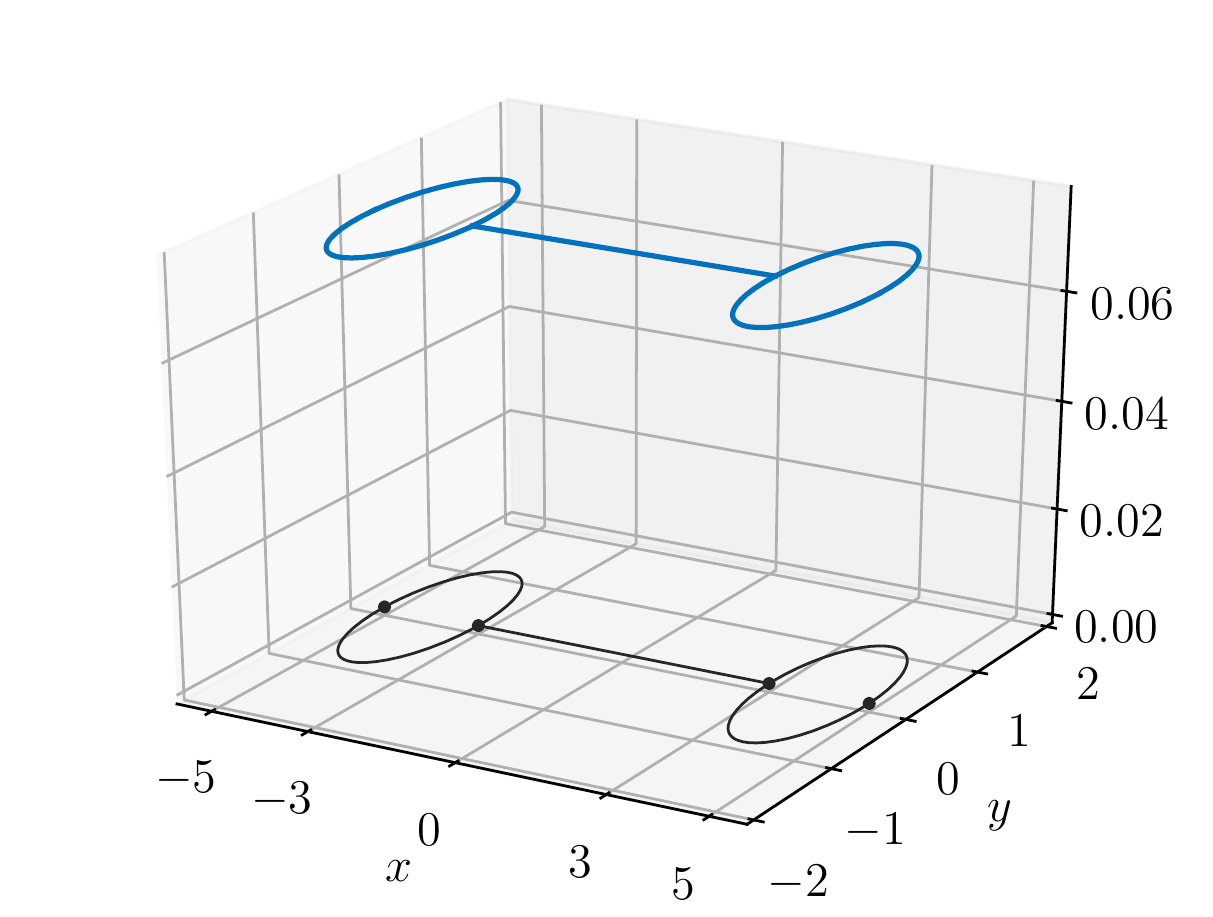}
  \caption{$m = 0.10$, $u_0=\psi_1$,\\ $E=-2.6930411461.10^{-4}$}
\end{subfigure}
  \begin{subfigure}[t]{0.31\textwidth}
    \centering
    \includegraphics[width=\textwidth]{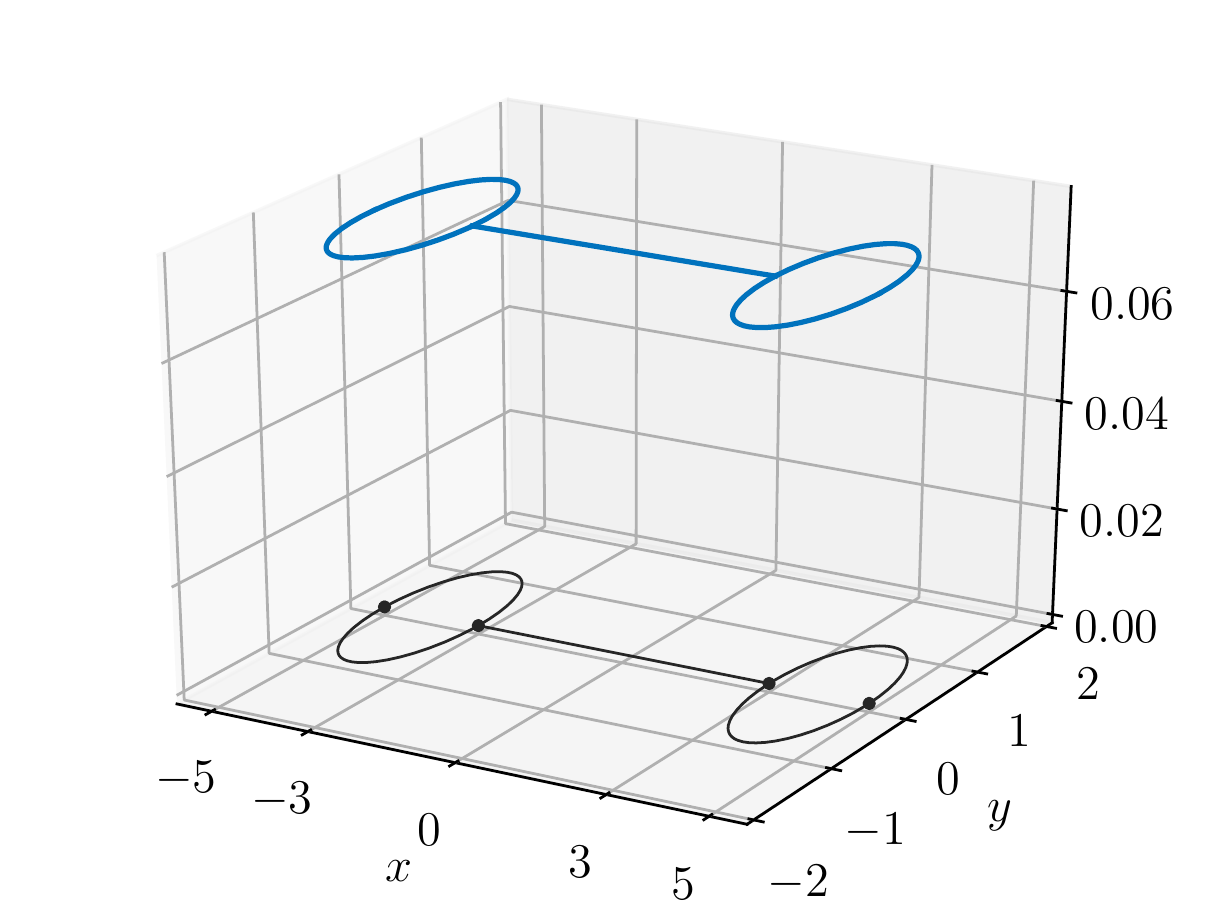}
  \caption{$m = 0.10$, $u_0=\psi_2$,\\ $E=-2.6930411103.10^{-4}$}
\end{subfigure}
  \begin{subfigure}[t]{0.31\textwidth}
  \centering
  \includegraphics[width=\textwidth]{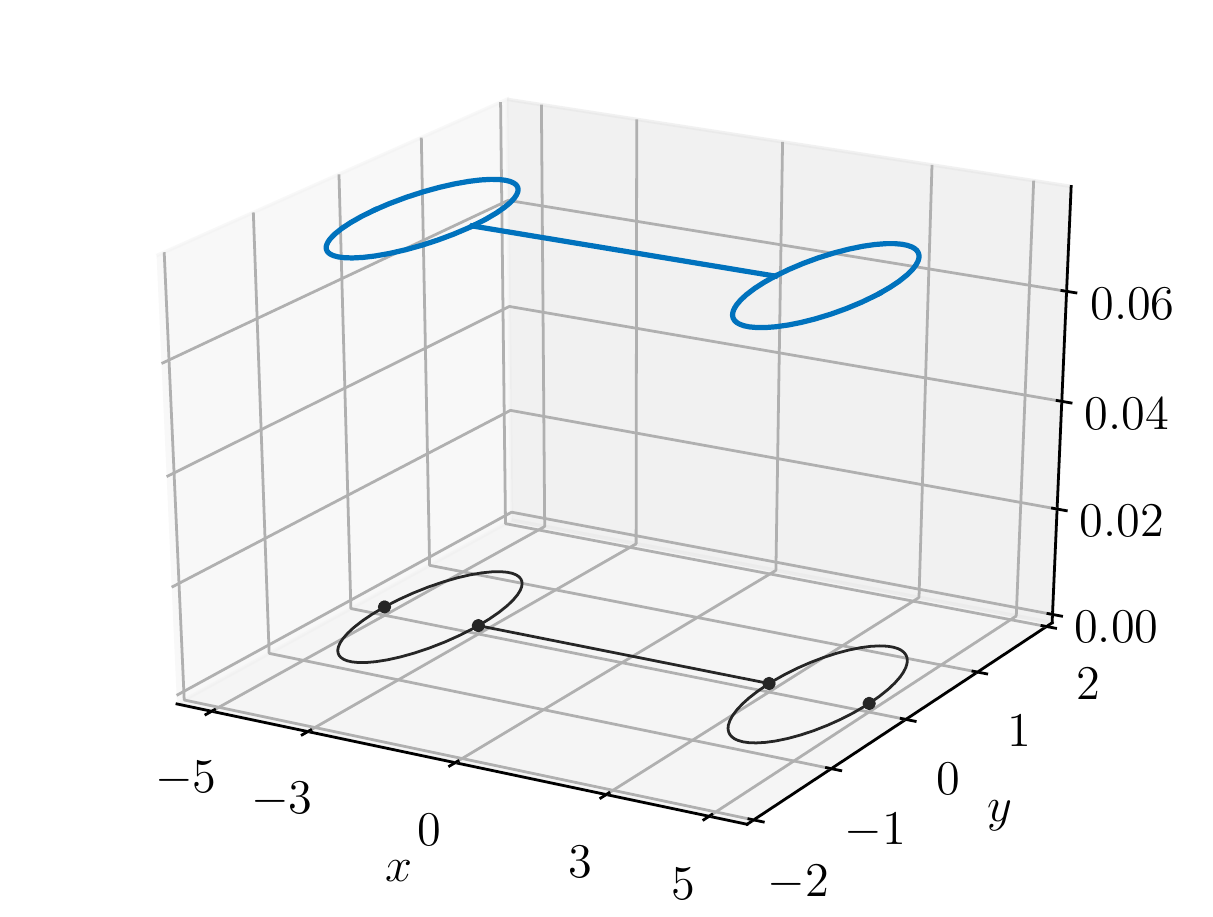}
  \caption{$m = 0.10$, $u_0=\psi_3$,\\ $E=-2.6930411193.10^{-4}$}
\end{subfigure}
\\
  \begin{subfigure}[b]{0.31\textwidth}
    \centering
    \includegraphics[width=\textwidth]{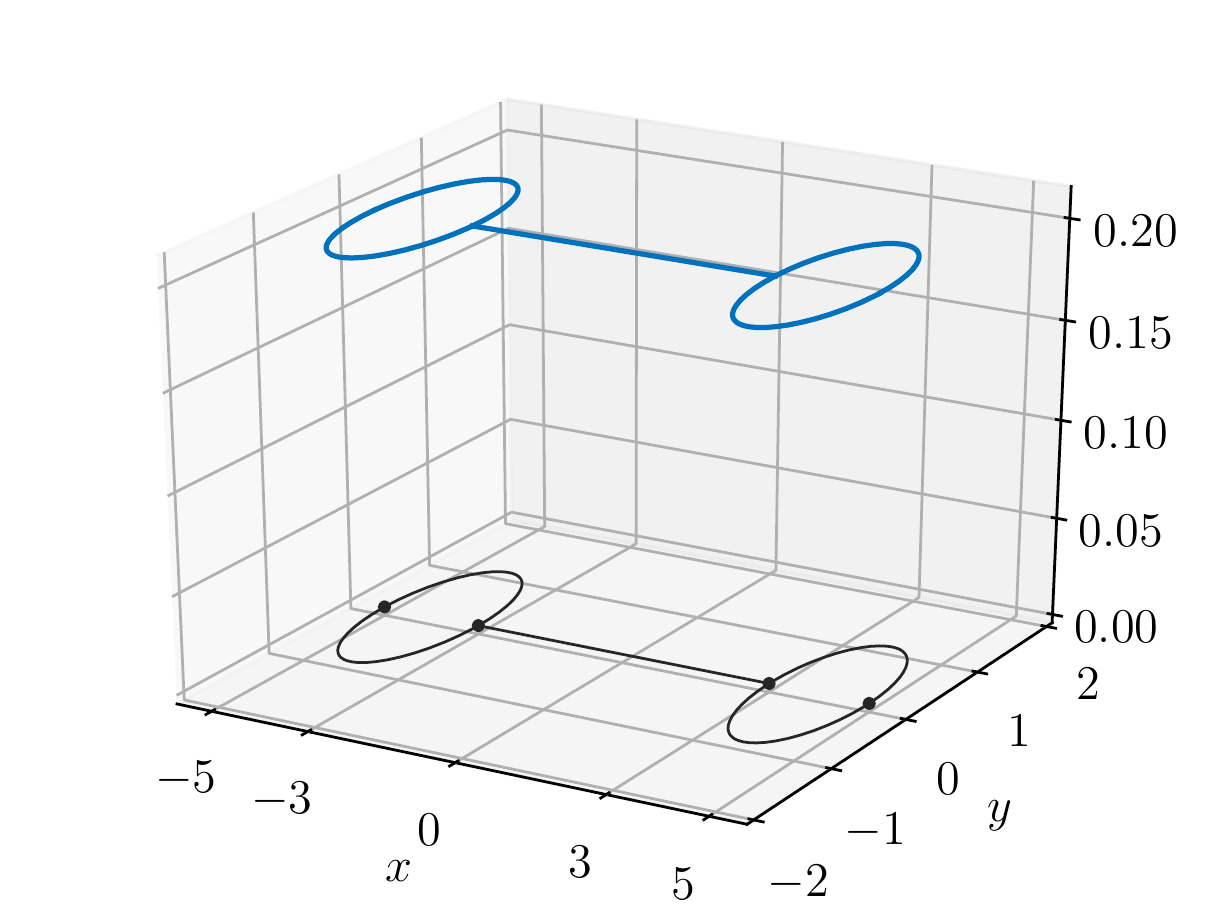}
  \caption{$m = 0.75$, $u_0=\psi_1$,\\  $E=-1.5148356447.10^{-2}$}
\end{subfigure}
  \begin{subfigure}[b]{0.31\textwidth}
  \centering
  \includegraphics[width=\textwidth]{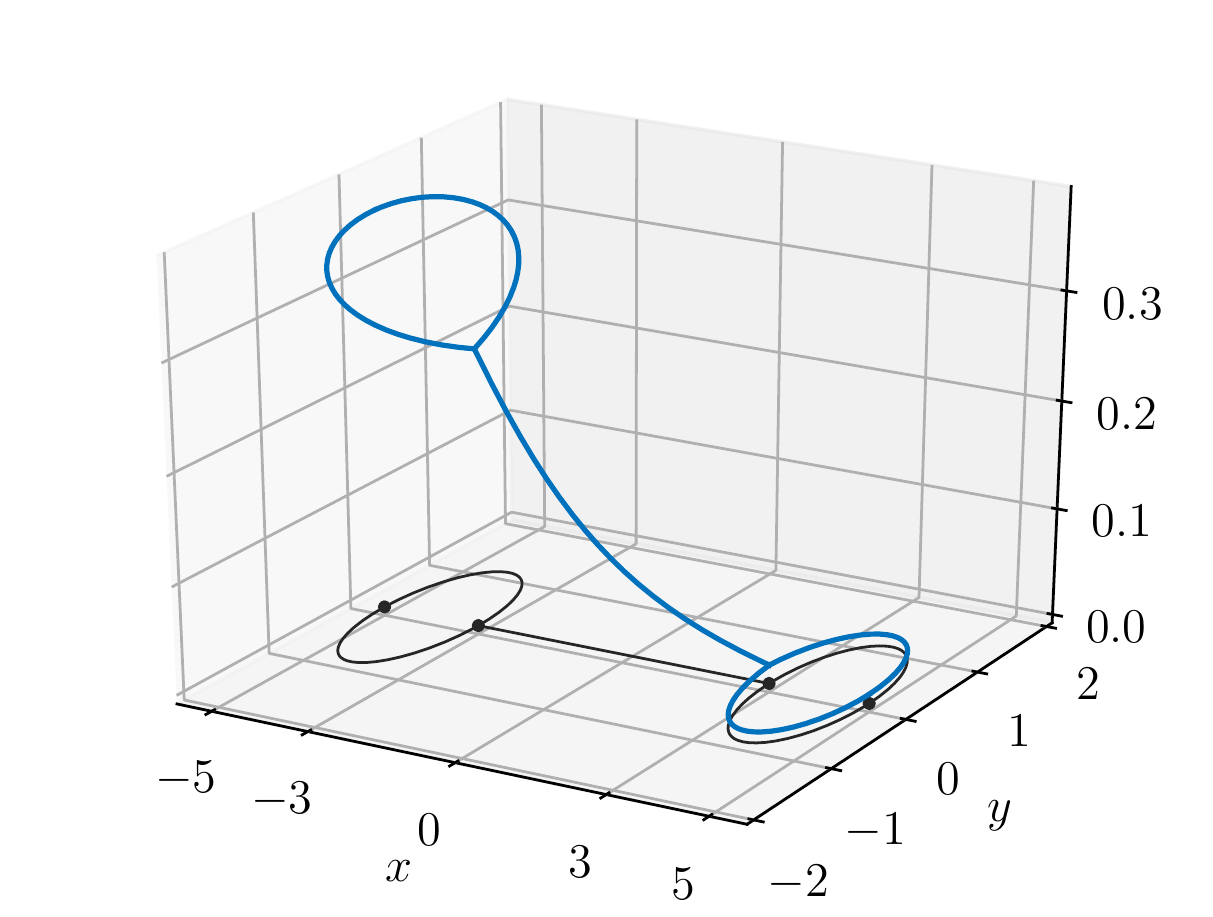}
  \caption{$m = 0.75$, $u_0=\psi_2$,\\ $E=-2.7205037742.10^{-2}$}
\end{subfigure}
\begin{subfigure}[b]{0.31\textwidth}
  \centering
  \includegraphics[width=\textwidth]{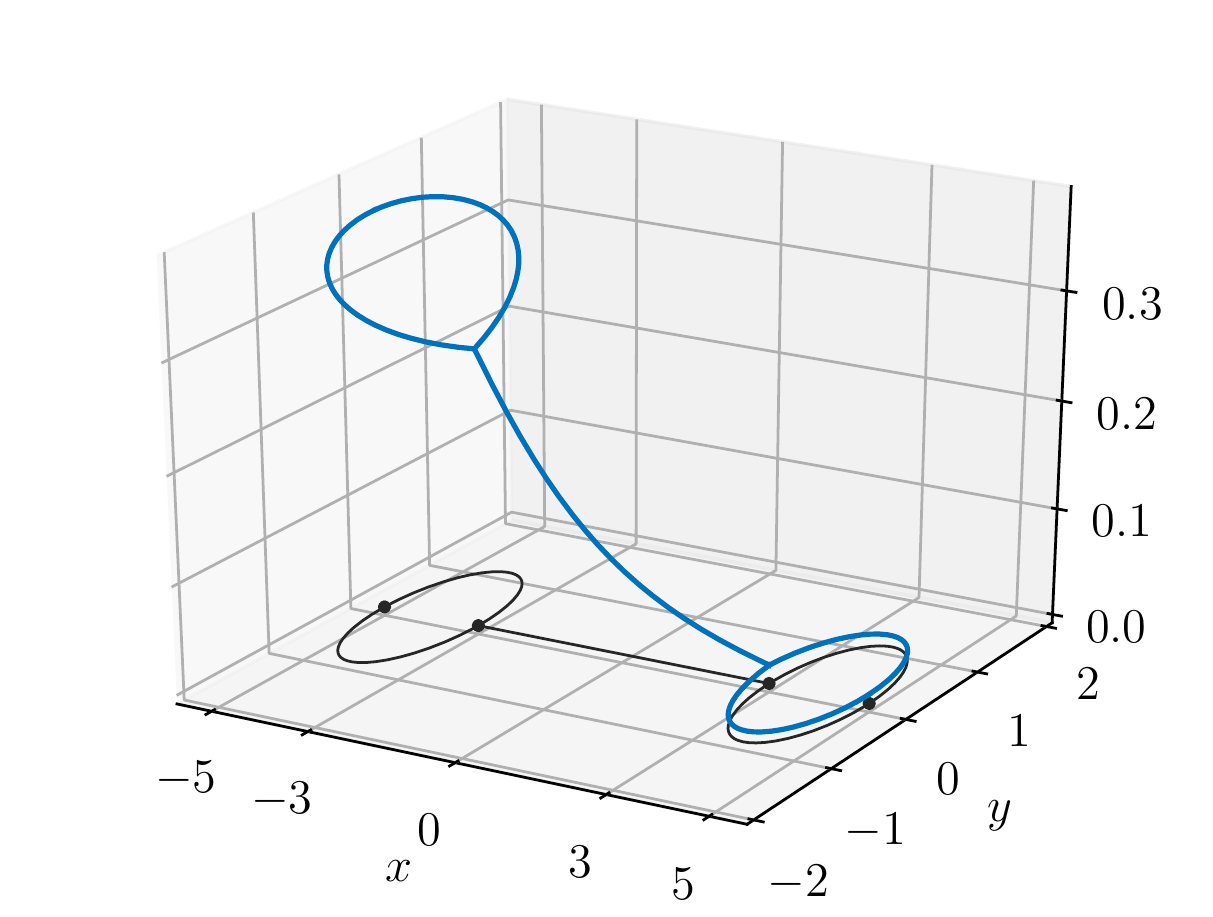}
  \caption{$m = 0.75$, $u_0=\psi_3$,\\ $E=-2.7205037743.10^{-2}$}
\end{subfigure}
\\
  \begin{subfigure}[b]{0.31\textwidth}
    \centering
    \includegraphics[width=\textwidth]{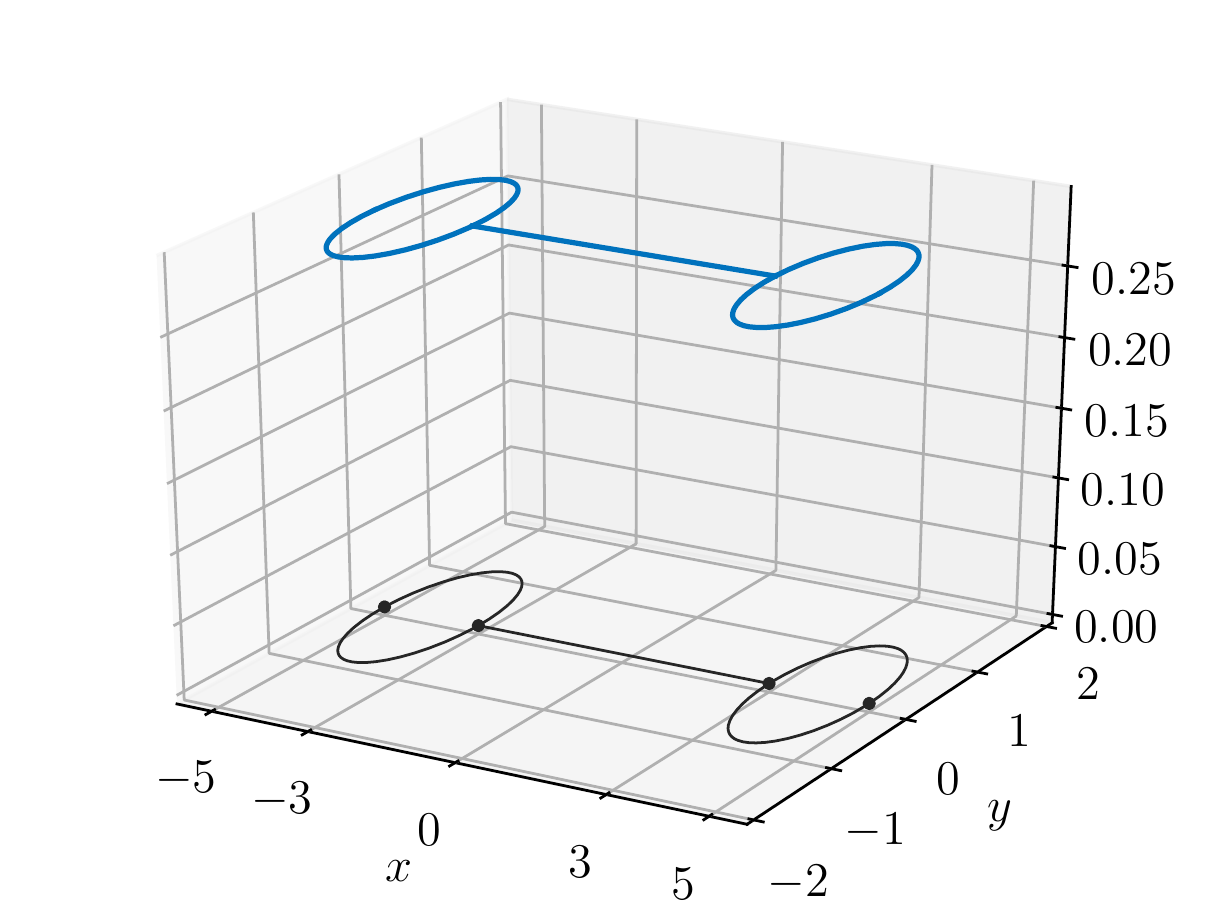}
  \caption{$m = 1.50$, $u_0=\psi_1$,\\  $E=-6.0593425789.10^{-2}$}
\end{subfigure}
  \begin{subfigure}[b]{0.31\textwidth}
  \centering
  \includegraphics[width=\textwidth]{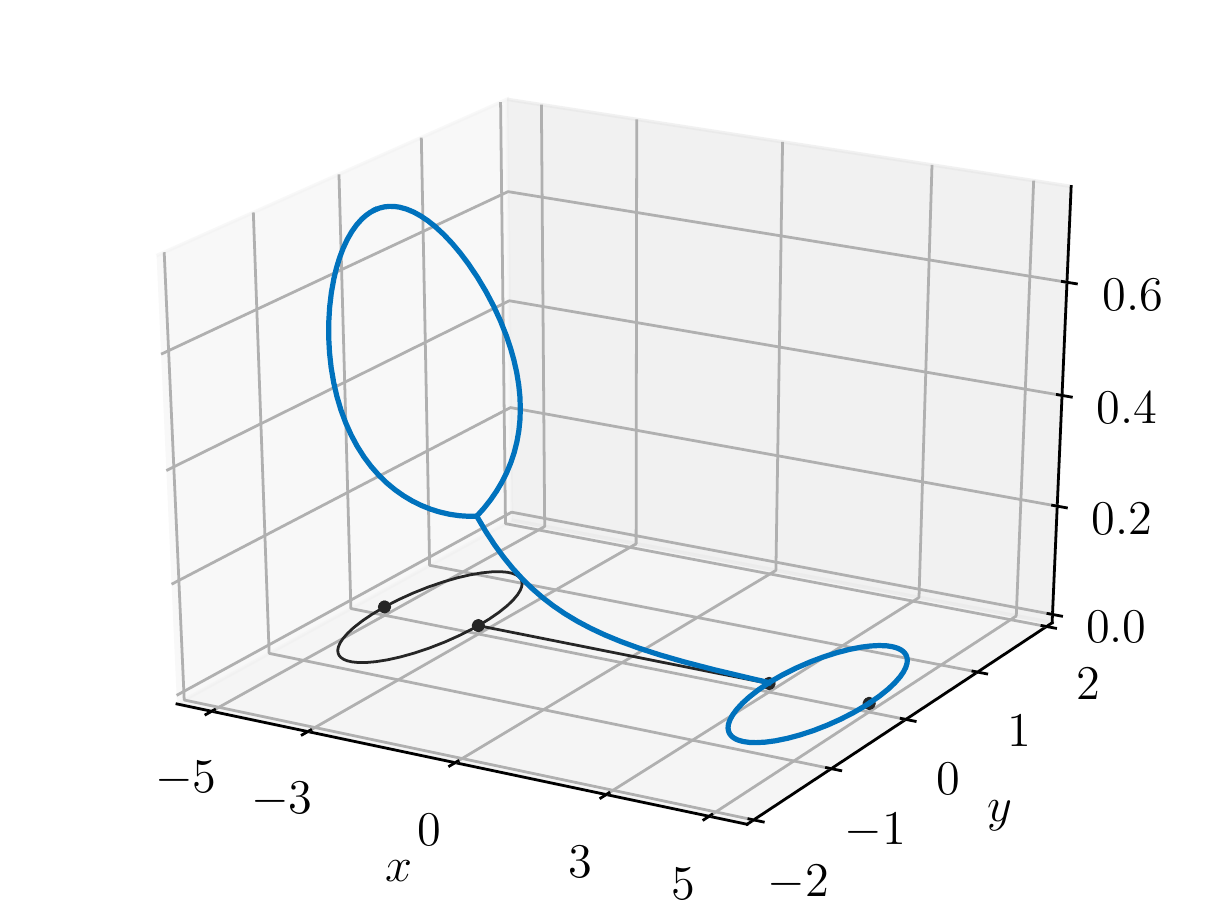}
  \caption{$m = 1.50$, $u_0=\psi_2$,\\ $E=-1.5097807829.10^{-1}$}
\end{subfigure}
\begin{subfigure}[b]{0.31\textwidth}
  \centering
  \includegraphics[width=\textwidth]{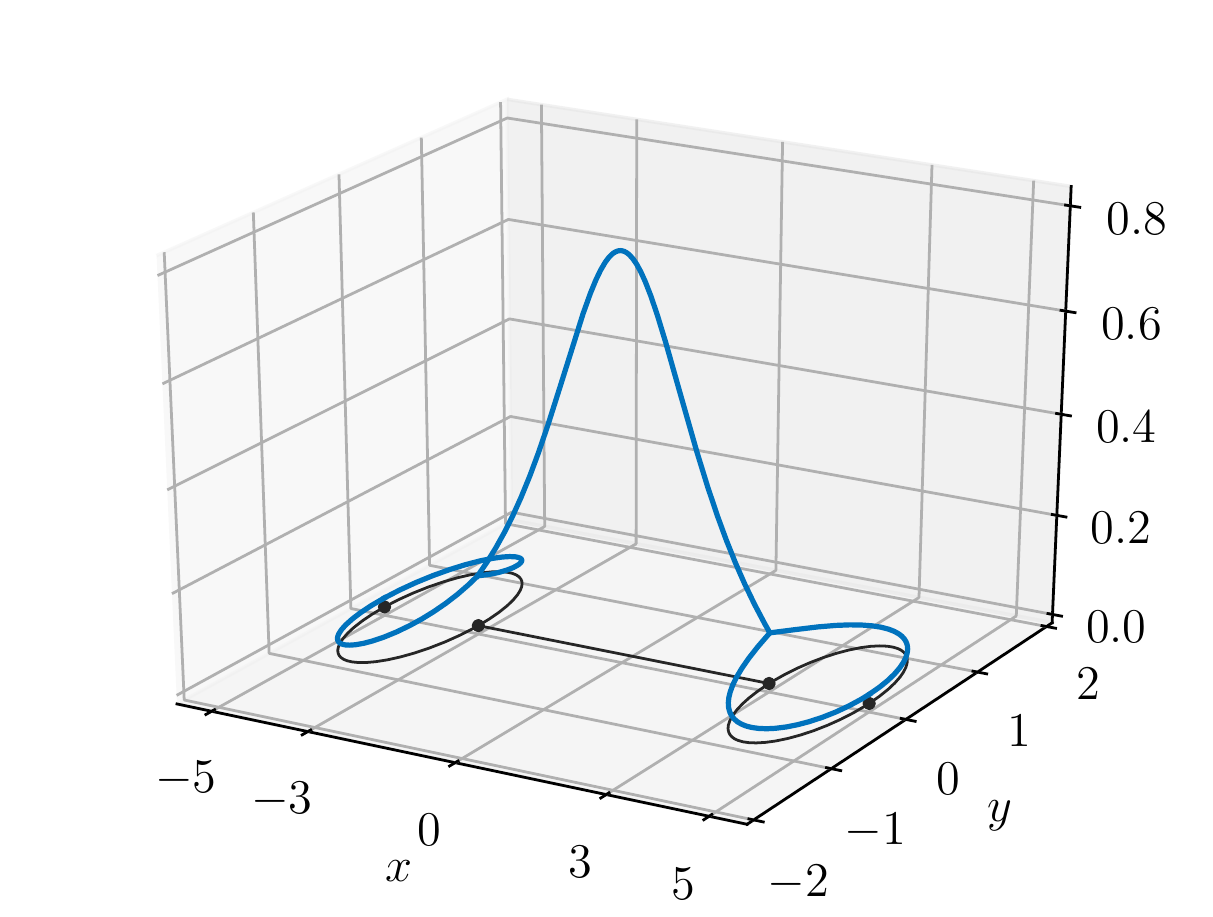}
  \caption{$m = 1.50$, $u_0=\psi_3$,\\  $E=-1.2925753851.10^{-1}$}
\end{subfigure}
\caption{Outcomes of the CNGF Algorithm~\ref{alg:CNGF} on the dumbbell graph for three remarkable values of the mass (one for each row) and three possible initial data (one for each column). In each subcaption, the quantity $E$ given is the final energy.}
\label{fig:very-dumb-bell}
\end{figure}

We will run the normalized gradient flow for each of these initial data for the three following masses:
\[
0<m_1=0.10<m^*,\quad m^*<m_2=0.75<m^{**},\quad m^{**}<m_3=1.50.
\]
The parameters of the algorithm are set as follows.  The total number of discretization nodes is $N=1000$ and $\delta t=10^{-2}$. The stopping criterion is set to $10^{-8}$, and the maximal number of iteration is set to $50000$ (which is large enough so that it is never reached in our experiences). The results are in perfect agreement with the theoretical results, as shown in Figure~\ref{fig:very-dumb-bell}. In particular, one can see that for large mass $m_3>m^{**}$, it is indeed possible to recover the three bound states described theoretically, and comparison of the energies shows that the asymmetric bound state centered on a loop is indeed the ground state. For the smaller mass $m_2$, the algorithm selects the constant or the asymmetric state, and comparison of the energy shows that the later is indeed the ground state. And for $m_1$, the algorithm converges in each case towards the constant function. Very small differences in the final energies (after the eighth digit in the $m_1$ case) may be noted, which are due to our stopping criterion set at $10^{-8}$. 


Compacity Compactnessfor graphs may be violated in several ways: with a semi-infinite edge, or with an infinite number of edges (which may be arranged e.g. periodically or in tree form). We discuss these cases in the next sections.

\subsection{Graphs with a semi-infinite edge}

In this section, we consider graphs having a finite number of edges, one of which is of semi-infinite length. A typical example for this kind of graph is the $N$-star graph, consisting of a vertex to which $N$ semi-infinite edges are attached. We will discuss this example in Section~\ref{sec:N_star_graph}. Before that, we will recall in Section~\ref{sec:obstruction} some of the results obtained by Adami and co. concerning a topological obstruction leading to non existence of ground states on nonlinear quantum graphs. Another example, the tadpole graph, will be discussed in Section~\ref{sec:tadpole}. 

\subsubsection{The topological obstruction}
\label{sec:obstruction}

The existence of ground states with prescribed mass for the focusing nonlinear
Schr\"odinger equation~\eqref{eq:nls2} on non-compact finite graphs $\mathcal{G}$ equipped with
Kirchhoff conditions at the vertices is linked to the topology of the graph. Actually, a
topological hypothesis {(H)} can prevent a graph from having ground states for every
value of the mass (see~\cite{AdSeTi17b} for a review). For the sake of clarity,
we recall that a \textsl{trail} in a graph is a path made of adjacent edges, in
which every edge is run through exactly once. In a trail vertices can be run
through more than once. The assumption {(H)} has many formulations
(again, see~\cite{AdSeTi17b}) but we give here only one.
\begin{assumption}[Assumption (H)]
  \label{ass:H}
  Every $x\in \mathcal{G}$ lies in a trail that contains two half-lines.
\end{assumption}

If a finite non-compact graph with Kirchhoff conditions at the vertices verifies Assumption~\ref{ass:H}, then no ground state exists, unless the graph is isomorphic to a tower of bubbles (see Figure~\ref{fig:H-yes}). 
Examples of graphs verifying Assumption~\ref{ass:H} abound, some are drawn on Figure~\ref{fig:H-yes}.
\tikzset{solid node/.style={circle,draw,inner sep=1,fill=black}}
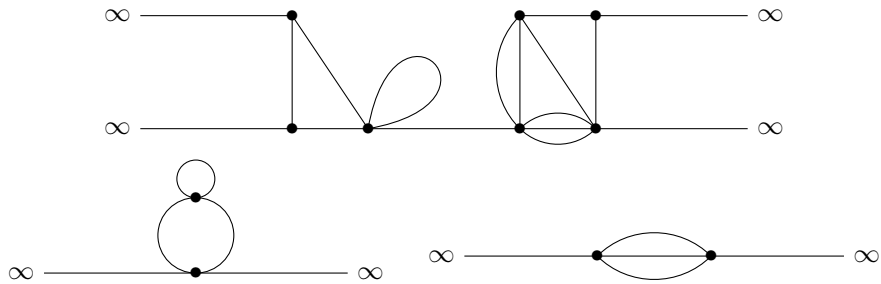
\begin{figure}[htbp!]
  \centering
    \begin{tabular}{c}
  \begin{tikzpicture}
    \node[left] at (-4,1.5) {$\infty$};
    \node[left] at (-4,0) {$\infty$};
    \node[right] at (4,0) {$\infty$};
    \node[right] at (4,1.5) {$\infty$};

    \node at (-2,1.5) {$\bullet$};
    \node at (-2,0) {$\bullet$};
    \node at (-1,0) {$\bullet$};
    \node at (1,0) {$\bullet$};
    \node at (1,1.5) {$\bullet$};
    \node at (2,0) {$\bullet$};
    \node at (2,1.5) {$\bullet$};

    \draw (-4,0) to (4,0);
    \draw (-2,0) -- (-2,1.5);
    \draw (-4,1.5) -- (-2,1.5) -- (-1,0);
    \draw (1,0) -- (1,1.5) -- (2,1.5) -- (2,0) -- (1,1.5);
    \draw (1,0) to [bend left=45] (1,1.5);
    \draw (1,0) to [bend left=45] (2,0);
    \draw (1,0) to [bend right=45] (2,0);
    \draw (2,1.5) -- (4,1.5);

    \draw[scale=4] (-0.25,0)  to[out=10,in=80,loop] (-0.25,0);
  \end{tikzpicture}\\
      \begin{tabular}{cc}
                         \begin{tikzpicture}
      \node[left] at (-2,0) {$\infty$};
      \node[right] at (2,0) {$\infty$};
      \draw (-2,0) -- (0,0)  -- (2,0);
      \node at (0,0) {$\bullet$};
      \draw (0,0.5) circle (0.5);
      \node at (0,1) {$\bullet$};
      \draw (0,1.25) circle (0.25);
    \end{tikzpicture}&
                        \begin{tikzpicture}
    \node[left] at (-2.5,0) {$\infty$};
    \node[right] at (2.5,0) {$\infty$};
    \draw (-2.5,0) to (2.5,0);
    \draw (-.75,0) to [bend right=45] (.75,0);
    \draw (-.75,0) to [bend left=45] (.75,0);
    \node at (-.75,0) {$\bullet$};
    \node at (.75,0) {$\bullet$};
  \end{tikzpicture}
       \end{tabular}
      \end{tabular}
  \caption{Graphs satisfying Assumption (H) : a generic graph (top), a tower of bubble (bottom left), a triple bridge (bottom right)}
  \label{fig:H-yes}
\end{figure}
Fortunately, graphs not satisfying Assumption~\ref{ass:H} and for which ground states exist also abound, some are shown on Figure~\ref{fig:H-no}.


\begin{figure}[htbp!]
  \centering
  \begin{tabular}{cc}
    \begin{tikzpicture}
      \node[left] at (-2,0) {$\infty$};
      \node at (1,0) {$\bullet$};
      \draw (-2,0) -- (1,0);
      \draw (1.5,0) circle (0.5);
    \end{tikzpicture}&
    \begin{tikzpicture}
      \node[left] at (-2,0) {$\infty$};
      \node at (1,0) {$\bullet$};
      \draw (-2,0) -- (1,0);
      \draw (1.5,0.75) -- (1,0) -- (2,0) -- (1,0) -- (1.3,-0.6);
      \node at (1.5,0.75) {$\bullet$};
      \node at (2,0) {$\bullet$};
      \node at (1.3,-0.6) {$\bullet$};
    \end{tikzpicture}\\
    \begin{tikzpicture}
      \node[left] at (-2,0) {$\infty$};
      \node[right] at (2,0) {$\infty$};
      \draw (-2,0) -- (0,0) -- (0,1) -- (0,0) -- (2,0);
      \node at (0,0) {$\bullet$};
      \node at (0,1) {$\bullet$};
    \end{tikzpicture}&
    \begin{tikzpicture}
      \node[left] at (-2,0) {$\infty$};
      \node[right] at (2,0) {$\infty$};
      \draw (-2,0) -- (0,0) -- (0,1) -- (0,0) -- (2,0);
      \node at (0,0) {$\bullet$};
      \node at (0,1) {$\bullet$};
      \draw (0,1.25) circle (0.25);
    \end{tikzpicture}
    \\
  \end{tabular}
    \caption{Graphs not satisfying Assumption (H) : a tadpole (top left), a $3$-fork (top right), a line with a pendant (bottom left), a sign-post (bottom right)}
  \label{fig:H-no}
\end{figure}
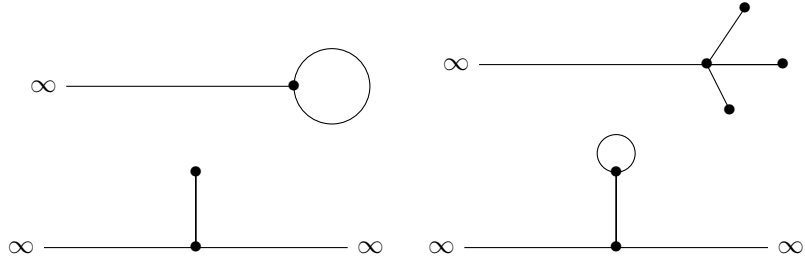


\subsubsection{Star graphs}
\label{sec:N_star_graph}

Star-graphs provide typical examples for nonlinear quantum graphs, as they are non-trivial graphs but retain many features of the well-studied half-line.  As star-graphs with Kirchhoff condition at the vertex verify Assumption~\ref{ass:H} and therefore do not possess a ground state, one usually studies star graphs with other vertex conditions such as $\delta$ or $\delta'$ conditions.

\begin{figure}[htpb!]
  \centering
\begin{tikzpicture}[scale=0.065]
  \node[label=above:$A$] at (0,0) {$\bullet$};
  \draw[<->,>=latex] (-40,0) -- (40,0);
  \draw[<->,>=latex] (-20,-34.64) -- (20,34.64);
  \draw[<->,>=latex] (-20,34.64) -- (20,-34.64);
  \node at (-24,-34.64) {$x_1$};
  \node at (24,-34.64) {$x_2$};
  \node at (44,0) {$x_3$};
  \node at (24,34.64) {$x_4$};
  \node at (-24,34.64) {$x_5$};
  \node at (-44,0) {$x_6$};
\end{tikzpicture}    
  \caption{Star-graph with $N=6$ edges}
  \label{fig:N_star_graph}
\end{figure}
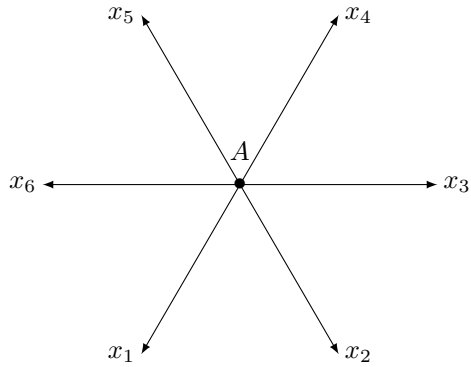

In this section, we are interested in the computation of ground state solutions for a general
$N$-edges star-graph $\mathcal{G}$ with a central vertex denoted by $A$ with a 
$\delta$ vertex condition at $A$. 
Each edge will be numbered with a
label $i=1,\dots,N$ (see Figure~\ref{fig:N_star_graph}) and will be identified when necessary with the right half-line $[0,\infty)$.
The unknown $\psi$ is the collection of the functions
$\psi_i$  living on every edge: $\psi=(\psi_1,\cdots,\psi_N)^T$. The total
mass is defined by $M_N(\psi)=\sum_{i=1}^N
\int_{\mathbb{R}^+}|\psi_i(x_i)|^2\, dx_i$.

\begin{figure}[htbp!]
 \centering
\begin{subfigure}[htbp!]{0.49\textwidth}
  \centering
  \includegraphics[width=0.6\textwidth]{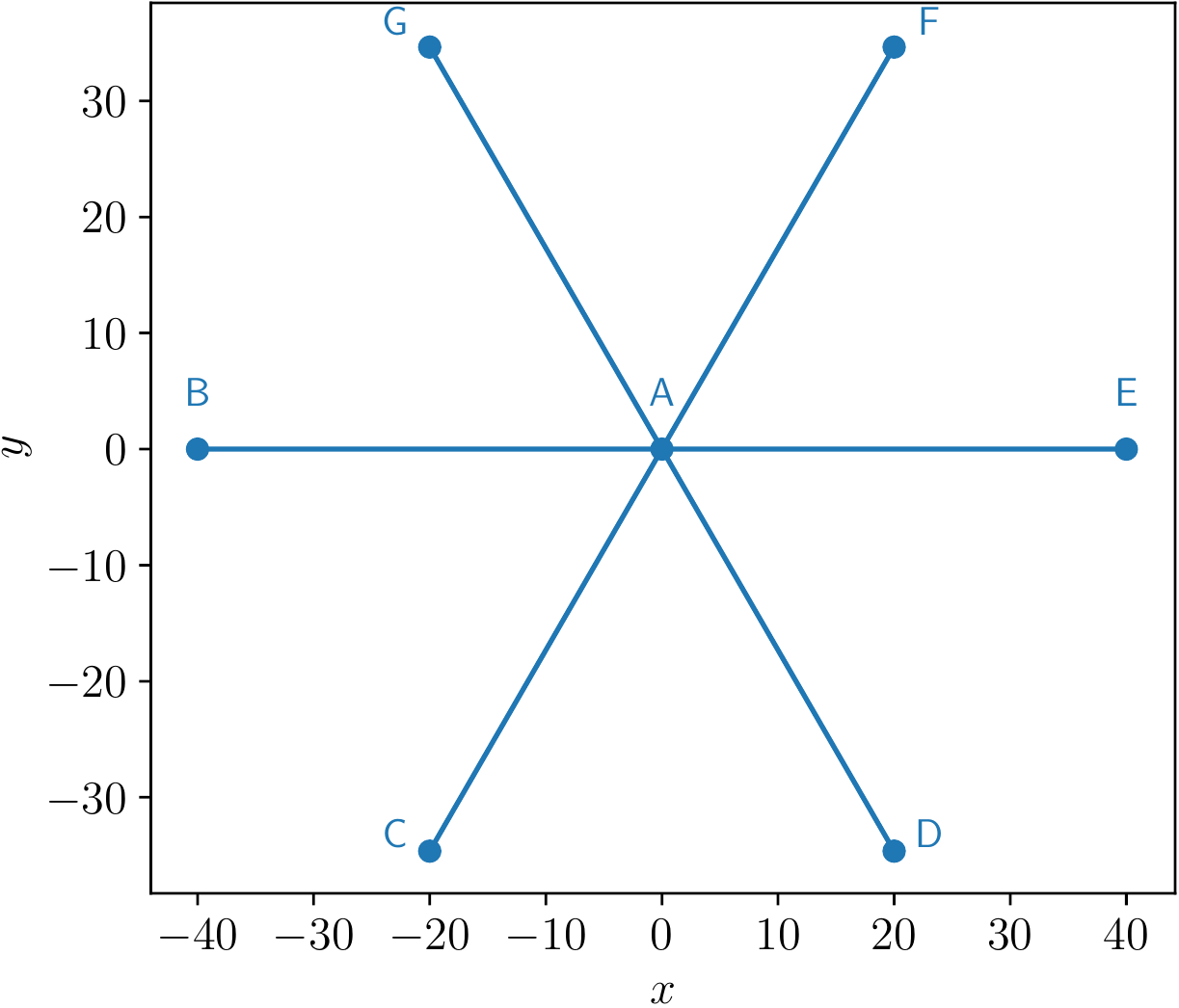}
  \caption{The approximated star-graph}
  \label{fig:approx_6graph}
\end{subfigure}
\begin{subfigure}[htbp!]{0.49\textwidth}
  \centering
  \includegraphics[width=0.6\textwidth]{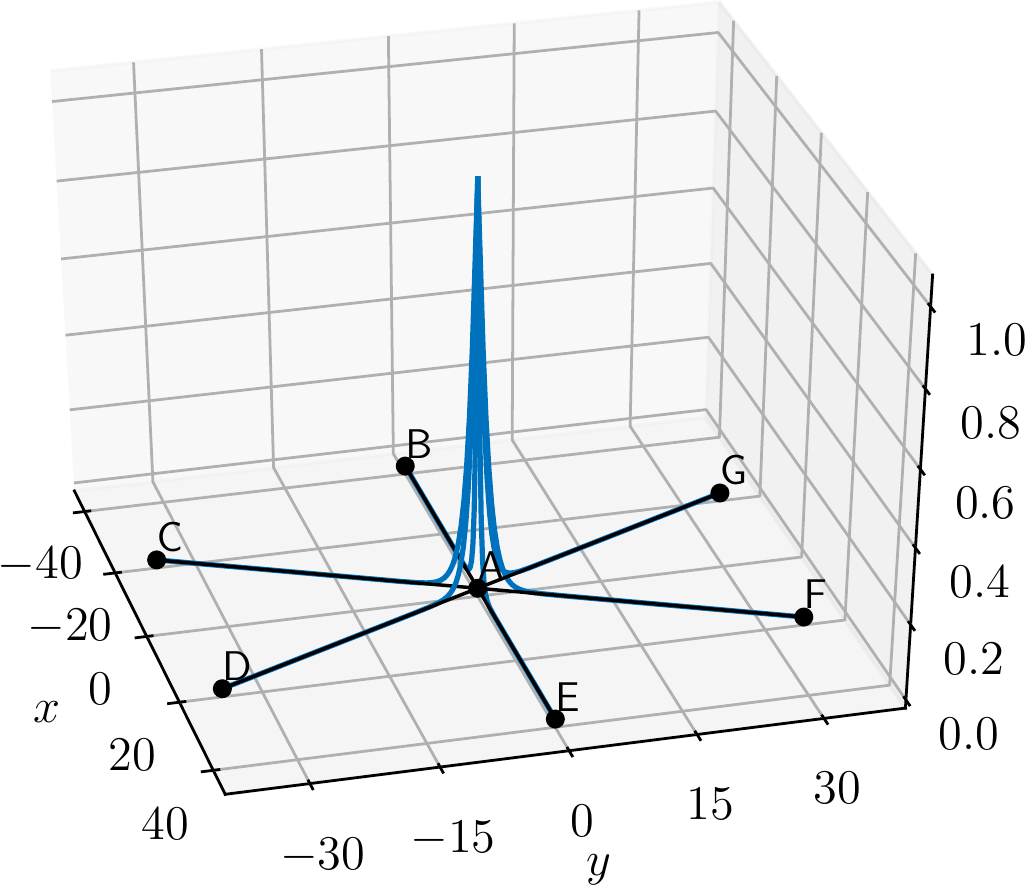}
  \caption{The numerical ground state when $\alpha=-4$ and $\omega=1$}
  \label{fig:num_sol_delta_6edges}
\end{subfigure}
\caption{Star
  graph with $6$ edges and $\delta$-condition}
\label{fig:6star-graph-num}
\end{figure}

The $\delta$ boundary conditions at $A$ are the generalization for $N>2$ of the $\delta$ potential on the line (i.e. the $2$-star graph, see e.g.~\cite{IaLeRo17,LeFuFiKsSi08} for studies of ground in this case):
\begin{equation*}
  \psi_j(0)=\psi_k(0)=:\psi(0),\quad 1\leq j,k\leq N,\qquad \sum_{j=1}^N\psi'_j(0)=\alpha
  \psi(0).
\end{equation*}
Ground states exist only for attractive $\delta$ potential, therefore we assume that
\[
  \alpha<0.
\]
We set $\lambda=1$. The energy is given by
\[
  E_{N,\delta}(\psi)=\sum_{i=1}^N\left[\frac12
  \int_{\mathbb{R}^+}|\psi'_i(x_i)|^2\, dx_i - \frac14 \int_{\mathbb{R}^+}|\psi_i(x_i)|^4\, dx_i\right]+\frac{\alpha}{2}|\psi_1(0)|^2.
\]
Let $m^*=4|\alpha|/N$. It was proved in~\cite{AdCaFiNo14a} that there exists a
ground state minimizing $E_{N,\delta}$ when $M_{N}=m$ if $m<m^*$ (there is no
constraint if $N=2$). The ground state is explicitly given in ~\cite{adami2012stationary}
and~\cite{AdCaFiNo14a} as follows. Let $\omega$ be implicitly given by
\[
m=2N\sqrt{\omega}-2\alpha.
  \]
Let $\bar{x}$ be defined by
\[
  \bar{x}=\frac{1}{\sqrt{\omega}} \text{arctanh}\left(\frac{|\alpha|}{N\sqrt{\omega}}\right).
\]
Then, the energy reaches its minimum when
$\psi=\psi_{\delta,\omega}$ (up to a phase factor)
where each component of $\psi_{\delta,\omega}$ is given by
\[
  {\psi}_{\delta,\omega,i}(x_i)=\frac{\sqrt{2\omega}}{\cosh(\sqrt{\omega}(x_i+\bar{x}))},\quad
  1\leq i\leq N,
\]
with $\omega \in (\alpha^2/N^2,+\infty)$. The mass of
$\psi_{\delta,\omega}$ is indeed 
\[
  M_N(\psi_{\delta,\omega})=2N\sqrt{\omega}-2\alpha=m,
\]
and its energy is given by
\[
  E_{N,\delta}(\psi_{\delta,\omega})=-\frac{N}{3}\omega^{3/2}-\frac{\alpha^3}{3N^2}.
\]
In order to compute numerically the ground state, each edge of the approximated
graph (see Figure~\ref{fig:6star-graph-num}~(\subref{fig:approx_6graph})) is of length $40$ and discretized with
$N_e=800$ nodes.
We add homogeneous Dirichlet boundary conditions at the terminal end of each edge. The
gradient step is $\delta t=10^{-2}$ and we perform $3000$ iterations. Each component
of the initial data $\psi_0$ is a Gaussian $\psi_{0,i}=\rho_i
e^{-10x_i^2}$ and $\rho_i$ is computed in such a way that the mass of
$\psi_0$ is $m$. We set
$\alpha=-4$ and $\omega=1$. The outcome is plotted on Figure~\ref{fig:6star-graph-num}~(\subref{fig:num_sol_delta_6edges}).

\begin{figure}[htbp!]
  \centering
  \begin{tabular}{cc}
    \includegraphics[width=.31\textwidth]{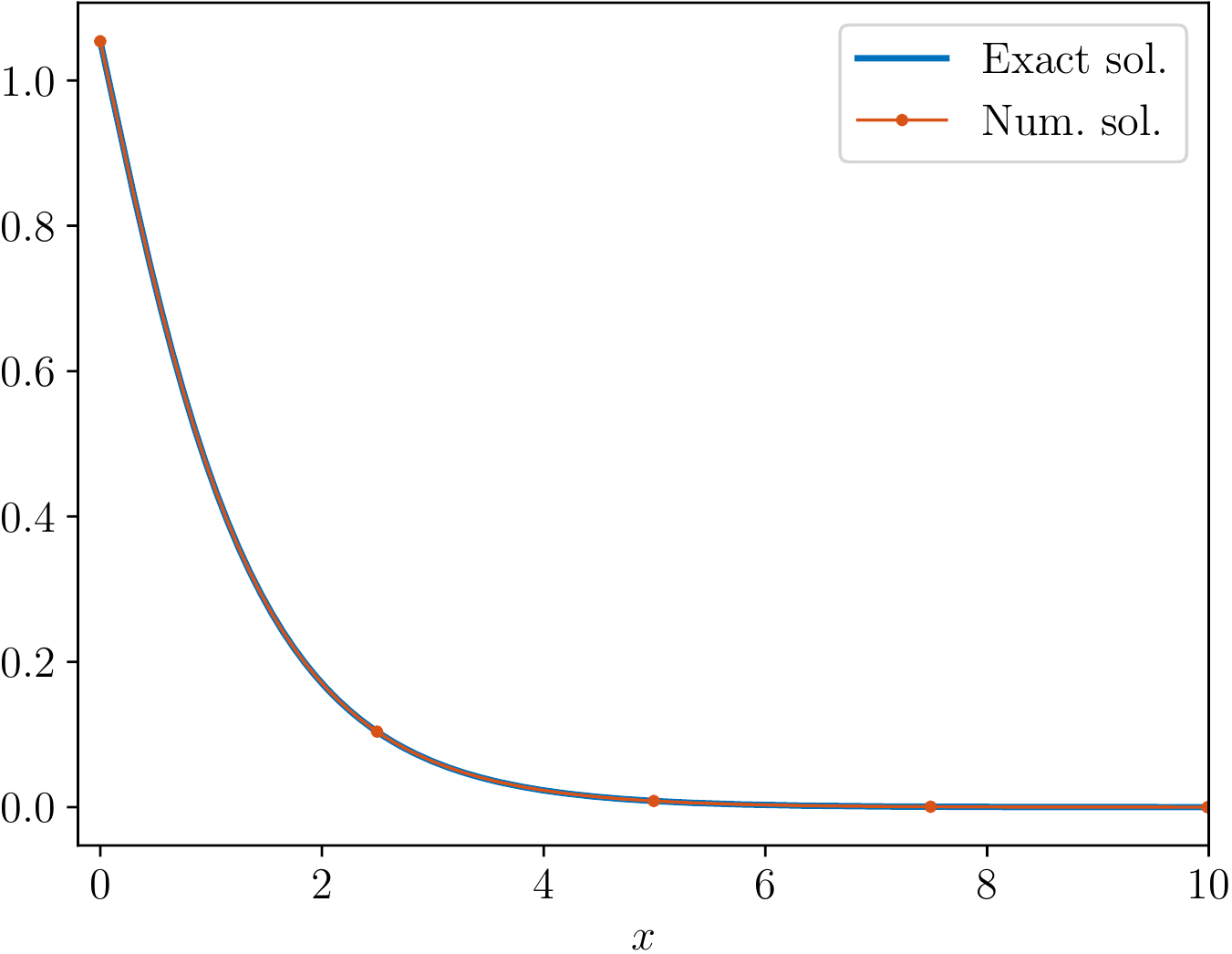} &
    \includegraphics[width=.31\textwidth]{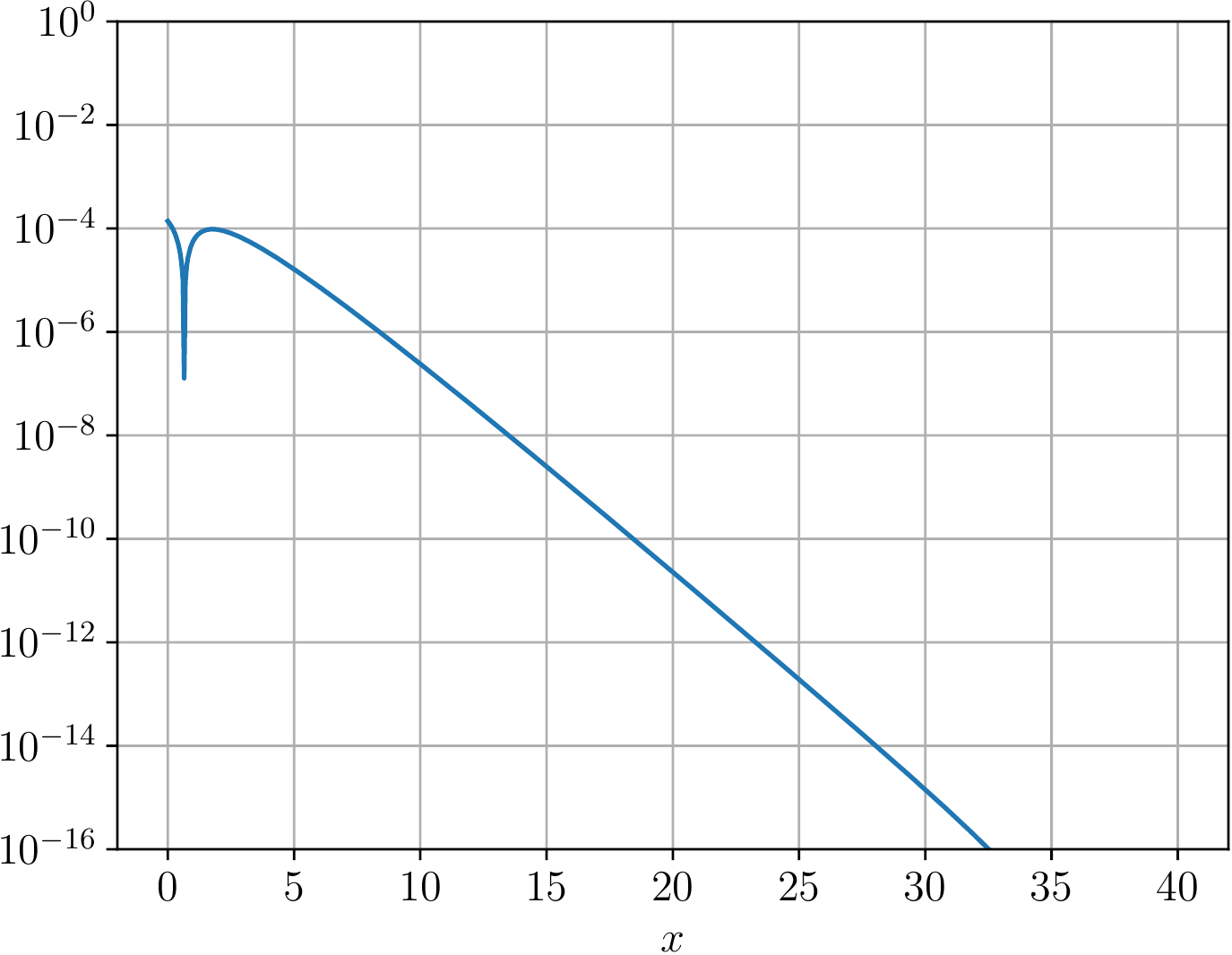}
  \end{tabular}
  \caption{Comparison between exact and numerical solutions for
    $\delta$-condition on a star-graph with $6$ edges (where the $\log10$ of the difference is depicted on the right).}
  \label{fig:comp_sol_delta_6edges}
\end{figure}

We plot on Figure~\ref{fig:comp_sol_delta_6edges} the comparison between the
exact solution and the numerical one on an edge (left) and the modulus of the
difference in log scale (right), thereby showing the very good agreement of our numerical computations with the theory.



\subsubsection{The tadpole
}
\label{sec:tadpole}

The classical tadpole graph consists of one loop with a half line attached to it and was considered in the subcritical case $1<p<5$ in~\cite{AdSeTi16,CaFiNo15,NoPeSh15}. The existence of a ground state for any given mass was established in~\cite[p 214]{AdSeTi16}, and the loop-centered bound state is the good candidate for the ground state. A classification of standing waves was performed in the cubic case $p=3$ by Cacciapuoti, Finco and Noja~\cite{CaFiNo15}, and was later extended to the whole subcritical range $1<p<5$ by Noja, Pelinovsky, and Shaikhova~\cite{NoPeSh15}, with some orbital stability results.

The generalized tadpole graph consists of one loop with $K$ half-lines attached at the same vertex (see e.g. Figure~\ref{fig:examples-tadpoles}) and was treated in~\cite{BeMaPe19}. When $K=2$, it is a particular case of the tower of bubbles on the line, with one bubble, and the ground state is known to be the soliton of the real line, folded on the bubble (see~\cite[Example 2.4]{AdSeTi15b}). For $K\geq 3$, there is no ground state (as Assumption~\ref{ass:H} is verified). 

Noja-Pelinovski~\cite{NoPe20} recently analyzed in details the standing waves on the tadpole graph for the critical quintic nonlinearity, with an alternative variational technique (minimization of the $H^1$ norm at fixed $L^6$-norm). In particular, they established the existence of a branch of standing waves for which three regimes exist, depending of the frequency $\omega$ of the wave. There exist $0<\omega_0<\omega_1$ such that standing waves are ground states if $0<\omega<\omega_0$, local minimizers of the energy at fixed mass if $\omega_0<\omega<\omega_1$, and saddle points for the energy at fixed mass if $\omega>\omega_1$.

\tikzset{solid node/.style={circle,draw,inner sep=1,fill=black}}
\begin{figure}[htbp!]
  \centering
  \begin{subfigure}[t]{0.31\textwidth}
      \centering
  \begin{tikzpicture}
      \draw (-0.5,0) circle (0.5);
    \node[solid node] at (0,0) {};
      \draw (0,0) -- (1,0);
      \draw[dotted] (1,0) -- (1.5,0);
    \end{tikzpicture}
  \caption{Classical tadpole graph}
  \label{fig:tadpole-1}
\end{subfigure}
  \begin{subfigure}[t]{0.31\textwidth}
  \centering
  \begin{tikzpicture}
      \draw (-0.5,0) circle (0.5);
    \node[solid node] at (0,0) {};
      \draw[rotate=15] (0,0) -- (1,0);
      \draw[rotate around={15:(0,0)},dotted] (1,0) -- (1.5,0);
      \draw[rotate=-15] (0,0) -- (1,0);
      \draw[rotate around={-15:(0,0)},dotted] (1,0) -- (1.5,0);
    \end{tikzpicture}
  \caption{Generalized tadpole graph with $2$ branches}
  \label{fig:tadpole-2}
\end{subfigure}
  \begin{subfigure}[t]{0.31\textwidth}
  \centering
  \begin{tikzpicture}
      \draw (-0.5,0) circle (0.5);
    \node[solid node] at (0,0) {};
\draw (0,0) -- (1,0);
      \draw[dotted] (1,0) -- (1.5,0);
      \draw[rotate=30] (0,0) -- (1,0);
      \draw[rotate around={30:(0,0)},dotted] (1,0) -- (1.5,0);
      \draw[rotate=-30] (0,0) -- (1,0);
      \draw[rotate around={-30:(0,0)},dotted] (1,0) -- (1.5,0);
    \end{tikzpicture}
  \caption{Generalized tadpole graph with $3$ branches}
  \label{fig:tadpole-3}
\end{subfigure}
\caption{Examples of tadpole graphs}
\label{fig:examples-tadpoles}
\end{figure}
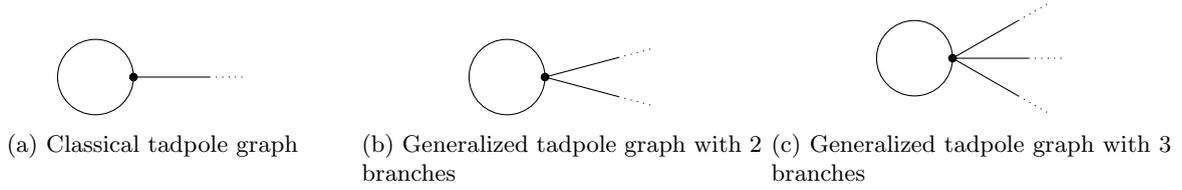

\begin{figure}[htbp!]
  \centering
  \begin{subfigure}[t]{0.49\textwidth}
  \centering
    \includegraphics[width=.9\textwidth]{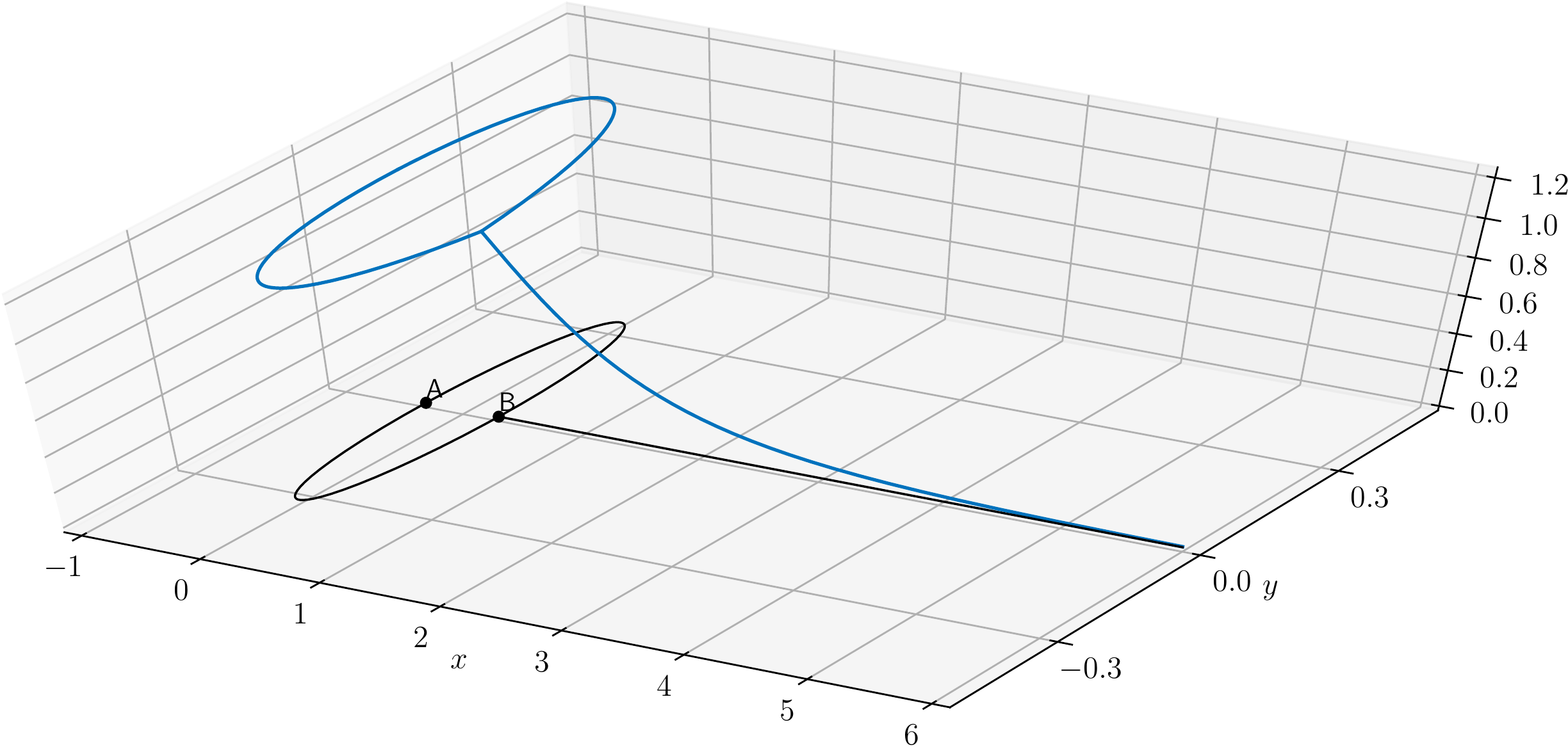}
  \caption{The numerical ground state}
  \label{fig:sol_tapole}
\end{subfigure}
  \begin{subfigure}[t]{0.49\textwidth}
  \centering
    \includegraphics[width=.9\textwidth]{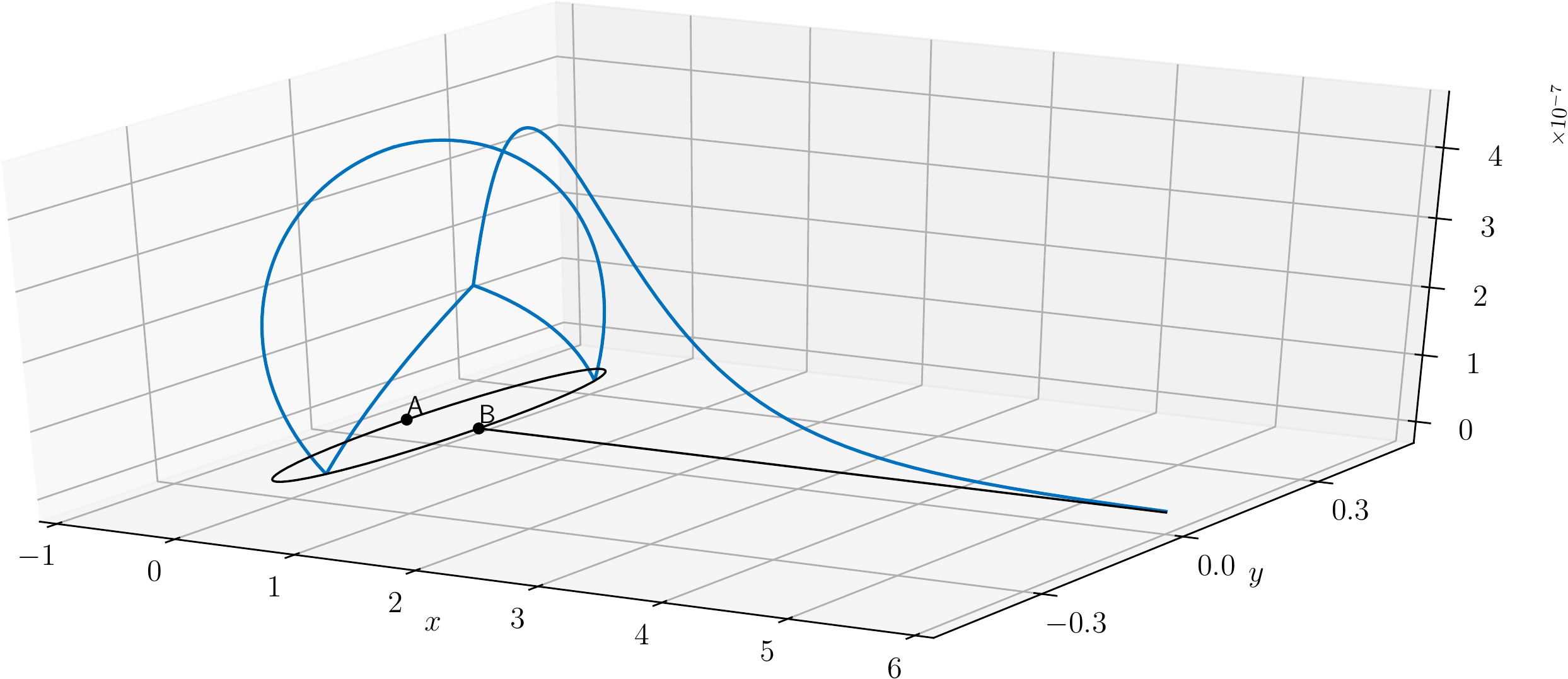}
  \caption{$|\psi_{\text{tadpole}} -
    \psi_{\text{num}}|$}
  \label{fig:err_tadpole}
\end{subfigure}
\caption{On the tadpole graph}
\label{fig:err-and-num-tadpole}
\end{figure}

In this section, we present the computation of the ground state to the NLS equation~\eqref{eq:nls2} with $\lambda=1$ on a classical tadpole graph. The graph $\mathcal{G}$ is made of a ring of perimeter
$2L$ and a semi-infinite line (tail) originated from a vertex with Kirchhoff condition. 
It is conjectured in~\cite{CaFiNo15} that the ground state
exists and is made of a dnoidal-type function on the ring and a sech-type
function on the tail. Its explicit formula on the ring is
\[
  \psi_{\text{ring}}(x)=\sqrt{\frac{2\omega}{2-k_*^2}} \dn\left(
    \sqrt{\frac{\omega}{2-k_*^2}}x;k^*\right),\quad 0< k_*<1,
\]
where $\dn$ is given by~\eqref{eq:dn} and $k_*\in(0,1)$ is the solution of
\[
  \frac{3k^4}{1-k^2} \text{cn}^2\left(\frac{L\sqrt{\omega}}{\sqrt{2-k^2}};k\right)\left[1-\text{cn}^2\left(\frac{L\sqrt{\omega}}{\sqrt{2-k^2}};k\right)\right]=1,
\]
where $\text{cn}$ denotes the cnoidal function defined in~\eqref{eq:sn_cn}.
The solution on the tail is
\[
  \psi_{\text{tail}}(x)=\frac{\sqrt{2\omega}}{\cosh(\sqrt{\omega}(x-b))},
\]
where $b$ is determined by the negative solution of
\[
  \frac{1}{\cosh^2(\sqrt{\omega} b)}=\frac{\psi_{\text{ring}}^2(L)}{2\omega}.
\]
We take a ring of radius $1/\pi$, so that $L=1$, and we approximate the tail by a segment of length
$30$. We add homogeneous Dirichlet boundary conditions at the terminal vertex. We
take $\omega=1$. With these quantities, the couple $(k_*,b)$ is given by
\[
  k_*=.81664827149276692790,\quad b=.89507479534736339894.
\]
The mass of the ground state
$\psi_{\text{tadpole}}=(\psi_{\text{ring}},\psi_{\text{tail}})$ is
\[M(\psi_{\text{tadpole}})=3.1727382562292.\]
The numerical solution
is plotted in Figure~\ref{fig:err-and-num-tadpole}~(\subref{fig:sol_tapole}). We also plot the difference in
absolute value between $\psi_{\text{tadpole}}$ and the numerical
solution on Figure~\ref{fig:err-and-num-tadpole}~(\subref{fig:err_tadpole}). The maximum value of the error is
$4.45\cdot 10^{-7}$.

\subsection{Periodic graphs}
Periodic graphs are graphs with an infinite number of (usually finite length) edges, for which an elementary structure, the periodicity cell, is repeated in one or more directions.

\begin{figure}[htbp!]
  \centering
  \begin{tikzpicture}
          \foreach \x in {0,2,...,6}
          {
            \node at (\x,0)  {$\bullet$};
            \node at (\x,1)  {$\bullet$};
            \draw (\x,0)--(\x,1);
             \draw (\x-1,0)--(\x+1,0);
      \draw (\x,1.5) circle (0.5);
    }
                 \draw[dotted] (7,0)--(7.5,0);
                 \draw[dotted] (-1.5,0)--(-1,0);
	\end{tikzpicture}
	\caption{A periodic graph not satisfying Assumption (H$_{per}$)}
	\label{fig:noHper}
      \end{figure}
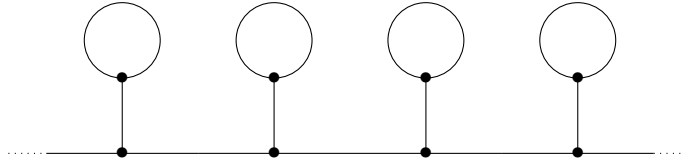

In the case of $1$-d periodic graphs (i.e. graphs for which the periodicity cell is copied in only one direction), Dovetta~\cite{Do19} proved that the situation is similar to the one of the real line: for $1<p<5$, there exists a ground state for every mass. The critical case $p=5$ is a bit more complicated. On one hand, if the graph satisfies the equivalent of the topological Assumption (H) adapted to the periodic setting (Assumption (H$_{per}$)), Dovetta~\cite{Do19} proved the non-existence of ground states. On the other hand, for graphs violating this topological assumption (see for example Figure~\ref{fig:noHper}), there may exist a whole interval of mass for which a ground state exists.

      In a somewhat different framework (including in particular periodic potentials in the problem), Pankov~\cite{Pa18} proved, under a spectral assumption on the underlying quantum graph, the existence of localized and periodic standing wave profile solutions. These profiles are obtained by minimizing the action (which in our case corresponds to $E+\omega M$ for a fixed $\omega$) on the corresponding Nehari manifold, but, as usual, it is unclear how and in which case these profiles could also be minimizers of the energy at fixed $L^2$ norm (recall that in the case of the real line minimizers are obtained on the Nehari manifold for any $1<p<\infty$, whereas on the mass constraint they exist only if $1<p<5$).

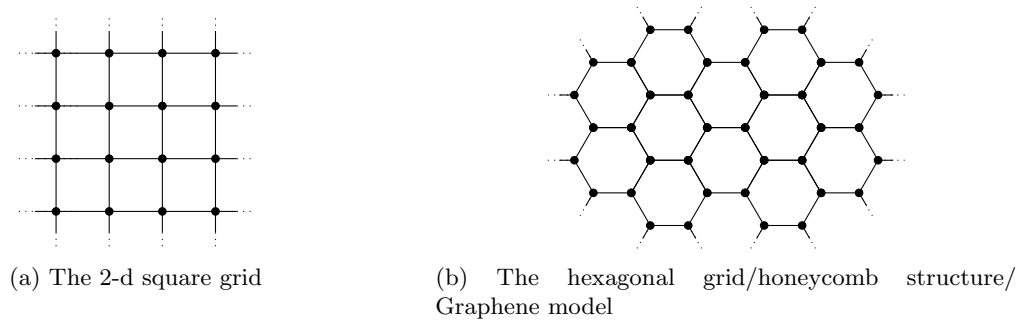
\begin{figure}
  \centering
                \tikzset{solid node/.style={circle,draw,inner sep=1,fill=black}}
    \begin{subfigure}[t]{0.45\textwidth}
  \centering
	\begin{tikzpicture}[scale=0.7]
	\draw (0,0) grid (3,3);
	\foreach \x in {0,1,...,3} \foreach \y in {0,1,...,3} \node[solid node] at (\x,\y)  {};
	\foreach \x in {0,1,...,3} 
	{\draw (\x,0)--(\x,-0.4);
          \draw[dotted] (\x,-0.4)--(\x,-0.7);
          \draw (\x,3)--(\x,3.4);
                    \draw[dotted] (\x,3.4)--(\x,3.7);
        }
        	\foreach \y in {0,1,...,3} 
                {\draw (-0.4,\y)--(0,\y);
                  \draw[dotted] (-0.7,\y)--(0.4,\y);
          \draw (3,\y)--(3.4,\y);
          \draw[dotted] (3.4,\y)--(3.7,\y);
        }
	\end{tikzpicture}
	\caption{The $2$-d square grid}
	\label{fig:2d-grid}
      \end{subfigure}
    \begin{subfigure}[t]{0.45\textwidth}
	\centering
	\begin{tikzpicture}[scale=0.5]
          \foreach \x/\y in
          {
            (-1.5,0),(1.5,0),
            (0,sin{60}),(-3,sin{60}),(3,sin{60}),
            (-1.5,2*sin(60)),(1.5,2*sin(60)),
            (0,-sin{60}),(-3,-sin{60}),(3,-sin{60}),
            (-1.5,-2*sin(60)),(1.5,-2*sin(60))
          }
          {
            \draw[shift={(\x,\y)}]
            (-1,0) node[solid node] {}
            --(cos{120},sin{120}) node[solid node] {}
            --(cos{60},sin{60}) node[solid node] {}
            --(1,0) node[solid node] {}
            --(cos{-60},sin{-60}) node[solid node] {}
            --(cos{-120},sin{-120}) node[solid node] {}
            -- cycle;
          }
          \foreach \x/\y in {(-1,3*sin{60}),(2, 3*sin{60}), (3.5, 2*sin{60})}
          {
            \draw[shift={(\x,\y)}] (0,0)--(cos{60}*0.4,sin{60}*0.4);
            \draw[shift={(\x,\y)},dotted] (0,0)--(cos{60}*0.7,sin{60}*0.7);
          }
          \foreach \x/\y in {(-1,-3*sin{60}),(2, -3*sin{60}), (3.5, -2*sin{60})}
          {
            \draw[shift={(\x,\y)}] (0,0)--(cos{60}*0.4,-sin{60}*0.4);
            \draw[shift={(\x,\y)},dotted] (0,0)--(cos{60}*0.7,-sin{60}*0.7);
          }
          \foreach \x/\y in {(1,3*sin{60}),(-2, 3*sin{60}), (-3.5, 2*sin{60})}
          {
            \draw[shift={(\x,\y)}] (0,0)--(-cos{60}*0.4,sin{60}*0.4);
            \draw[shift={(\x,\y)},dotted] (0,0)--(-cos{60}*0.7,sin{60}*0.7);
          }
          \foreach \x/\y in {(1,-3*sin{60}),(-2, -3*sin{60}), (-3.5,- 2*sin{60})}
          {
            \draw[shift={(\x,\y)}] (0,0)--(-cos{60}*0.4,-sin{60}*0.4);
            \draw[shift={(\x,\y)},dotted] (0,0)--(-cos{60}*0.7,-sin{60}*0.7);
          }
          \foreach \x/\y in { (4,sin{60}),  (4,-sin{60})}
          {
            \draw[shift={(\x,\y)}] (0,0)--(0.4,0);
            \draw[shift={(\x,\y)},dotted] (0,0)--(0.7,0);
          }
          \foreach \x/\y in {        (-4,sin{60}),(-4,-sin{60})}
          {
            \draw[shift={(\x,\y)}] (0,0)--(-0.4,0);
            \draw[shift={(\x,\y)},dotted] (0,0)--(-0.7,0);
          }
       	\end{tikzpicture}
	\caption{The hexagonal grid/honeycomb structure/ Graphene model}
	\label{fig:honey}
      \end{subfigure}
      \caption{Doubly periodic metric graphs}
      \label{fig:doubly-periodic}
      \end{figure}

That graphs periodic along only one direction essentially mimic the behavior of the real line is somewhat expected. However, if the periodicity occur in more than one direction, a new dimensionality of the problem may appear (which was also absent for non-compact graphs with a finite number of edges). At the microscopic level, periodic graphs remain clearly $1$-d structures. But at the macroscopic level, periodic graphs may be seen as higher dimensional structures, for examples the $2$-d grid (see Figure~\ref{fig:doubly-periodic} (\subref{fig:2d-grid})) or the honeycomb hexagonal grid (see Figure~\ref{fig:doubly-periodic} (\subref{fig:honey})) are clearly $2$-d structures at the macroscopic level. This dimensional transition is reflected in the range of critical exponents and masses. Non-compact graphs with a finite number of edges share the same critical exponent (from a nonlinear Schr\"odinger point of view) as the line, i.e. the graphs are subcritical for power nonlinearities with exponents $1<p<5$, and minimization of the energy under a fixed mass constraint is possible only if $1<p\leq 5$. On the other hand, it was revealed in~\cite{AdDoRu19,AdDoSeTi19} that a \emph{dimensional crossover} with a continuum of critical exponents occurs for the $2$-d grid and the hexagonal grid. More precisely, the following has been established in~\cite{AdDoRu19,AdDoSeTi19}. If $1<p<3$, then there exists a ground state for any possible value of the mass. If $3\leq p<5$, then there exists a critical value $m_c$ of the mass such that ground states exist if and only if $m\geq m_c$ (unless $p=3$, in which case the case $m=m_c$ is open). If $p\geq 5$, then a ground state never exists, no matter the value of the mass. Recall that $3$ (resp. $5$) is the critical exponent for the nonlinear Schr\"odinger equation on $\R^2$ (resp. on $\R$). Similar results have been obtained for the $3$-d grid by Adami and Dovetta~\cite{AdDo18}.

In what follows, we present some numerical experiments realized in two model cases: the necklace and the hexagonal grid.

\subsubsection{The necklace}

\begin{figure}[htbp!]
  \centering
  \begin{tikzpicture}
      \node at (-2.5,0) {$\bullet$};
      \node at (-1.5,0) {$\bullet$};
      \node at (-0.5,0) {$\bullet$};
      \node at (0.5,0) {$\bullet$};
      \node at (1.5,0) {$\bullet$};
      \node at (2.5,0) {$\bullet$};
      \draw (-2,0) circle (0.5);
      \draw (0,0) circle (0.5);
      \draw (2,0) circle (0.5);
      \draw (-1.5,0) -- (-0.5,0);
      \draw (0.5,0) -- (1.5,0);
      \draw (-3,0) -- (-2.5,0);
      \draw (2.5,0) -- (3,0);
      \draw[dotted] (-3.5,0) -- (-3,0);
      \draw[dotted] (3,0) -- (3.5,0);
    \end{tikzpicture}
  \caption{The necklace graph, a periodic graph with alternating loop and single edge}
  \label{fig:necklace}
\end{figure}
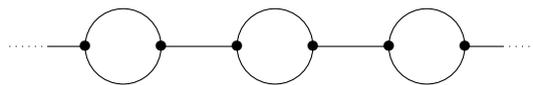

The necklace graph is a periodic graph consisting of a series of loops alternating with single edges (see Figure~\ref{fig:necklace}) and is probably one of the simplest non-trivial periodic graphs. 
The validity of the NLS approximation for periodic quantum graphs of necklace type was established by Gilg, Pelinovsky and Schneider~\cite{GiPeSc16}. Moreover, Pelinovsky and Schneider~\cite{PeSc17} showed the existence, at fixed sufficiently small frequency $\omega$, of two symmetric positive exponentially decaying bound states, one located at the center of the single edge and the other equally distributed with respect to the centers of each half-loop. It is conjectured in~\cite{PeSc17} that the state located on the single edge should be the ground state at small mass. On the other hand, for large masses,   
it was experimentally observed in~\cite{BeMaPe19} that their estimates on edge localized bound states could also be applied in the case of the necklace graph. The conclusion of this observation is that at large mass the ground state should be centered on the loop if the length of the internal edge is smaller that the length of the half-loop, and vice versa.

We have performed numerical calculations of the ground states on a necklace consisting of  loops of total length $\pi$ (i.e. each branch of the loop is of length $\pi/2$) and connecting edges of length $1$. The length of the necklace is chosen to be large, but obviously necessarily finite. In practice, the length needs to be adapted depending on the mass $m$ on which we are minimizing the Schr\"odinger energy. Indeed, it is expected (and appears to be so in practice) that the ground state will be decaying as $e^{-m|x|}$ from some central point on the graph (here, $|x|$ is referring to the (graph) distance with respect to this point). Therefore, the smaller the mass is, the larger the length of the necklace needs to be in order to fully capture the tail of the ground state. The conditions at the vertices are Kirchhoff conditions, apart from the end points where we have chosen to set Dirichlet conditions.

\begin{figure}[htbp!]
  \centering
  \begin{subfigure}[t]{0.31\textwidth}
    \centering
    \includegraphics[width=\textwidth]{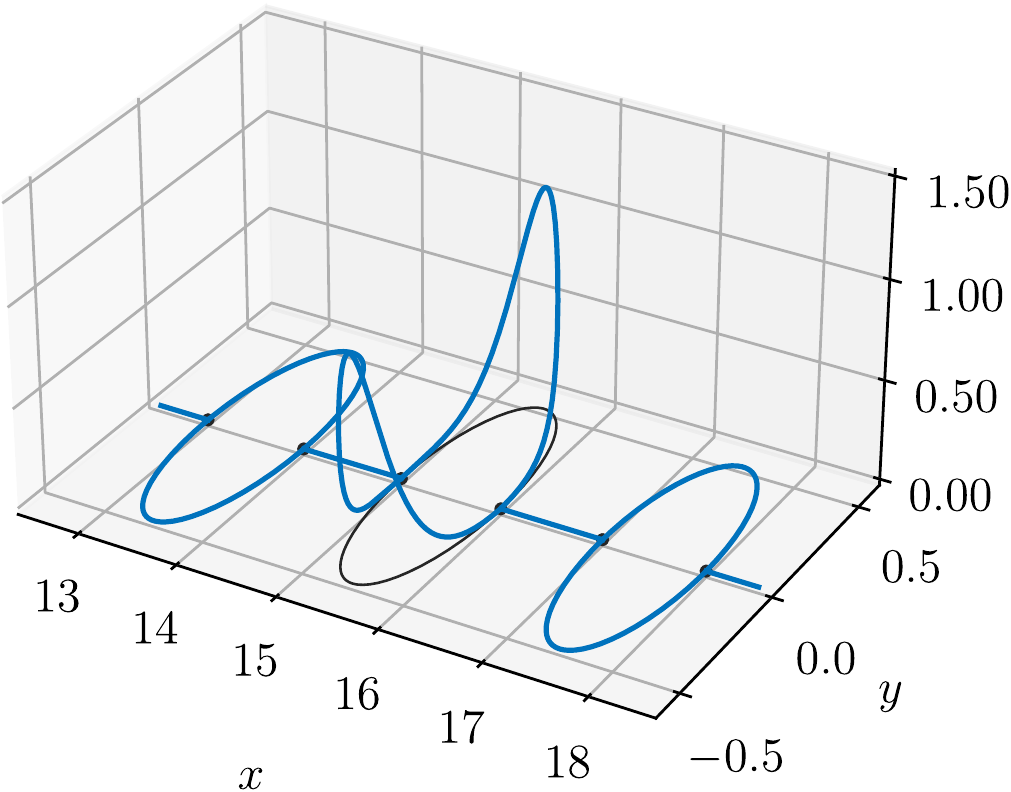}
    \caption{$Init\;2$}
    \label{fig:init2}
  \end{subfigure}
  \begin{subfigure}[t]{0.31\textwidth}
    \centering
    \includegraphics[width=\textwidth]{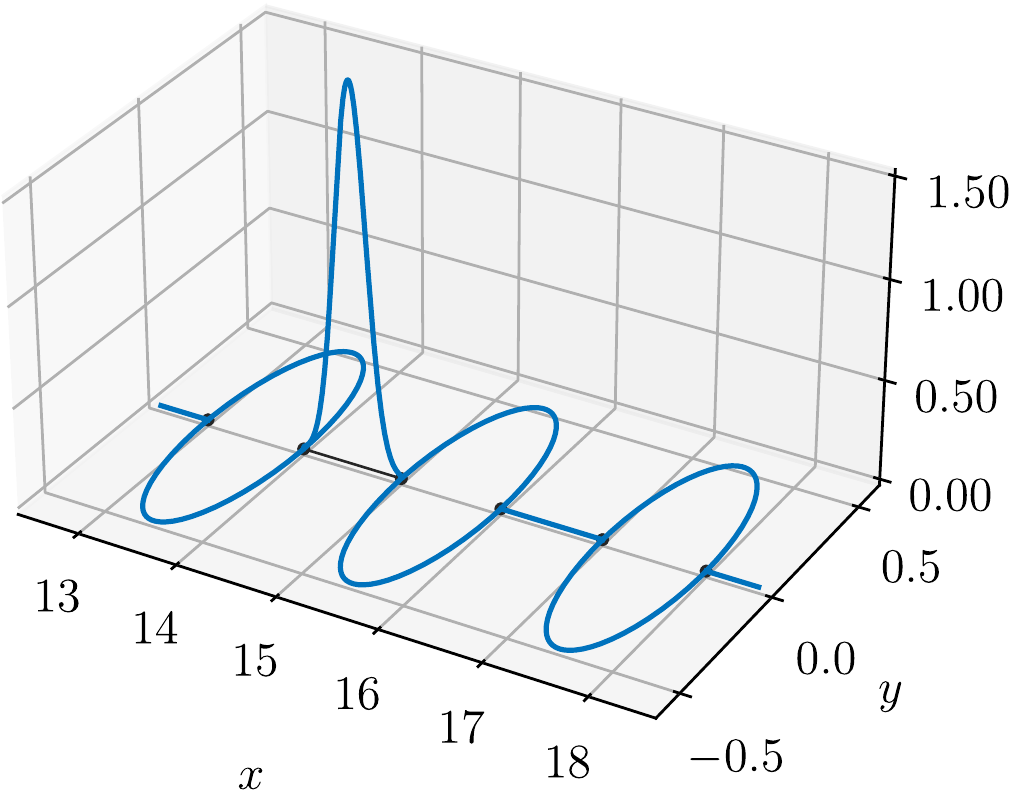}
    \caption{$Init\;3$}
    \label{fig:init3}
  \end{subfigure}
  \begin{subfigure}[t]{0.31\textwidth}
    \centering
    \includegraphics[width=\textwidth]{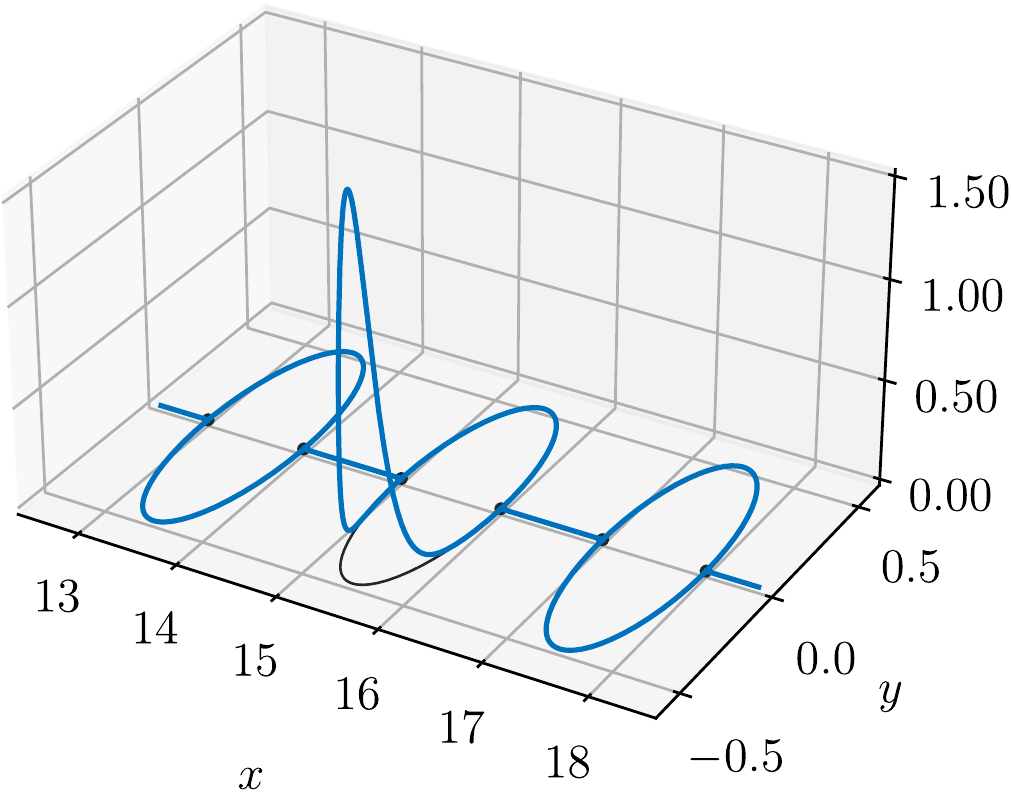}\\
    \caption{$Init\;4$}
    \label{fig:init4}
  \end{subfigure}
  \caption{The three initial data constructed with Gaussians of the
    form $Ce^{-10x^2}$, centered on edges, truncated at the end
    points, and with $C$ adjusted to satisfy the mass constraint.}
  \label{fig:necklace-init}
\end{figure}

We have chosen to perform a collection of experiments for masses
varying from very small to very large and with three different types
of initial data, all positioned on the periodicity cell at the middle
of the necklace: two gaussians concentred and centered  on each of the branches of the loop
(referred to as $Init\;2$, see Figure~\ref{fig:necklace-init} (\subref{fig:init2})), a
gaussian concentred and centered on the single connecting edge  (referred to as $Init\;3$, see Figure~\ref{fig:necklace-init} (\subref{fig:init3})), and 
a gaussian concentred and centered  on a branch of
the circle (referred to as $Init\;4$, see Figure~\ref{fig:necklace-init} (\subref{fig:init4})).

\begin{figure}[htpb!]
    \includegraphics[width=\textwidth]{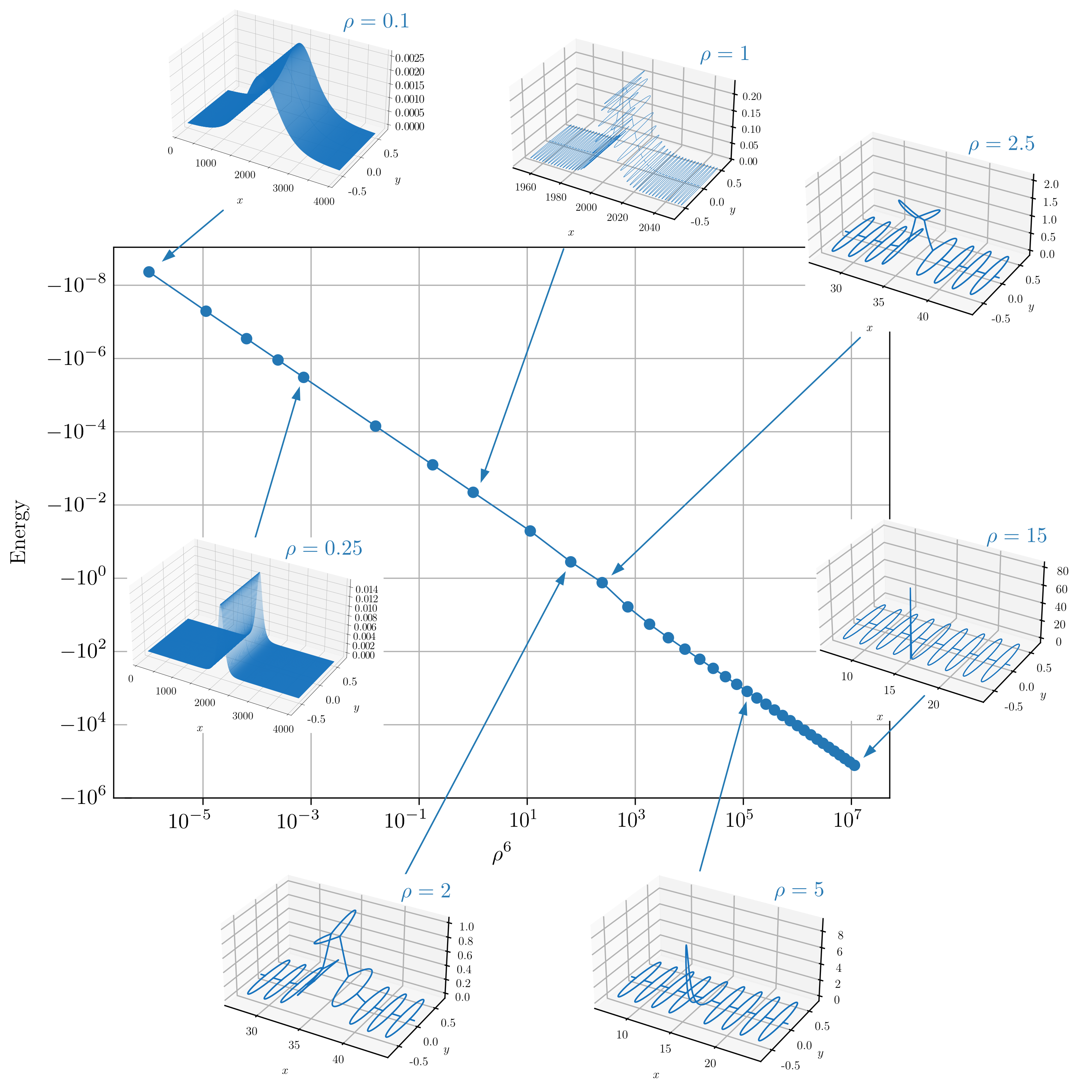} 
  \caption{Mass-Energy diagram with sample representations for the
    ground states on a necklace graph}
  \label{fig:necklace-diagram}
\end{figure}

We first present a global picture (see Figure~\ref{fig:necklace-diagram}) of the ground states for 
$L^2$ norms ranging from $\rho=0.1$ to $\rho=15$ (recall that $\rho=\sqrt{m}$). Since we expected the energy to be of order $m^3$, we have presented the mass-energy with $\rho^6=m^3$ log-scale on the horizontal axis.
Our expectation is confirmed by the representation which is indeed a straight line, with a slight shift around $\rho=2.5$ corresponding to a bifurcation (on which we will comment after).  We observe that for small masses, the ground state is scattered across many periodicity cells. As the mass increases, the ground state becomes more and more concentrated on a loop, first symmetrically on both branches of the loop, then on only one branch of the loop.

  \begin{figure}[htpb!]
    \includegraphics[width=0.5\textwidth]{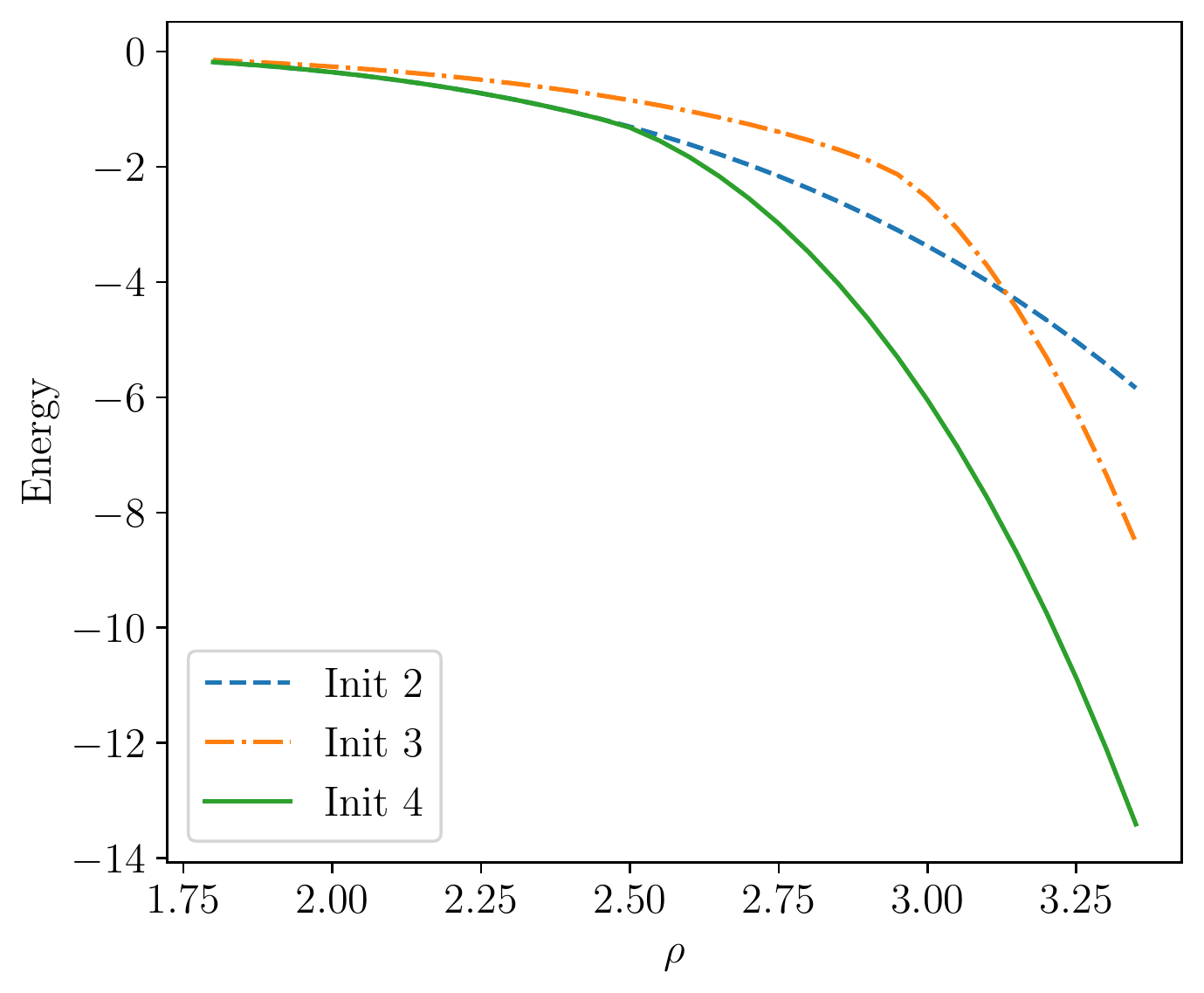} 
  \caption{Mass-Energy diagram for local minimizers computed with $Init\;2$, $Init\;2$, and $Init\;4$}
  \label{fig:necklace-bifurcation}
\end{figure}
Figure~\ref{fig:necklace-diagram} was devoted to the ground state. In fact, we may perform a
  more detailed analysis and obtain other branches of local minimizers of the
  energy at fixed mass. Indeed, provided the parameters of our algorithm are
  suitably chosen, starting from each of the initial data $Init\;j$, $j=1,2,3$, we should have
  convergence towards the closest local minimizer of the energy at fixed
  mass. The obtained minimizer should itself enjoy similar features as the
  initial data (e.g. the place of centering). We present the outcome of our
  simulations in Figure~\ref{fig:necklace-bifurcation}. Each initial data gives rise to a branch of
  local minimizers. For small mass, the branches corresponding to $Init\;2$ and
  $Init\;4$ coincide and correspond to the ground state, which is centered on a loop
  and symmetric with respect to both sides of the loop. At $\rho\simeq 2.5$, we observe a bifurcation and the branches corresponding to $Init\;2$ and $Init\;4$ separate, as the $Init\;4$ branch bifurcates with smaller energy and is formed of ground states peaked on one side of a loop, whereas the $Init\;4$ branch continues the branch of symmetric states on a loop (which are not anymore ground states). The $Init\;3$ branch is formed all along of states centered on a single edge. It is never a ground state branch, but it is meeting the $Init\;4$ branch at small and large mass, up to a point where they become indistinguishable numerically (for large mass, outside of Figure~\ref{fig:necklace-bifurcation}, at $\rho\simeq 10$).


\subsubsection{The honeycomb}

We now turn to the honeycomb grid. This is a graph which is built recursively using a hexagonal tessellation. As the necklace graph, it is a very simple periodic graph and we can see that it is two-dimensional on a large scale. In~\cite{AdDoRu19}, the existence of minimizers for the NLS energy functional is proved for $1<p<3$, for any mass. Here, we perform some numerical simulations in the case $p = 2$. To be more specific, we use the gradient methods to compute the ground state of the NLS energy functional under a specified mass. As noted in~\cite{AdDoRu19,AdDoSeTi19}, for low masses, we expect the ground state to display a $2$d structure due to the spreading on the graph. For large masses, on the contrary, the ground state should be more localized on the graph and we expect a $1$d structure. The goal of this numerical investigation is to describe the transition from the $2$d regime to the $1$d regime by varying the mass of the ground state from $1$ to $16$.

\begin{figure}[htpb!]
\includegraphics[width=0.8\textwidth]{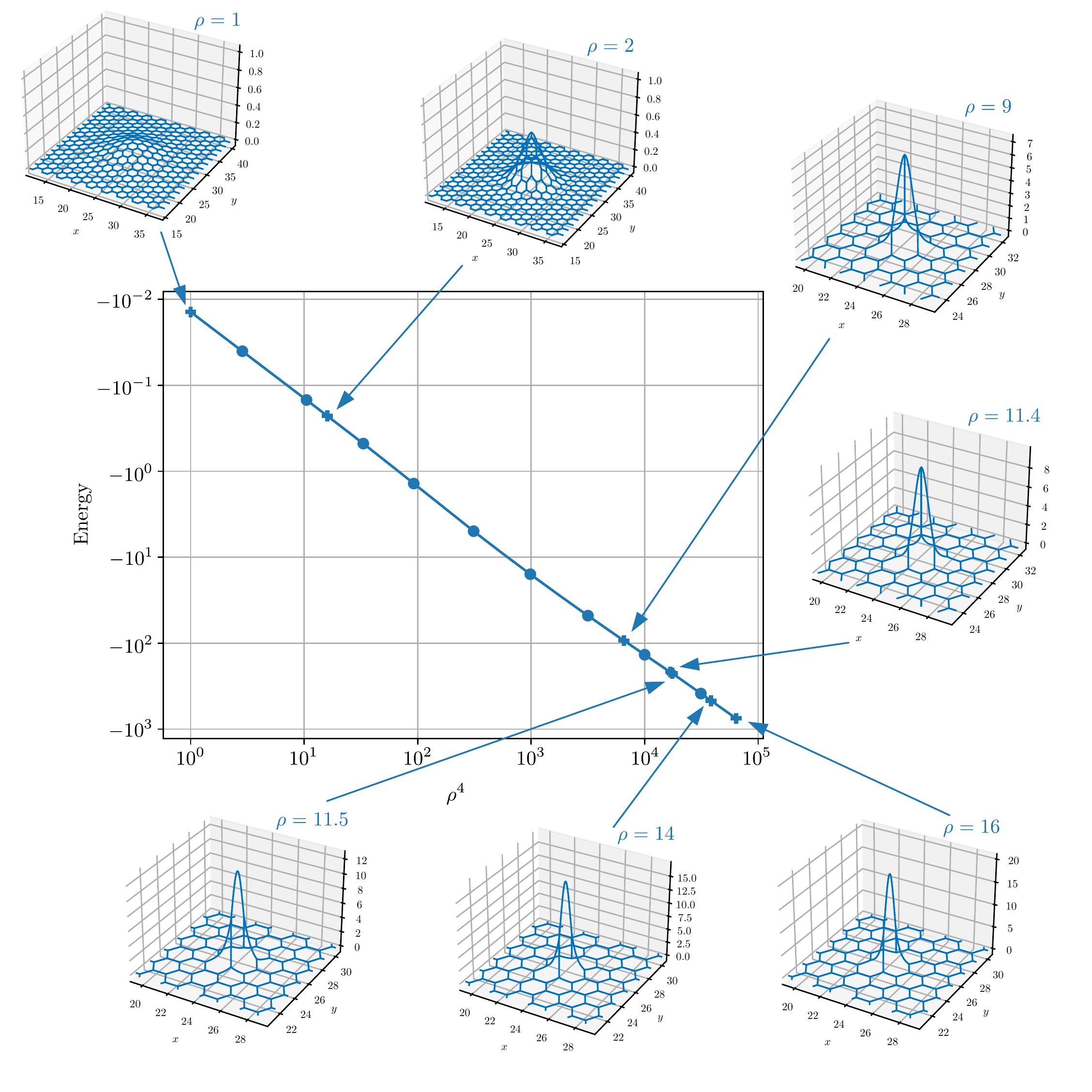} 
\caption{Mass-Energy diagram for the minimizers of NLS energy functional on a honeycomb grid.}
\label{fig:hexaglobal}
\end{figure}

The graph is set such that each edge has a length of $1$. We have obtained the Mass-Energy diagram which is depicted in Figure~\ref{fig:hexaglobal}. To begin with, we note that there is a linear relation between the energy and $\rho^4$. We can see that, for low masses, the ground state looks like a $2$d ground state in the Euclidean case. Furthermore, we remark that it is centered on a node (that is, its maximum is located on a node) and symmetric. As the mass grows larger, the ground state is more concentrated. Then, between a mass of $11.4$ and $11.5$, we observe a structural transition: the minimizer becomes centered on an edge (still symmetric). For larger masses, it keeps concentrating (slowly) on a single edge and, thus, it displays a $1$d regime.

\subsection{Metric Trees}

\begin{figure}[htpb!]
  \centering
  \begin{tabular}{cc}
    \includegraphics[width=.25\textwidth]{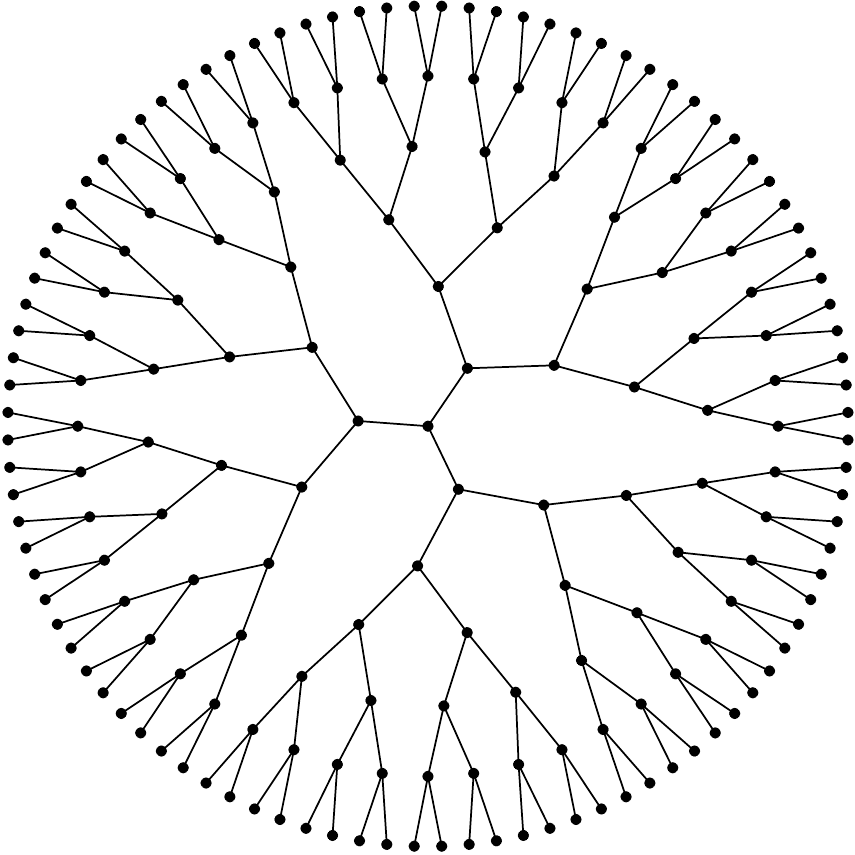}
    &     \includegraphics[width=.31\textwidth]{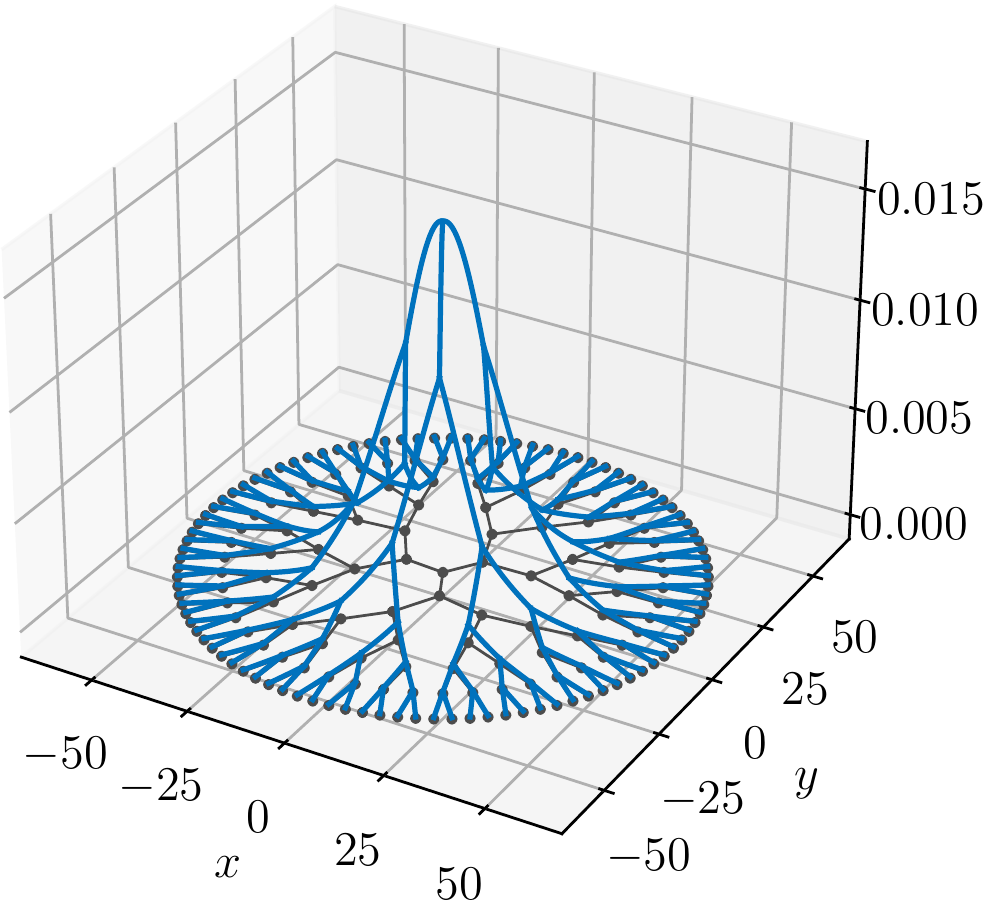}\\
  \end{tabular}
  \caption{Metric tree (left) and ground state with mass constraint $\rho=0.1$ (right) for $|u|u$ nonlinearity.}
  \label{fig:GS_on_tree}
\end{figure}

Metric trees are tree-type graphs endowed with a metric structure. In this section, we are interested in the case of binary trees, i.e. trees for which each vertex (except for the root, if any) has degree $3$ and all the edges share the same length. Dispersion of the Schr\"odinger group on trees (with $\delta$ conditions at the vertices) was investigated by Banica and Ignat~\cite{BaIg14}. Existence of ground states on metric trees (with Kirchhoff conditions at the vertices) has been considered by Dovetta, Serra and Tilli~\cite{DoSeTi20} in the case of binary trees, either rooted or non-rooted. Let $\mathcal G$ be a rooted or non-rooted binary tree with Kirchhoff vertices conditions and (following the notation of~\cite{DoSeTi20}), define the minimum of the Schr\"odinger energy at fixed mass $m$ by
\[
  \mathcal L_{\mathcal G}(m)=\min\{E(u):u\in H^1_D(\mathcal G),\;M(u)=m\}.
\]
It was proved in~\cite{DoSeTi20} that there exists a critical mass $m_{\mathcal G}^*\geq 0$ such that
      \[
        \left\{
        \begin{array}{lll}
          \mathcal L_{\mathcal G}(m)=\frac12\lambda_1m,&\text{ and there is no ground state,}&\text{ if }\mu\in(0,m_{\mathcal G}^*),\\
                    \mathcal L_{\mathcal G}(m)<\frac12\lambda_1m,&\text{ and a ground state exists,}&\text{ if }m>m_{\mathcal G}^*,
        \end{array}
      \right.
     \] 
     where $\lambda_1$ is the optimal constant for the Poincar\'e inequality on the graph.
The nonlinearity considered in~\cite{DoSeTi20} is any mass-subcritical power nonlinearity, i.e. $|u|^{p-1}u$ with $1<p<5$. 
If $3<p<5$ or if $1<p<5$ and minimization is done in the class of radially symmetric functions, the authors of~\cite{DoSeTi20} proved that $m_{\mathcal G}^*>0$. The case $1<p<3$ is open if no symmetry assumption is made, but the authors conjecture that  $m_{\mathcal G}^*$ is also positive in this case. Moreover, they conjecture that minimizers should be radial even when no symmetry assumption is made on the class of function in which minimization is done. This is confirmed by experiments that we conducted on a binary tree of depth $6$ with each branch of approximate length $10$ (we have arranged the vertices in such a way that they are on concentric circles). We give a sample result of our experiments in Figure~\ref{fig:GS_on_tree}.

\appendix
\section{Features of the Grafidi library}

The Grafidi library relies on the following Python libraries:
Matplotlib~\cite{Hunter:2007},
Networkx~\cite{hagberg-2008-exploring},
Numpy~\cite{harris2020array},
Scipy~\cite{2020SciPy-NMeth}.

\subsection{Methods from the Graph class}
\subsubsection{\textunderscore\textunderscore init\textunderscore\textunderscore}

The Graph class constructor builds an instance of a graph which is based on the graphs from the NetworkX library (Network Analysis in Python).

\vspace{0.5em}
\noindent\textbf{g = Graph(g\textunderscore nx,Np,user\textunderscore bc)}

\vspace{0.5em}
\textit{Parameters:}
\begin{description}[labelindent = 3em]
\item[g\textunderscore nx] an instance of a NetworkX graph that must have for each edge at least the attribute \texttt{'Length'} with value a positive scalar.
\item[Np] (optional) an integer corresponding to the total number of discretization points on the graph. By default, the number of discretization points is set to $100$ on each edge.
\item[user\textunderscore bc] (optional) a dictionary whose keys are the identifiers of vertices used to describe edges in \textbf{g\textunderscore nx} and whose values must be of the form: \texttt{['Dirichlet']} for Dirichlet boundary condition, \texttt{['Kirchhoff']} for Kirchhoff-Neumann boundary condition, \texttt{['Delta',val]} for a $\delta$ boundary condition  with a strength equal to \texttt{val} (which must be a scalar), \texttt{['Delta\textunderscore Prime',val]} for a $\delta'$ boundary condition with a strength equal to \texttt{val} (which must be a scalar) or \texttt{['UserDefined',[A\textunderscore v,B\textunderscore v]]} for a user-defined boundary condition with matrices $A_v$ and $B_v$ which must be $2$-dimensional \texttt{numpy.array} instances. For a full description of all boundary conditions (see~\cite{BeDuLe20} or~\cite{BeKu13}). By default, the boundary conditions for all vertices are Kirchhoff-Neumann boundary conditions.
\end{description}

\vspace{0.5em}
\textit{Return:}
\begin{description}[labelindent = 3em]
\item[g] an instance of the Graph class which contains the finite-differences discretization of the Laplace operator as well as the identity matrix corresponding to the identity operator on the graph. We can access these matrices with \texttt{g.Lap} and \texttt{g.Id} which are $2$-dimensional \texttt{scipy.sparse} instances. The sparse format of these instances is \texttt{csc} (Compressed Sparse Columns).
\end{description}

\subsubsection{Position} A method that enables the user to set the position (on the $x,y$-plane) of every vertex on the graph. This is only useful when drawing a graph or a wave-function on the graph.

\vspace{0.5em}
\noindent\textbf{Position(g,dict\textunderscore nodes)}

\vspace{0.5em}
\textit{Parameters:}
\begin{description}[labelindent = 3em]
\item[g] an instance of the Graph class whose edges' position will be set.
\item[dict\textunderscore nodes] a dictionary whose keys are the identifiers of the nodes used in \textbf{g} and whose values must be of the form \texttt{[posx,posy]} where \texttt{posx} and \texttt{posy} must be scalars corresponding to the desired $x$ and $y$ coordinates associated to the key node.
\end{description}

\subsubsection{draw} A method to plot the graph in the Matplotlib figure named
\verb+'QGraph'+. Each vertex is represented as a dot and its associated label is
displayed.

\vspace{0.5em}
\noindent\textbf{draw(g,AxId,Color,Text,TextSize,LineWidth,MarkerSize,FigName)}

\vspace{0.5em}
\textit{Parameters:}
\begin{description}[labelindent = 3em]
\item[g] an instance of the Graph class.
\item[AxId] (optional) an \verb+Axes+ instance of the Matplotlib library. Allow
  to draw the graph \verb+g+ in an already existing axes.
\item[Color] (optional) by default, the color of the graph is blue. It allows to
  specify an alternative color. The user must follow the standard naming color
  of Matplotlib library.
\item[Text] (optional) Logical variable. This option allows to control the
  display of the vertices labels. By default, \verb+Text=True+. To avoid the display of labels, set \verb+Text=False+.
\item[TextSize] (optional) a float variable. This allows to control the text
  size to display vertices labels. By default, the text size parameter is set to 12.
\item[LineWidth] (optional) a float variable. This allows to control the width
  of the curve representing an edge. By default, the value is set to 1.
\item[MarkerSize] (optional) a float variable. This allows to control the size
  of the marker representing the vertices of the graph. The default value is 20.
\item[FigName] (optional) a string variable. By default, the name of the figure
  is \verb+'QGraph'+. The user can change the name of the figure.
\end{description}

\vspace{0.5em}
\textit{Return:}
\begin{description}[labelindent = 3em]
\item[fig] the figure Matplotlib instance containing the axes \verb+ax+.
\item[ax] the axes Matplotlib instance containing the plot of \verb+g+.
\end{description}

\subsubsection{Diag} A method constructing a diagonal matrix with respect to the discretization points on the graph. The diagonal is explicitly prescribed.

\vspace{0.5em}
\noindent\textbf{M = Diag(g,diag\textunderscore vect)}

\vspace{0.5em}
\textit{Parameters:}
\begin{description}[labelindent = 3em]
\item[g] an instance of the Graph class.
\item[diag\textunderscore vect] either an instance of WFGraph or a $1$-dimensional \texttt{numpy.array} corresponding to the desired diagonal.
\end{description}

\vspace{0.5em}
\textit{Return:}
\begin{description}[labelindent = 3em]
\item[M] a matrix whose diagonal corresponds to \texttt{diag\textunderscore vect}. It is a $2$-dimensional \texttt{scipy.sparse} instance. The sparse format of this instance is \texttt{csc} (Compressed Sparse Columns).
\end{description}

\subsection{Methods from the WFGraph class}

\subsubsection{\textunderscore\textunderscore init\textunderscore\textunderscore}

The WFGraph class constructor builds a discrete function that is described on a discretized graph (given by an instance of the Graph class).

\vspace{0.5em}
\noindent\textbf{psi = WFGraph(initWF,g,Dtype)}

\vspace{0.5em}
\textit{Parameters:}
\begin{description}[labelindent = 3em]
\item[initWF] either a dictionary whose keys are the identifiers of edges of \texttt{g} and whose values are \texttt{lambda} functions with a single argument (say \texttt{x}) describing the desired function in an analytical way on the corresponding edge or a $1$-dimensional \texttt{numpy.array} instance which corresponds to the discretized function on the discretization points of \texttt{g}. Note that, in the first case, the variable \texttt{x} will take values between $0$ and the length of the edge (starting at the node corresponding to the first coordinate of the edge's identifier).
\item[g] (optional) an instance of \texttt{Graph} on which the function is
  described. If it has already been set in a previous instance of
  \texttt{WFGraph}, it does not need to be prescribed again. 
\item[Dtype] (optional) a string set by default to
  \verb+'float'+. The default data type for \texttt{numpy.arrays} is
  \texttt{np.float64}. It is possible to switch to complex arrays by setting \texttt{Dtype = 'complex'}.
\end{description}

\vspace{0.5em}
\textit{Return:}
\begin{description}[labelindent = 3em]
\item[psi] an instance of the WFGraph class which contains \texttt{vect}, a $1$-dimensional \texttt{numpy.array} associated to the discretization points of the graph \texttt{g}.\end{description}

\subsubsection{norm} This method enables to compute the $\ell^p$-norm of a discrete function on a graph. It is computed with a trapezoidal rule on each vertex of the graph.

\vspace{0.5em}
\noindent\textbf{a = norm(psi,p)}

\vspace{0.5em}
\textit{Parameters:}
\begin{description}[labelindent = 3em]
\item[psi] an instance of the WFGraph class whose norm is computed on its associated graph.
\item[p] a scalar value that corresponds to the exponent of the $\ell^p$ space.
\end{description}

\vspace{0.5em}
\textit{Return:}
\begin{description}[labelindent = 3em]
\item[a] a scalar value that is the $\ell^p$-norm of \textbf{psi} on its graph.
\end{description}

\subsubsection{dot} A method that computes the $\ell^2$ (hermitian) inner product between two discrete functions on a graph.

\vspace{0.5em}
\noindent\textbf{a = dot(psi,phi)}

\vspace{0.5em}
\textit{Parameters:}
\begin{description}[labelindent = 3em]
\item[psi] an instance of the WFGraph class.
\item[phi] an instance of the WFGraph class.
\end{description}

\vspace{0.5em}
\textit{Return:}
\begin{description}[labelindent = 3em]
\item[a] a (complex) scalar value that corresponds to the inner product between \textbf{psi} and \textbf{phi} on their associated graph.
\end{description}

\subsubsection{draw} A method to plot an instance \verb+f+ of the WFGraph class
and the graph \verb+g+ (instance of the Graph class) in the Matplotlib figure named
\texttt{'Wave function on the graph'}. Each vertex of \verb+g+ is represented as a
dot and its associated label is displayed.

\vspace{0.5em}
\noindent\textbf{draw(f,data\textunderscore plot,fig\textunderscore
  name,Text,AxId,ColorWF,ColorG,TextSize,LineWidth,MarkerSize,...}\\ \hspace*{10mm}\textbf{LineWidthG,AlphaG,xlim,ylim)}

\vspace{0.5em}
\textit{Parameters:}
\begin{description}[labelindent = 3em]
\item[f] an instance of the WFGraph class.
\item[data\textunderscore plot] (optional) a list of elements (matplotlib
  primitives) representing the plot of \verb+f+ already existing in the
  figure. This variable allows to efficiently update the figure containing the
  plot of the wave function \verb+f+ without redrawing all the scene.
\item[fig\textunderscore name] (optional) a string variable. By default, the name of the figure
  is \texttt{'Wave function on the graph'}. The user can change the name of the figure.
\item[Text] (optional) Logical variable. This option allows to control the
  display of the labels of the vertices of \verb+g+. By default, \verb+Text=True+. To avoid the display of labels, set \verb+Text=False+.
\item[AxId] (optional) an \verb+Axes+ instance of the Matplotlib library. Allow
  to draw \verb+f+ and \verb+g+ in an already existing axes.
\item[ColorWF](optional) by default, the color of \texttt{f} is blue. It allows to
  specify an alternative color. The user must follow the standard naming color
  of Matplotlib library.
\item[ColorG] (optional) by default, the color of \texttt{g} is dark gray. It allows to
  specify an alternative color. The user must follow the standard naming color
  of Matplotlib library.
\item[TextSize] (optional) a float variable. This allows to control the text
  size to display labels of vertices of \texttt{g}. By default, the text size parameter is set to 10.
\item[LineWidth] (optional) a float variable. This allows to control the width
  of the curve representing \texttt{f}. By default, the value is set to $1.5$.
\item[MarkerSize] (optional) a float variable. This allows to control the size
  of the marker representing the vertices of \texttt{g}. The default value is 10.
\item[LineWidthG] (optional) a float variable. This allows to control the width
  of the curves representing the edges of \texttt{g}. By default, the value is set to $0.8$.
\item[AlphaG] (optional) a float variable belonging to $[0,1]$. It allows to
  adjust the transparency (alpha property) of the graph \texttt{g} (both the
  edges, markers and labels). By defaults,
  the value is set to $1$. If the user chooses \texttt{AlphaG=0}, the graph
  \texttt{g} is completely transparent and does not appear.
\item[xlim] (optional) a two-components list instance allowing to adjust the x-axis view limits.
\item[ylim] (optional) a two-components list instance allowing to adjust the y-axis view limits.
\end{description}

\vspace{0.5em}
\textit{Return:}
\begin{description}[labelindent = 3em]
\item[K] a list of elements (matplotlib
  primitives) representing the plot of \verb+f+ in \verb+ax+.
\item[fig] the figure Matplotlib instance the axes \verb+ax+.
\item[ax] the axes Matplotlib instance containing the plot of \verb+f+ and \verb+g+.
\end{description}

\subsubsection{Arithmetic operations: \texttt{+}, \texttt{-}, \texttt{*} and \texttt{/}} The basic arithmetic operations can be applied to two instances of \texttt{WFGraph}. As a matter of fact, these operations are carried pointwise on the \texttt{vect} associated to each instance. The output is an instance of \texttt{WFGraph} with the resulting \texttt{vect} associated.

\subsubsection{Mathematical functions:  \texttt{abs}, \texttt{Real}, \texttt{Imag}, \texttt{**}, \texttt{exp}, \texttt{cos}, \texttt{sin} and \texttt{log}}
Some basic mathematical functions can be applied to an instance of \texttt{WFGraph}. It turns out that the function is applied pointwise on the \texttt{vect} associated to the instance. The output is an instance of \texttt{WFGraph} with the resulting \texttt{vect} associated.

\subsubsection{Lap} This method applies the (finite-differences) discretization of the Laplace operator to a discrete function on a graph.

\vspace{0.5em}
\noindent\textbf{phi = Lap(psi)}

\vspace{0.5em}
\textit{Parameters:}
\begin{description}[labelindent = 3em]
\item[psi] an instance of the WFGraph class on which the discrete Laplace operator is applied (specifically, on its associated \texttt{vect}).
\end{description}

\vspace{0.5em}
\textit{Return:}
\begin{description}[labelindent = 3em]
\item[phi] an instance of the WFGraph class.
\end{description}

\subsubsection{Solve} A method that solves a linear system where the matrix is understood as a discrete operator and the right-hand-side is understood as a discrete function on a graph.

\vspace{0.5em}
\noindent\textbf{phi = Solve(M,psi)}

\vspace{0.5em}
\textit{Parameters:}
\begin{description}[labelindent = 3em]
\item[M] a $2$-dimensional \texttt{scipy.sparse} instance which is associated to a discrete operator on the graph of \textbf{psi} and that we formally want to inverse.
\item[psi] an instance of the WFGraph class which correspond to the right-hand-side of the linear system.
\end{description}

\vspace{0.5em}
\textit{Return:}
\begin{description}[labelindent = 3em]
\item[phi] an instance of the WFGraph class.
\end{description}

\bibliographystyle{abbrv}
\bibliography{master}

\def\cprime{$'$}
\begin{thebibliography}{10}

\bibitem{adami2012stationary}
R.~Adami, C.~Cacciapuoti, D.~Finco, and D.~Noja.
\newblock {Stationary states of NLS on star graphs}.
\newblock {\em EPL (Europhysics Letters)}, 100(1):10003, 2012.

\bibitem{AdCaFiNo14a}
R.~Adami, C.~Cacciapuoti, D.~Finco, and D.~Noja.
\newblock Constrained energy minimization and orbital stability for the {NLS}
  equation on a star graph.
\newblock {\em Ann. Inst. H. Poincar\'e Anal. Non Lin\'eaire},
  31(6):1289--1310, 2014.

\bibitem{AdDo18}
R.~Adami and S.~Dovetta.
\newblock One-dimensional versions of three-dimensional system: ground states
  for the {NLS} on the spatial grid.
\newblock {\em Rend. Mat. Appl. (7)}, 39(2):181--194, 2018.

\bibitem{AdDoRu19}
R.~Adami, S.~Dovetta, and A.~Ruighi.
\newblock Quantum graphs and dimensional crossover: the honeycomb.
\newblock {\em Commun. Appl. Ind. Math.}, 10(1):109--122, 2019.

\bibitem{AdDoSeTi19}
R.~Adami, S.~Dovetta, E.~Serra, and P.~Tilli.
\newblock Dimensional crossover with a continuum of critical exponents for
  {NLS} on doubly periodic metric graphs.
\newblock {\em Anal. PDE}, 12(6):1597--1612, 2019.

\bibitem{AdSeTi15b}
R.~Adami, E.~Serra, and P.~Tilli.
\newblock N{LS} ground states on graphs.
\newblock {\em Calc. Var. Partial Differential Equations}, 54(1):743--761,
  2015.

\bibitem{AdSeTi16}
R.~Adami, E.~Serra, and P.~Tilli.
\newblock Threshold phenomena and existence results for {NLS} ground states on
  metric graphs.
\newblock {\em J. Funct. Anal.}, 271(1):201--223, 2016.

\bibitem{AdSeTi17b}
R.~Adami, E.~Serra, and P.~Tilli.
\newblock Nonlinear dynamics on branched structures and networks.
\newblock {\em Riv. Math. Univ. Parma (N.S.)}, 8(1):109--159, 2017.

\bibitem{AdSeTi19}
R.~Adami, E.~Serra, and P.~Tilli.
\newblock Multiple positive bound states for the subcritical {NLS} equation on
  metric graphs.
\newblock {\em Calc. Var. Partial Differential Equations}, 58(1):Art. 5, 16,
  2019.

\bibitem{Al94}
F.~Ali~Mehmeti.
\newblock {\em Nonlinear waves in networks}, volume~80 of {\em Mathematical
  Research}.
\newblock Akademie-Verlag, Berlin, 1994.

\bibitem{AlVoNi01}
F.~Ali~Mehmeti, J.~von Below, and S.~Nicaise, editors.
\newblock {\em Partial differential equations on multistructures}, volume 219
  of {\em Lecture Notes in Pure and Applied Mathematics}. Marcel Dekker, Inc.,
  New York, 2001.

\bibitem{antoine2017efficient}
X.~Antoine, A.~Levitt, and Q.~Tang.
\newblock {Efficient spectral computation of the stationary states of rotating
  Bose--Einstein condensates by preconditioned nonlinear conjugate gradient
  methods}.
\newblock {\em Journal of Computational Physics}, 343:92--109, 2017.

\bibitem{BaIg14}
V.~Banica and L.~I. Ignat.
\newblock Dispersion for the {S}chr\"odinger equation on the line with multiple
  {D}irac delta potentials and on delta trees.
\newblock {\em Anal. PDE}, 7(4):903--927, 2014.

\bibitem{BeKu13}
G.~Berkolaiko and P.~Kuchment.
\newblock {\em Introduction to quantum graphs}, volume 186 of {\em Mathematical
  Surveys and Monographs}.
\newblock American Mathematical Society, Providence, RI, 2013.

\bibitem{BeMaPe19}
G.~Berkolaiko, J.~L. Marzuola, and D.~E. Pelinovsky.
\newblock Edge-localized states on quantum graphs in the limit of large mass,
  2019.

\bibitem{besse2004relaxation}
C.~Besse.
\newblock {A relaxation scheme for the nonlinear Schr{\"o}dinger equation}.
\newblock {\em SIAM Journal on Numerical Analysis}, 42(3):934--952, 2004.

\bibitem{BeDuLe20}
C.~Besse, R.~Duboscq, and S.~Le~Coz.
\newblock {Gradient Flow Approach to the Calculation of Ground States on
  Nonlinear Quantum Graphs}.
\newblock arXiv:2006.04404, June 2020.

\bibitem{Grafidi}
C.~Besse, R.~Duboscq, and S.~Le~Coz.
\newblock Grafidi.
\newblock {\em PLMlab repository},
  \url{https://plmlab.math.cnrs.fr/cbesse/grafidi}, 2021.

\bibitem{BhBoHe20}
K.~Bhandari, F.~Boyer, and V.~Hern{\'a}ndez-Santamar{\'i}a.
\newblock {Boundary null-controllability of 1-D coupled parabolic systems with
  Kirchhoff-type condition}.
\newblock working paper or preprint, hal-02748405, July 2020.

\bibitem{CaDoSe18}
C.~Cacciapuoti, S.~Dovetta, and E.~Serra.
\newblock Variational and stability properties of constant solutions to the
  {NLS} equation on compact metric graphs.
\newblock {\em Milan J. Math.}, 86(2):305--327, 2018.

\bibitem{CaFiNo15}
C.~Cacciapuoti, D.~Finco, and D.~Noja.
\newblock Topology-induced bifurcations for the nonlinear {S}chr\"{o}dinger
  equation on the tadpole graph.
\newblock {\em Phys. Rev. E (3)}, 91(1):013206, 8, 2015.

\bibitem{danaila2017computation}
I.~Danaila and B.~Protas.
\newblock {Computation of ground states of the Gross--Pitaevskii functional via
  Riemannian optimization}.
\newblock {\em SIAM Journal on Scientific Computing}, 39(6):B1102--B1129, 2017.

\bibitem{delfour1981finite}
M.~Delfour, M.~Fortin, and G.~Payr.
\newblock Finite-difference solutions of a non-linear schr{\"o}dinger equation.
\newblock {\em Journal of computational physics}, 44(2):277--288, 1981.

\bibitem{Do18}
S.~Dovetta.
\newblock Existence of infinitely many stationary solutions of the
  {$L^2$}-subcritical and critical {NLSE} on compact metric graphs.
\newblock {\em J. Differential Equations}, 264(7):4806--4821, 2018.

\bibitem{Do19}
S.~Dovetta.
\newblock Mass-constrained ground states of the stationary {NLSE} on periodic
  metric graphs.
\newblock {\em NoDEA Nonlinear Differential Equations Appl.}, 26(5):Paper No.
  30, 30, 2019.

\bibitem{DoGhMiPi20}
S.~Dovetta, M.~Ghimenti, A.~M. Micheletti, and A.~Pistoia.
\newblock Peaked and low action solutions of nls equations on graphs with
  terminal edges.
\newblock {\em SIAM Journal on Mathematical Analysis}, 52(3):2874--2894, 2020.

\bibitem{DoSeTi20}
S.~Dovetta, E.~Serra, and P.~Tilli.
\newblock Nls ground states on metric trees: existence results and open
  questions.
\newblock {\em Journal of the London Mathematical Society}, n/a(n/a), 2020.

\bibitem{GiPeSc16}
S.~Gilg, D.~Pelinovsky, and G.~Schneider.
\newblock Validity of the {NLS} approximation for periodic quantum graphs.
\newblock {\em NoDEA Nonlinear Differential Equations Appl.}, 23(6):Art. 63,
  30, 2016.

\bibitem{GnWa16}
S.~Gnutzmann and D.~Waltner.
\newblock Stationary waves on nonlinear quantum graphs: general framework and
  canonical perturbation theory.
\newblock {\em Phys. Rev. E}, 93(3):032204, 19, 2016.

\bibitem{Go19}
R.~H. Goodman.
\newblock N{LS} bifurcations on the bowtie combinatorial graph and the dumbbell
  metric graph.
\newblock {\em Discrete Contin. Dyn. Syst.}, 39(4):2203--2232, 2019.

\bibitem{GuLeTs17}
S.~Gustafson, S.~Le~Coz, and T.-P. Tsai.
\newblock Stability of periodic waves of 1{D} cubic nonlinear {S}chr\"odinger
  equations.
\newblock {\em Appl. Math. Res. Express. AMRX}, 2:431--487, 2017.

\bibitem{hagberg-2008-exploring}
A.~A. Hagberg, D.~A. Schult, and P.~J. Swart.
\newblock Exploring network structure, dynamics, and function using {NetworkX}.
\newblock In {\em Proceedings of the 7th Python in Science Conference
  (SciPy2008)}, pages 11--15, Pasadena, CA USA, Aug. 2008.

\bibitem{harris2020array}
C.~R. Harris, K.~J. Millman, S.~J. van~der Walt, R.~Gommers, P.~Virtanen,
  D.~Cournapeau, E.~Wieser, J.~Taylor, S.~Berg, N.~J. Smith, R.~Kern, M.~Picus,
  S.~Hoyer, M.~H. van Kerkwijk, M.~Brett, A.~Haldane, J.~F. del R{'{\i}}o,
  M.~Wiebe, P.~Peterson, P.~G{'{e}}rard-Marchant, K.~Sheppard, T.~Reddy,
  W.~Weckesser, H.~Abbasi, C.~Gohlke, and T.~E. Oliphant.
\newblock Array programming with {NumPy}.
\newblock {\em Nature}, 585(7825):357--362, Sept. 2020.

\bibitem{HuTrMa11}
N.~V. Hung, M.~Trippenbach, and B.~A. Malomed.
\newblock {Symmetric and asymmetric solitons trapped in $\mathsf{H}$-shaped
  potentials}.
\newblock {\em Phys. Rev. A}, 84:053618, Nov 2011.

\bibitem{Hunter:2007}
J.~D. Hunter.
\newblock {Matplotlib: A 2D graphics environment}.
\newblock {\em Computing in Science \& Engineering}, 9(3):90--95, 2007.

\bibitem{IaLeRo17}
I.~Ianni, S.~Le~Coz, and J.~Royer.
\newblock On the {C}auchy problem and the black solitons of a singularly
  perturbed {G}ross-{P}itaevskii equation.
\newblock {\em SIAM J. Math. Anal.}, 49(2):1060--1099, 2017.

\bibitem{KaMaPeXi20}
A.~Kairzhan, R.~Marangell, D.~E. Pelinovsky, and K.~L. Xiao.
\newblock Standing waves on a flower graph, 2020.

\bibitem{KuSh20}
K.~Kurata and M.~Shibata.
\newblock Least energy solutions to semi-linear elliptic problems on metric
  graphs.
\newblock {\em Journal of Mathematical Analysis and Applications},
  491(1):124297, 2020.

\bibitem{LeFuFiKsSi08}
S.~Le~Coz, R.~Fukuizumi, G.~Fibich, B.~Ksherim, and Y.~Sivan.
\newblock Instability of bound states of a nonlinear {S}chr\"odinger equation
  with a {D}irac potential.
\newblock {\em Phys. D}, 237(8):1103--1128, 2008.

\bibitem{MaPe16}
J.~L. Marzuola and D.~E. Pelinovsky.
\newblock Ground {S}tate on the {D}umbbell {G}raph.
\newblock {\em Appl. Math. Res. Express. AMRX}, 2016(1):98--145, 2016.

\bibitem{MeAmNi15}
F.~A. Mehmeti, K.~Ammari, and S.~Nicaise.
\newblock Dispersive effects and high frequency behaviour for the
  {S}chr\"{o}dinger equation in star-shaped networks.
\newblock {\em Port. Math.}, 72(4):309--355, 2015.

\bibitem{MeAmNi17}
F.~A. Mehmeti, K.~Ammari, and S.~Nicaise.
\newblock Dispersive effects for the {S}chr\"{o}dinger equation on the tadpole
  graph.
\newblock {\em J. Math. Anal. Appl.}, 448(1):262--280, 2017.

\bibitem{No14}
D.~Noja.
\newblock Nonlinear {S}chr\"odinger equation on graphs: recent results and open
  problems.
\newblock {\em Philos. Trans. R. Soc. Lond. Ser. A Math. Phys. Eng. Sci.},
  372(2007):20130002, 20, 2014.

\bibitem{NoPeSh15}
D.~Noja, D.~Pelinovsky, and G.~Shaikhova.
\newblock Bifurcations and stability of standing waves in the nonlinear
  {S}chr\"{o}dinger equation on the tadpole graph.
\newblock {\em Nonlinearity}, 28(7):2343--2378, 2015.

\bibitem{NoPe20}
D.~Noja and D.~E. Pelinovsky.
\newblock {Standing waves of the quintic NLS equation on the tadpole graph},
  2020.

\bibitem{Pa18}
A.~Pankov.
\newblock Nonlinear {S}chr\"{o}dinger equations on periodic metric graphs.
\newblock {\em Discrete Contin. Dyn. Syst.}, 38(2):697--714, 2018.

\bibitem{PeSc17}
D.~Pelinovsky and G.~Schneider.
\newblock Bifurcations of standing localized waves on periodic graphs.
\newblock {\em Ann. Henri Poincar\'{e}}, 18(4):1185--1211, 2017.

\bibitem{PiSoVe20}
D.~Pierotti, N.~Soave, and G.~Verzini.
\newblock {Local minimizers in absence of ground states for the critical NLS
  energy on metric graphs}.
\newblock {\em Proceedings of the Royal Society of Edinburgh: Section A
  Mathematics}, page 1–29, 2020.

\bibitem{SaSoBaMa13}
K.~K. Sabirov, Z.~A. Sobirov, D.~Babajanov, and D.~U. Matrasulov.
\newblock Stationary nonlinear {S}chr\"{o}dinger equation on simplest graphs.
\newblock {\em Phys. Lett. A}, 377(12):860--865, 2013.

\bibitem{SoBaMa17}
Z.~Sobirov, D.~Babajanov, and D.~Matrasulov.
\newblock Nonlinear standing waves on planar branched systems: shrinking into
  metric graph.
\newblock {\em Nanosystems: Physics, Chemistry, Mathematics}, 8(1):29, 2017.

\bibitem{strang1968construction}
G.~Strang.
\newblock On the construction and comparison of difference schemes.
\newblock {\em SIAM journal on numerical analysis}, 5(3):506--517, 1968.

\bibitem{ToOsDe08}
A.~Tokuno, M.~Oshikawa, and E.~Demler.
\newblock {Dynamics of One-Dimensional Bose Liquids: Andreev-Like Reflection at
  $Y$ Junctions and the Absence of the Aharonov-Bohm Effect}.
\newblock {\em Phys. Rev. Lett.}, 100:140402, Apr 2008.

\bibitem{2020SciPy-NMeth}
P.~Virtanen, R.~Gommers, T.~E. Oliphant, M.~Haberland, T.~Reddy, D.~Cournapeau,
  E.~Burovski, P.~Peterson, W.~Weckesser, J.~Bright, S.~J. {van der Walt},
  M.~Brett, J.~Wilson, K.~J. Millman, N.~Mayorov, A.~R.~J. Nelson, E.~Jones,
  R.~Kern, E.~Larson, C.~J. Carey, {\.I}.~Polat, Y.~Feng, E.~W. Moore,
  J.~{VanderPlas}, D.~Laxalde, J.~Perktold, R.~Cimrman, I.~Henriksen, E.~A.
  Quintero, C.~R. Harris, A.~M. Archibald, A.~H. Ribeiro, F.~Pedregosa, P.~{van
  Mulbregt}, and {SciPy 1.0 Contributors}.
\newblock {{SciPy} 1.0: Fundamental Algorithms for Scientific Computing in
  Python}.
\newblock {\em Nature Methods}, 17:261--272, 2020.

\bibitem{weideman1986split}
J.~Weideman and B.~Herbst.
\newblock {Split-step methods for the solution of the nonlinear Schr{\"o}dinger
  equation}.
\newblock {\em SIAM Journal on Numerical Analysis}, 23(3):485--507, 1986.

\end{thebibliography}



\end{document}